\documentclass{article}
\pdfoutput=1
\usepackage{arxiv}

\usepackage[utf8]{inputenc} 
\usepackage[T1]{fontenc}    
\usepackage{url}            
\usepackage{booktabs}       
\usepackage{amsfonts}       
\usepackage{nicefrac}       
\usepackage{microtype}      
\usepackage{lipsum}
\usepackage{xcolor}
\usepackage{subfigure}
\usepackage{placeins}
\usepackage{graphicx}
\usepackage{enumerate}

\usepackage[numbers]{natbib}     
\usepackage{color}             
\usepackage{footnote}          
\usepackage{amssymb,amsmath,amsfonts}
\usepackage{listings}
\usepackage{mathrsfs}
\usepackage[utf8]{inputenc}
\usepackage{indentfirst}

\def\avg#1{\langle #1 \rangle}

\def\avg#1{\left \langle #1 \right \rangle}

\def\fe2{FE$^2$}
\def\DeepBCn{DeepBND}
\def\DeepBC{DeepBND }
\def\DeepBCNN{DeepBND-NN }

\def\d{\mathrm{d}}
\def\bs#1{\boldsymbol{#1}}
\def\Vec#1{\mathbf{#1}}
\def\Mat#1{\boldsymbol{\mathsf{#1}}}

\def\ten#1{\mathbf{#1}}
\def\Ten4#1{\mathbb{#1}}

\def\Real{\mathbb{R}}

\DeclareMathOperator*{\Span}{span}
\renewcommand{\div}{\operatorname{div}}

\def\x{\Vec{x}}

\def\y{\Vec{y}}

\newcommand{\bmat}{\begin{bmatrix}}
\newcommand{\emat}{\end{bmatrix}}

\def\uf{\tilde{\Vec{u}}_{\mu}}
\def\vf{\tilde{\Vec{v}}_{\mu}}

\def\umu{\Vec{u}_{\!\mu}}

\def\uf{\tilde{\Vec{u}}_{\!\mu}}

\def\u{\Vec{u}}

\def\H1V{\mathbf{H}^{1}}

\newtheorem{problem}{Problem} 
\newtheorem{remark}{Remark} 

\usepackage{float}

\floatstyle{boxed}

\newfloat{algorithmbox}{thp}{lop}
\floatname{algorithmbox}{Algorithm}

\def\moved{\color{black}}
\def\coupling{\color{black}}
\def\added{\color{black}}
\def\addedSimone{\color{black}}

\def\figdir{figures}

\title{DeepBND: a Machine Learning approach to enhance Multiscale Solid Mechanics} 

\author{
  Felipe Rocha \thanks{Corresponding author.}\\
  \texttt{felipe.figueredorocha@epfl.ch} \\
   \And
 Simone Deparis \\
  \And
  Pablo Antolin \\
  \And
  Annalisa Buffa \\
  \And \\
  \begin{minipage}{.8\textwidth}
  \centering
  	\textit{Ecole polytechnique f{\'e}d{\'e}rale de Lausanne ‐ SB MATH MNS MA Station 8 CH-1015 Lausanne, Switzerland}
  \end{minipage}
}

\begin{document}
\maketitle

\begin{abstract}
Effective properties of materials with random heterogeneous structures are typically determined by homogenising the mechanical quantity of interest in a window of observation. The entire problem setting encompasses the solution of a local PDE and some averaging formula for the quantity of interest in such domain. There are relatively standard methods in the literature to completely determine the formulation except for two choices: i) the local domain itself and the ii) boundary conditions. Hence, the modelling errors are governed by the quality of these two choices. The choice i) relates to the degree of representativeness of a microscale sample, i.e., it is essentially a statistical characteristic. Naturally, its reliability is higher as the size of the observation window becomes larger and/or the number of samples increases. On the other hand, excepting few special cases there is no automatic guideline to handle ii). Although it is known that the overall effect of boundary condition becomes less important with the size of the microscale domain, the computational cost to simulate such large problem several times might be prohibitive even for relatively small accuracy requirements. 
Here we introduce a machine learning procedure to select most suitable boundary conditions for multiscale problems, particularly those arising in solid mechanics. We propose the combination Reduced-Order Models and Deep Neural Networks in an offline phase, whilst the online phase consists in the very same homogenisation procedure plus one (cheap) evaluation of the trained model for boundary conditions. Hence, the method allows an implementation with minimal changes in existing codes and the use of relatively small domains without losing accuracy, which reduces the computational cost by several orders of magnitude.
\end{abstract}



\section{Introduction}

Multiscale methods can be regarded as theories that link the macroscopic behaviour of continua to phenomena occurring at smaller spatial scales, a situation found in a large range of applications such as porous media or composites, to name a few. Early developments in solid mechanics date back at least from the mid-twentieth century, with the estimation of macroscopic properties of heterogeneous materials, using the self-consistent approach \citep{Hashin1963,Hill1965,Mandel1972}, and methods based on the asymptotic analysis of partial differential equations with periodic coefficients \citep{Bensoussan1978,SanchezPalencia1983}. Despite methodological differences, common across the range of different approaches is the fact that macroscopic continuum quantities, so-called homogenised, are invariably linked to their microscale counterpart fields upon a solution of an auxiliary local problem defined at microscale domain followed by some kind of averaging process.

As for of the local problem, also known as microscale or corrector problem, besides very specific situations, the choice of the boundary conditions (BCs) cannot be exactly determined. Depending on the BCs assumed, homogenised quantities can be under or overestimated, {\added yielding to the so-called Reuss' and Voigt's estimates, the lower and upper bounds, respectively \citep{Steinmann2016}. In problems involving strain localisation in solids, such choice can be even more critical, affecting the directions of nucleation bands \citep{Rocha2021}}. The most commonly known approach to mitigate the sensibility with the BCs choice is by artificially assuming periodic BCs (even for non-periodic cells) and/or by increasing the microscopic domain size \citep{Yue2007}. However, such an approach also entails larger computational costs in the context of numerical approximations. Recently, some improvement attempts have been proposed in this regard, as for example by adding specific parabolic and source terms to the original correct problem \citep{Abdulle2019}, and by a generalisation of periodic BCs of stochastic media \citep{Steinmann2019}. 

In situations of limited physical guidance, as such case of choosing the most suitable BC, data-driven and machine learning (ML) approaches are known to be useful in providing insights. Among several ML techniques, Deep Neural Networks (DNN) is known to be a universal approximator for several activation functions \citep{Guhring2019,Pinkus1999}, also established in the context of parametrised PDEs \citep{Kutyniok2019}. {\added In this contribution, we present an attempt to design a ML strategy to accurately predict BCs for the local problems in multiscale modelling. In particular our approach relies on coupling the predictor ability of a deep neural network {\addedSimone to Reduced-Order Modelling (ROM) using the Reduced Basis Method (RB-ROM) approach \citep{Quarteroni2016,Hesthaven2015}}. Indeed RB-ROM is used to reduce the output size of the DNN and to guide the choice of the loss function. The combination of RB-ROM and DNN has already been proved effective in the context of physically-based ML techniques for PDEs \citep{Hesthaven2018,Regazzoni2019,dalSanto2020}.
}

%

{\added 
	The manuscript is organised as follows: in Section \ref{sec:preliminaries} we introduce necessary concepts for the sake of completeness purposes, in particular the basic background on multiscale solid mechanics and reduced-basis for parametrised PDEs are provided in Section \ref{sec:multiscaleSolidMechanics} and \ref{sec:reducedBasisPDE}, respectively; Section \ref{sec:formulationMethod} deals with the formulation of the proposed method, where special attention is devoted to cast the problem of Multiscale Solid Mechanics into the Parametrised PDE format  in Section \ref{sec:multiscaleRB} that, among other steps, presents the key idea of the method of finding the reduced-order description of boundary (Section \ref{sec:fluctuationSpace_parametrisation}). The section ends with an useful overview of the method is Section  \ref{sec:overviewMethod}. In
	Section \ref{sec:designingNeuralNet} we propose a novel physically inspired DNN architecture {\addedSimone tailored to exploit particular physical symmetries in the reduced-basis generation (Section \ref{sec:RBfindingStrategy}).}} Finally, Section \ref{sec:dataset} and Section \ref{sec:NumericalExamples} are dedicated to the numerical implementation of the proposed methodology, concerning dataset construction, DNN model training, and relevant numerical examples in coupled micro-macro simulations, respectively.

\section{Preliminaries} \label{sec:preliminaries}

As already commented, one of the main goals of this section is to set the general context of the computational homogenisation for random media, here specialised for multiscale solid mechanics (Section \ref{sec:multiscaleSolidMechanics}), which is instrumental to the presentation of our method (Section \ref{sec:formulationMethod}). Moreover, the second goal is to review the reduced-basis setting for parametrised PDEs (Section \ref{sec:reducedBasisPDE}), which can skipped by an expert reader.

\subsection{Multiscale Solid Mechanics}\label{sec:multiscaleSolidMechanics} 

In this section, we review the basic equations of computational homogenisation in multiscale solid mechanics. The literature in this realm is vast and dates back to the '70s with the postulation of the cornerstone Hill-Mandel Macro-homogeneity Principle \citep{Hill1965,Mandel1972}. More recently, emerging from the applied mathematical community, the so-called Heterogeneous Multiscale Method (HMM) \citep{Abdulle2012} has brought significant new insights into the field. Although these two methods have independent points of departure at the first sight, it has been shown that the resulting models are coincident \citep{Fischer2018,Fischer2019}. Hence our method can be applied to both methodologies.  Our presentation follows closely \citep{SouzaNeto2015,Blanco2016a}, that provides a natural mathematical framework for the formulation of our method. Additional details that are not necessary to follow the manuscript are included in Appendix \ref{sec:PMVP}, for the sake of completeness.

Before digging into the formulation itself, let us consider the following conventions. As usual, the unit vectors $\{\Vec{e}_i\}_{i=1}^d$ refer to the canonical basis for $\Real^d$. We also use the font convention $\ten{a} \in \Real^{d}$, $\ten{A}\in \Real^{d,d}$, $\mathbb{A} \in \Real^{d,d,d,d}$ to discern between vectors, second- and fourth-order tensors, respectively. Regular type face fonts are always used to denote components (or scalars) in cartesian coordinates, e.g. $\Vec{a} = a_i \Vec{e}_i$, $\Vec{A} = A_{ij} \Vec{e}_i \otimes \Vec{e}_j$ and so-forth (Einstein's summation convention applied). Space $\Real^{\d,\d}_{sym}$ denotes symmetric real-valued second-order tensors. Finally, we also consider the following operations:
\begin{subequations}
	\begin{align}
	\Vec{a} \cdot \Vec{b} &:= a_i b_i, \\
	\ten{A} \cdot \ten{B} &:= A_{ij} B_{ij}, \\
	\mathbb{A} \ten{B} &:= A_{ijkl}  B_{kl} \Vec{e}_i \otimes \Vec{e}_j, \\
	\Vec{a} \otimes^s \Vec{b} &:= \frac{1}{2} ( \Vec{a} \otimes \Vec{b} + \Vec{b} \otimes \Vec{a} ). 
	\end{align}
\end{subequations}  

In the following, we present the multiscale solid mechanics computational homogenisation framework, starting with the macroscale followed by the microscale model.

\subsubsection{{\added Setting of the homogenisation problem}} \label{sec:generalSettingHomogenisation}

\begin{figure}[!h]
	\centering
	\includegraphics[width=1.0\linewidth]{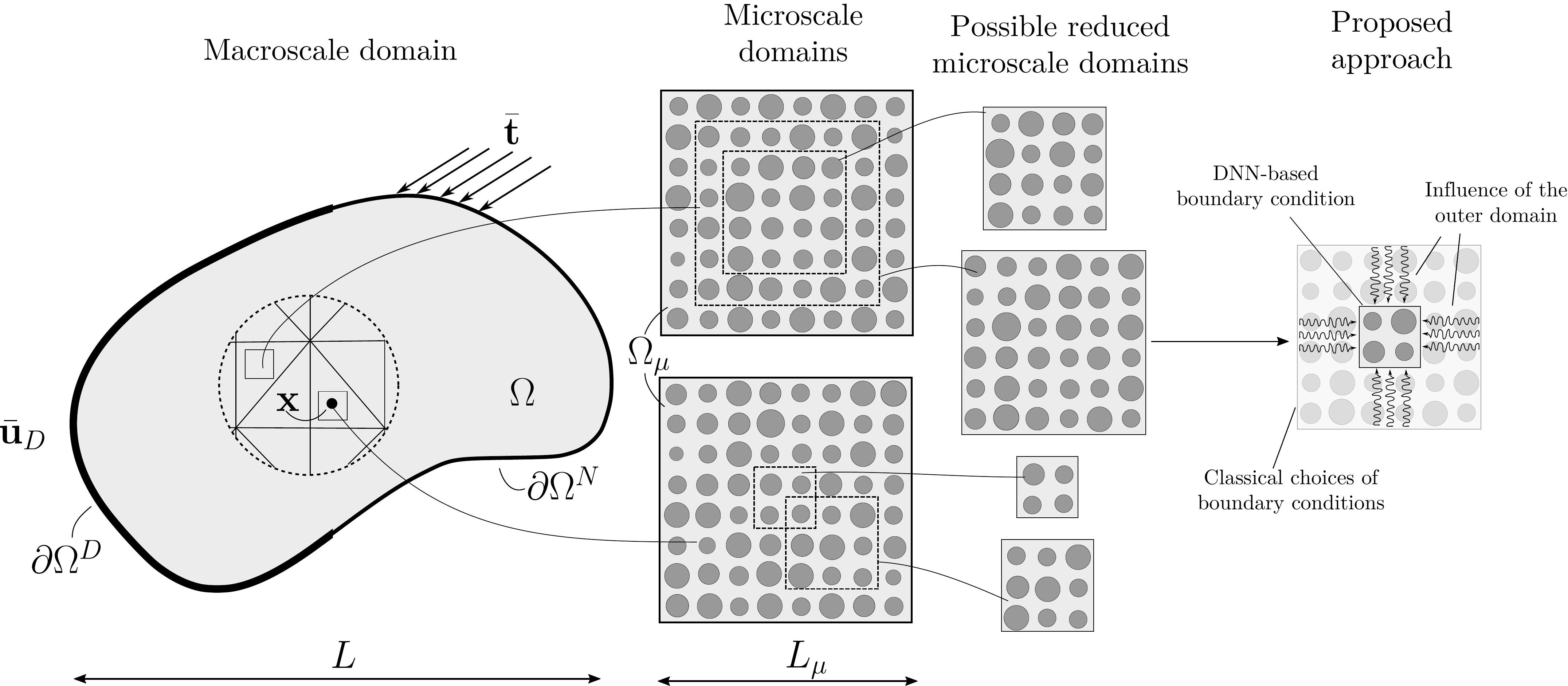}
	\caption{Homogenisation in random media and motivation of the proposed approach.}
	\label{fig:randommedia}
\end{figure}

Let us consider $\Omega \subset \Real^\d$, $\d=2$ or $3$, the domain in which the physical problem of interest is defined as shown in the left-hand side of Figure \ref{fig:randommedia}. {\coupling Thereafter, we refer to it as the \textit{macroscale domain}, whose characteristic size is $L$, as opposed the to \textit{microscale domains}, displayed on the remaining columns of Figure \ref{fig:randommedia}, whose size $L_{\mu}$ is in the same other of the magnitude of the length-scale of material properties heterogeneities. For the sake of simplicity, the example depicted is a composite material with two phases, the matrix (in light grey) and circular inclusions (in dark grey), although more general structures for the heterogeneity are also allowed.

	It is important to remark that the length of heterogeneities can be arbitrary small{\addedSimone, posing} computational difficulties if one is interested in solving such problem up to that scale. Notwithstanding, it is well established in the literature that in the scale separation limit, i.e. $\frac{L_{\mu}}{L} \ll 1$, the original problem with {(fast\nobreakdash-)varying} parameters can be replaced by the so-called homogenised problem \citep{Bensoussan1978}}. Indeed, the homogenised problem have the very same mathematical nature as the original at hand, but with effective, so-called homogenised, coefficients. The very process of finding such coefficients is called \textit{homogenisation}, which encompasses the solution of an auxiliary problem, known as \textit{corrector problem} or simply microscale problem, solved in the microscale domain, and a posterior averaging scheme using its solution, both established in Section \ref{sec:microscaleModel} for our problem.   

In the specific setting of this work, we consider as macroscale (or homogenised) problem a standard solid continuum mechanics model in the infinitesimal strain regime. We are interested in solving a quasi-static mechanical equilibrium problem, in which the displacement vector field $\u : \Omega \to \Real^d$ is obtained as the solution to the equilibrium problem once suitable BCs and constitutive equations are provided. Associated with the field $\u$, the symmetric gradient tensor field $\bs{\varepsilon}:=\nabla^s \u : \Omega \to \Real^{\d,\d}_{sym}$ is the strain measure of choice. We adopt the linear elasticity constitutive model, hence the Cauchy stress tensor $\bs{\sigma} \in \Real^{d,d}_{sym}$ is given through the Hookean law  $\bs{\sigma} = \mathbb{C} \bs{\varepsilon}$, where $\mathbb{C} \in L^2(\Omega)^{d,d,d,d}$ is a fourth-order positive-definite elastic tensor, with symmetries $C_{ijkl} = C_{jikl} = C_{ijlk} = C_{klij}$. {\coupling In the context of this work, we aim at obtaining $\mathbb{C}$ via the homogenisation procedure.  Using the notation of the left part of Figure \ref{fig:randommedia}, the macroscale problem is read as below}\footnote{For the sake of simplicity{\coupling, but without loss of generality, we are disregarding} body forces (see \citep{SouzaNeto2015}).}
{\coupling
	\begin{align} \label{prob:macroscaleProblemStrong}
	\begin{cases}
	\div \bs{\sigma} = \Vec{0} &\quad \text{in } \Omega \\
	\bs{\sigma} = \mathbb{C} \bs{\varepsilon} &\quad \text{in } \\
	\Vec{u} = \bar{\Vec{u}}_D &\quad \text{on } \partial \Omega^D \\
	\bs{\sigma} \Vec{n} = \bar{\Vec{t}} &\quad \text{on } \partial \Omega^N \\
	\end{cases},
	\end{align}
	where $\partial \Omega^D$ and $\partial \Omega^N$ are the Dirichlet and Neumann boundaries, respectively, forming partition of $\partial \Omega$, with the associated prescribed displacement $\bar{\Vec{u}}_D$ and traction $\bar{\Vec{t}}$ . It is worth noticing that} the computational implementation via FEM of {\coupling the above problem} follows standard procedures, apart from the evaluation of the $\bs{\sigma}$ (or $\mathbb{C}$) at the integration point level. For each integration point, there is an associated microstructure which should be coupled kinematically and energetically with the macroscale quantities at that point. {\added This defines the corrector problem, which however lacks boundary condition information that needs to be chosen properly. This choice lies behind the very idea of the proposed method, to be explored in Section \ref{sec:formulationMethod}.}

%
%

\subsubsection{{\coupling Micro-macro coupling}} \label{sec:microscaleModel}

For each macroscale {\added (integration)} point $\x \in \Omega$ we associate a microscale domain (MD) $\Omega_{\mu} \subset \Real^d $ \footnote{{\added We decided to use the acronym MD, instead of Representative Volume Element (RVE), to avoid discussions of statistical nature that go beyond the scope of this work. Notwithstanding, for most practical considerations MD and RVE can be understood as synonyms. Other possible nomenclature, preferred by some authors to avoid this confusion, is SVE (Statistical Volume Element) \citep{Ostoja2006}}.}. It is considered that $L_{\mu}/L \ll 1$, being $L$ and $L_{\mu}$ the characteristic lengths related to the macroscale and microscale, respectively. The characteristic size $L_{\mu}$ is chosen such that $\Omega_{\mu}$ should have  meaningful numbers and type of material heterogeneities but also keeping the compromise in terms of computational cost. Thereafter, for the sake of simplicity, we use the notation $\u = \u(\x)$, $\bs{\varepsilon} = \bs{\varepsilon}(\x)$ and $\bs{\sigma} =\bs{\sigma}(\x)$,  to represent point-valued vectors or tensors at a point $\x\in\Omega$, rather than vector or tensor fields in $\Omega$. On the other hand, objects associated with $\Omega_{\mu}$ are denoted by the subindex $\mu$.

We also consider the notation $\avg{\cdot}_{\Omega_{\mu}} = \frac{1}{|\Omega_{\mu}|} \int_{\Omega_{\mu}} (\cdot) \, \d \Omega_{\mu}$, 
$\avg{\cdot}_{\partial \Omega_{\mu}} = \frac{1}{|\Omega_{\mu}|} \int_{\partial \Omega_{\mu}}  (\cdot) \, \d \partial \Omega_{\mu}$ to denote microscale domain and boundary averages respectively. Also,
we consider a local coordinate system at $\Omega_{\mu}$, with points denoted $\y \in \Omega_{\mu}$ and such that $\int_{\Omega_{\mu}} \y \d \Omega_{\mu}= \Vec{0}$ (centred at the centroid). Accordingly, the symmetric gradient operator in the microscale coordinates becomes
$\bs{\varepsilon}_{\mu}(\cdot):= \nabla^s_\y (\cdot)$.

\begin{figure}
	\centering
	\includegraphics[width=0.7\linewidth]{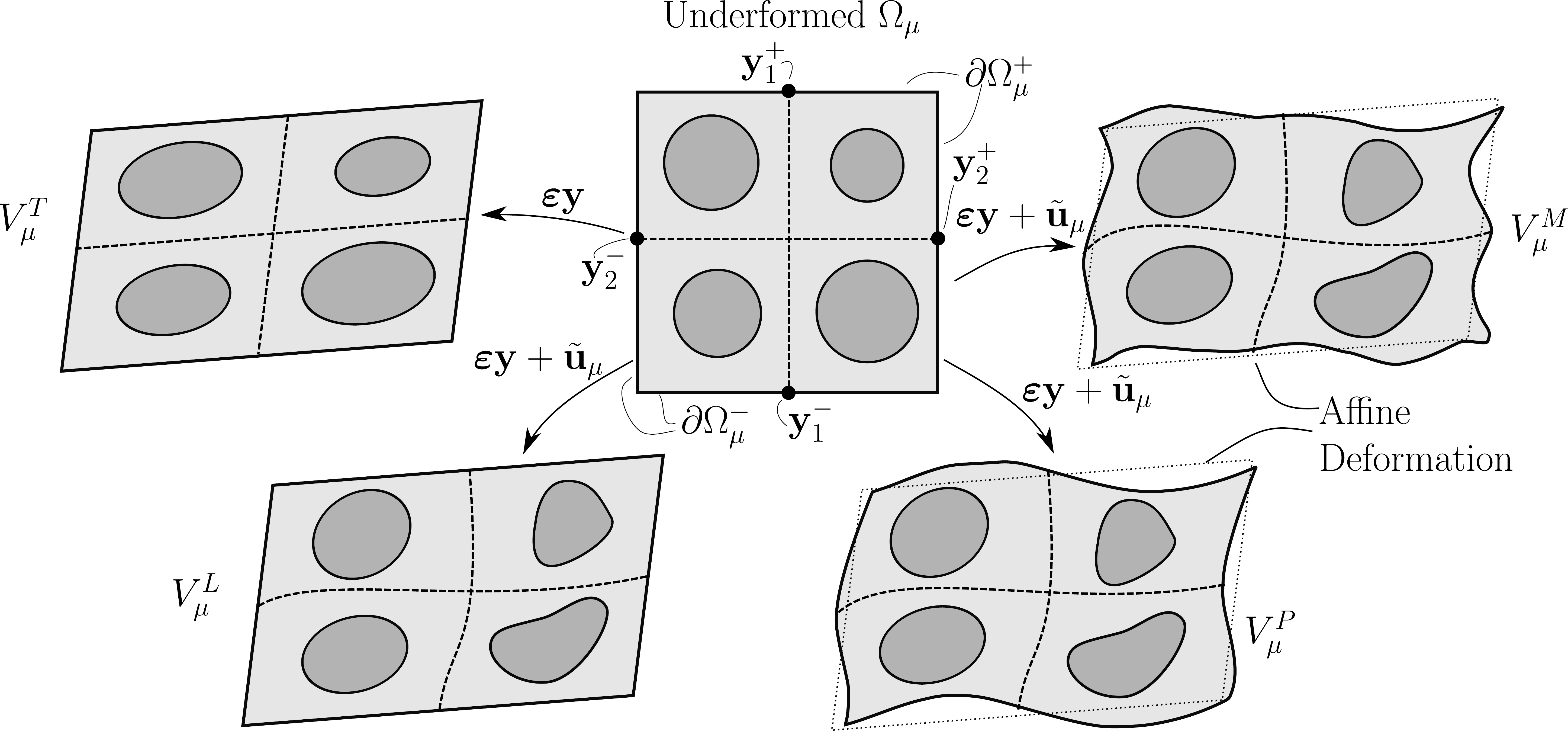}
	\caption{Comparison of classical multiscale models of kinematical constraints.}
	\label{fig:spacesboundary}
\end{figure}

The basic idea of the multiscale mechanics approach is to be able to describe the microscopic displacement $\umu :\Omega_{\mu} \to \Real^d$ which is related to the macroscale kinematics but not entirely determined by it. We assume that $\u$ and $\bs{\varepsilon}$ dictate affine deformations at the microscale while the so-called microscale displacement fluctuations field $\uf :\Omega_{\mu} \to \Real^d$ (supposedly  $\uf \in [H^1(\Omega_{\mu})]^d$) is responsible for higher-order terms. This yields to the decomposition
\begin{align} \label{eq:insertion}
\umu(\y) = \u + \bs{\varepsilon} \y + \uf(\y), \qquad \forall \y \in \Omega_{\mu}.
\end{align}
{\added Following some natural physical constraints (see Appendix \ref{sec:PMVP} for the interested reader), fluctuations $\uf$ should live in}  
\begin{align} \label{eq:VM}
V^M_{\mu} := \left \{\bs{\eta} \in [H^1(\Omega_{\mu})]^d ; \avg{\bs{\eta}}_{{\Omega}_{\mu}} = \Vec{0} ; \avg{ \bs{\eta} \otimes^s \Vec{n} }_{\partial {\Omega}_{\mu}} = \Vec{0} \right \},
\end{align}
the so-called Minimally Constrained Space for fluctuations, also known as uniform traction model. Variations of this model are possible considering subspaces of $V_{\mu}^M$. Remarkably, popular models (subspaces) in the literature are
\begin{subequations}
	\begin{align}
	V_{\mu}^T &= \{\Vec{0}\} \quad &\text{(Taylor Model)}, \\
	V_{\mu}^L &= \{\bs{\eta} \in V_{\mu}^M ; \bs{\eta}|_{\partial \Omega_{\mu}} = \Vec{0} \} \quad &\text{(Linear Boundary Model)}, \\
	V_{\mu}^P &= \{\bs{\eta} \in V_{\mu}^M ; \bs{\eta}(\y^+) = \bs{\eta}(\y^-) \; \forall \y^\pm \in \partial \Omega_{\mu}^\pm \} \quad &\text{(Periodic Model)},
	\end{align}
\end{subequations}
where the notation $(\cdot)^\pm$ denotes opposite boundary/points, e.g., left(-)/right(+) or bottom(-)/top(+) in a unit square $\Omega_{\mu}$. Figure \ref{fig:spacesboundary} depicts the kinematical effect of the different classical models, where $V_{\mu}^T$ constrain the admissible displacements to affine transformations, whilst in $V_{\mu}^L$ just the boundary deforms affinely and models $V_{\mu}^M$ and $V_{\mu}^P$ allow non-affine displacements on $\partial \Omega_{\mu}$. Excepting the space $V_{\mu}^T$, the other three classical spaces will be useful for different reasons in the proposed method. It is worth noticing that $V_{\mu}^L \subset V_{\mu}^P \subset V_{\mu}^M$, where $V_{\mu}^L$ and $V_{\mu}^M$ are known to deliver upper and lower bounds in terms of homogenised stress respectively \citep{Yue2007}. In turn, in the absence of physical guidance, $V_{\mu}^P$ is the standard choice even if the microstructure is not periodic since it delivers a more balanced response. In this work, we seek an affine subspace $V_{\mu} \subset V_{\mu}^M$ that outperforms $V_{\mu}^P$, by exploiting \textit{a priori} knowledge provided by RB-ROM and ML techniques. For convenience, thereafter we denote simply as $V_{\mu} \subset V_{\mu}^M$, the choice of one of these subspaces. The notation $V_{\mu}(\Omega_{\mu})$ may be used to avoid ambiguity in the domain of definition.  

It is also helpful for the subsequent sections to introduce an auxiliary space $V_{\mu}^{ZA}$ (from zero-average) as 
\begin{align}
V_{\mu}^{ZA} := \left \{\bs{\eta} \in [H^1(\Omega_{\mu})]^d ; \avg{\bs{\eta}}_{{\Omega}_{\mu}} = \Vec{0} \right \},
\end{align}
which differs from $V_{\mu}^M$ only by relaxing the restriction on the boundary. From the definition $V_{\mu}^{M} \subset V_{\mu}^{ZA}$.

In terms of material constitutive behaviour at microscale, we also consider a Hookean material characterised by the fourth-order micro-elasticity tensor $\mathbb{C}_{\mu} \in L^2(\Omega_{\mu})^{d,d,d,d}$ positive-definite and with the symmetries $(C_{\mu})_{ijkl} = (C_{\mu})_{jikl} = (C_{\mu})_{ijlk} = (C_{\mu})_{klij}$, which defines the microscopic stress state as $\bs{\sigma}_{\mu} = \mathbb{C}_{\mu} \bs{\varepsilon}_{\mu}(\umu)$.

{\added As result of balance of virtual powers in macro- and microscale (see details in Appendix \ref{sec:PMVP}), the corrector (or microscale equilibrium) problem and the homogenisation read as follows.} 

\begin{problem}[Microscale] \label{prob:correctorProblemSolidMechanics}
	Given $\bs{\varepsilon} \in \Real^{d,d}_{sym}$, solve:
	\begin{enumerate}[C(I)]	
		\item Corrector (or microscale equilibrium) : find $\uf \in V_{\mu}$ such that
		\begin{align} \label{eq:micro_problem}
		\avg{ \mathbb{C}_{\mu} \bs{\varepsilon}_{\mu}(\uf)  \cdot \bs{\varepsilon}_{\mu}(\tilde{\Vec{v}}_{\mu}) }_{\Omega_{\mu}} = - \avg{ \mathbb{C}_{\mu} \bs{\varepsilon}  \cdot \bs{\varepsilon}_{\mu}(\tilde{\Vec{v}}_{\mu}) }_{\Omega_{\mu}}
		\qquad \forall \tilde{\Vec{v}}_{\mu} \in V_{\mu}.
		\end{align}
		\item Stress homogenisation: let $\uf \in V_{\mu}$ be the solution of  \eqref{eq:micro_problem}, then
		\begin{align} \label{eq:homogenisation}
		\bs{\sigma} = \avg{\mathbb{C}_{\mu}(\bs{\varepsilon} + \bs{\varepsilon}_{\mu}(\uf))}_{\Omega_{\mu}}.
		\end{align}
	\end{enumerate}
\end{problem} 

Instead of computing the stress homogenisation through \eqref{eq:homogenisation}, it is convenient to obtain the homogenised elasticity tensor $\mathbb{C}$ and then evaluate $\bs{\sigma} = \mathbb{C} \bs{\varepsilon}$. This can be easily performed by solving Problem \ref{prob:correctorProblemSolidMechanics} replacing $\bs{\varepsilon}$ by the unitary strains $\ten{E}_{kl} = \Vec{e}_k \otimes^s \Vec{e}_l$, with $k,l=1,\dots,d$, and then averaging as follows
\begin{align}
\mathbb{C} = \avg{ \mathbb{C}_{\mu} ( \Vec{E}_{kl} + \bs{\varepsilon}_{\mu}(\uf^{\ten{E}_{kl}}))   }_{\Omega_{\mu}} \otimes \ten{E}_{kl}, 
\end{align}
where  and $\uf^{\ten{E}_{kl}}$ denotes the solution of \eqref{eq:micro_problem} for a given $\bs{\varepsilon} = \Vec{E}_{kl}$.

\begin{remark} 
	As already commented, the theory developed so far is consistent for any choice of $V_{\mu} \subset V_{\mu}^M$. Moreover, variationally, we are also allowed to choose different subsets of $V_{\mu}^M$ to work as trial set or test spaces. {\addedSimone This property is explored in Section \ref{sec:multiscaleRB} to reduce the MD size of problem}.
\end{remark}

{\moved 
	\begin{remark} \label{rmk:periodicity}
		{\added As already commented, the choice of BCs is not readily defined in the case of random media. The only case in which this choice is straightforward if for periodic cells, i.e. periodic BCs.}
		Particularly, {\added the material is periodic if we can identify a microscale domain which is repeated throughout (unique up to a translation)}.  If microstructures are not periodic, the corrector problem, although valid and useful, is inexact for several reasons as discussed next.
		
	\end{remark} 
}

{\added 
	\begin{remark} \label{rmk:translation} It is worth noticing that Problem \ref{prob:correctorProblemSolidMechanics} does not depend on the macroscale displacement $\u$ {\added itself, but rather on its gradient, already taken into account with $\bs{\varepsilon}$}. As consequence, in terms of macroscale kinematical dependence, $\bs{\sigma}$ is invariant to $\u$ and only depends on $\bs{\varepsilon}$. In other words, the zero-average constraint in \eqref{eq:VM} has merely the role of eliminating rigid-body movements and can be replaced by any other equivalent strategy. We use this property in Section \ref{sec:multiscaleRB}. The dependence with $\u$ is important in dynamics problems and for non-uniform body forces at the macroscale \citep{SouzaNeto2015}. 
	\end{remark}
}

\begin{remark} 
	Instead of adopting displacement fluctuations $\uf$ as primal variables, some authors prefer adopting the full microscale displacement $\umu$ \citep{Steinmann2016}, which is an equivalent framework. Nevertheless, the formulation in terms of $\uf$ is convenient since the kinematical constraints become homogeneous.
\end{remark}

\FloatBarrier

%
%

\subsection{Reduced Basis Method for Parametrised PDEs} \label{sec:reducedBasisPDE} 

In this section, we aim at providing {\addedSimone fundamental concepts of the RB-ROM} technique for parametrised PDEs. The presentation focuses on a generic problem to set a suitable notation before linking it to the specific PDE concerned, which is the focus of Section \ref{sec:multiscaleRB}. For a deeper understanding of the method, we refer to \citep{Hesthaven2015,Quarteroni2016}.

Consider a vector space $\mathcal{W}$, where the solutions of a parametrised PDE live (in our case a multiscale solid mechanics boundary value problem described in Section \ref{sec:multiscaleSolidMechanics}), equipped with the inner product $(\cdot,\cdot)_{\mathcal{W}}$ with the induced norm $\| \cdot \|_{\mathcal{W}}$.  Let $\mathcal{P} \subset \Real^{N_p}$ be a manifold of admissible parameters and let $\mathcal{L} : \mathcal{P} \to \mathcal{W}$ be the mapping from parameters onto the solutions of the  parametrised PDE at hand. Let $\mathbb{P} = \{\Vec{p}^{(i)}; i=1,\dots,N_s \in \mathcal{P}\}$ be a set containing $N_s$ samples of $\mathcal{P}$, so-called \textit{sampling} set, and the associated images upon $\mathcal{L}$, so-called \textit{snapshots}, collected in the set $\mathbb{S} := \mathcal{L}(\mathbb{P}) =  \{ \Vec{w}(\Vec{p}^{(i)})  = \mathcal{L}(\Vec{p}^{(i)}) \in \mathcal{W} ; \Vec{p}^{(i)} \in \mathbb{P}, i = 1,\dots, N_s \}$.

Note that $\Span \mathbb{S}$ is an approximation of $\mathcal{L}(\mathcal{P})$, which becomes better accordingly to how well $\mathbb{P}$ samples $\mathcal{P}$ and {\added depends upon the} complexity of the mapping $\mathcal{L}$. We are interested in finding the orthonormal basis $\mathcal{B}_{N_{rb}} := \{ \bs{\xi}_i \in \Span \mathbb{S}  \}_{i=1}^{N_{rb}}$, with $N_{rb}<<N_s$, so-called {\added Reduced Basis (RB)}, which is optimal in the sense that among all possible basis choices, $\mathcal{B}_{N_{rb}}$ is chosen such that $
\frac{1}{N_s} \sum_{i=1}^{N_s} \| \Vec{w}^{(i)} - \Pi_{N_{rb}} \Vec{w}^{(i)}\|_{\mathcal{W}}^2
$ is minimised, where $\Pi_{N_{rb}} (\cdot) = \sum_{i=1}^{N_{rb}} (\bs{\xi}_i, \cdot ) \bs{\xi}_i$ is {\added the orthogonal projection of a function on the space spanned by the RB}.

This task can be performed by the Proper Orthogonal Decomposition (POD) procedure detailed in Algorithm \ref{box:PODprocedure}, yielding the POD-error
\begin{align} \label{eq:errorPOD}
\mathcal{E}_{POD}^2(N_{rb})= \frac{1}{N_s} \sum_{i=1}^{N_s} \| \Vec{w}^{(i)} - \Pi_{N_{rb}} \Vec{w}^{(i)}\|_{\mathcal{W}}^2 =  \sum_{j=N_{rb}+1}^{N_s} \lambda_j,
\end{align}
which is related the ordered eigenvalues $\lambda_j$, for $j=1,\cdots,N_s$,  of the correlation matrix defined in \eqref{eq:correlationMatrix}. The size $N_{rb}$ is selected to attain a given tolerance $\mathtt{tol_{POD}} \in (0,1)$ in terms of relative POD-errors.

\begin{algorithmbox}[h!]
	Given $\mathbb{S}$, $\mathtt{tol_{POD}} \in (0,1)$, do:
	\begin{enumerate}
		\item Correlation matrix:
		\begin{align} \label{eq:correlationMatrix}
		\Mat{C} = \frac{1}{N_s} \sum_{i,j=1}^{N_s}(\Vec{w}(\Vec{p}^{(i)}) , \Vec{w}(\Vec{p}^{(j)}))_{\mathcal{W}} \Vec{e}_i \otimes \Vec{e}_j  
		\end{align}
		\item Eigenvalues analysis:
		\begin{align} \label{eq:eigenvalue}
		\Mat{C} \Vec{v}^i = \lambda_i \Vec{v}^i \quad \forall i= 1,2, \dots, N_s,
		\end{align}
		where the order $\lambda_1 \geq \lambda_2 \geq \dots \geq \lambda_{N_s}$ is assumed.
		\item Basis determination: The $i$-th basis vector is constructed as
		\begin{align}
		\bs{\xi}^i = \frac{1}{\sqrt{ \lambda_i  N_s}} \sum_{j=1}^{N_s} \Vec{v}^i_j \Vec{w}(\Vec{p}^{(j)}) \quad \forall i=1,2, \dots, N_s.
		\end{align}
		\item RB construction: choose the first $N_{rb} << N_s$ such that {\added $\displaystyle  1 - \frac{\sum_{i=1}^{N_{rb}} \lambda_i }{\sum_{i=1}^{N_s} \lambda_i} < \mathtt{tol_{POD}}$} and define
		\begin{align}
		\mathcal{B}_{N_{rb}} = \left \{  \bs{\xi}^i  \right\}_{i=1}^{N_{rb}}. 
		\end{align}
	\end{enumerate}
	\caption{{\added Proper Orthogonal Decomposition procedure to obtain the Reduced Basis with a given error tolerance.}}.
	\label{box:PODprocedure}
\end{algorithmbox}

In the remaining of the manuscript, {\addedSimone we target an RB representation of the traces of solutions on an internal boundary $\Gamma$ (properly specified later). Therefore we always set $\mathcal{W} = [L^2(\Gamma)]^d$ with the $[L^2(\Gamma)]^d$-inner product as our ambient Hilbert space.}

\begin{remark}
	Note that POD procedure format in Algorithm \ref{box:PODprocedure} is independent of discretisation. As result, the RB can be obtained even using snapshots with heterogeneous discretisation. Just for computational purposes, e.g., speed-up calculations, we use coincident meshes on $\partial \Omega_{\mu}^R$ (see Section \ref{sec:dataset}). 
\end{remark}

\FloatBarrier

\FloatBarrier

\vskip 2cm

\section{{\added Formulation of the method}} \label{sec:formulationMethod}

{\moved 
	{\added Before formulating the method itself, it is worth remembering that in order to have a computable microscopic cell problem, we need to truncate a problem posed essentially in an infinite domain, intrinsically connected to the choices below:}  
	\begin{enumerate}
		\item Boundary condition: As already commented, there is no inherently correct boundary condition to the corrector problem apart from periodic geometries. Notwithstanding, the error arising from the artificial imposition of BCs fades according to the MD size \citep{Yue2007}. 
		{\added 
			\item Microscale domain size: regarding the window of observation as a realisation of an underlying random variable that dictates the microstructure, the larger is such window, the better the heterogeneities are sampled. As result, as the microscale domain gets larger, the homogenised coefficients tend to the infinite infinite cell values.
		}
	\end{enumerate} 
	Concerning the microscale domain size decision, we have to leverage two competing characteristics, namely the computational cost and the accuracy in the two senses above. {\added The remaining option to increase the accuracy without interfering the computational cost is improving the BCs choice. Therefore,} the main contribution of this work is to bring errors due to BCs choices in small MDs down to the same levels of larger domains.

	The very central rationale behind the method is the reduction of the computational domain in the local problem, as depicted on the rightmost column of Figure \ref{fig:randommedia}. First, according to the final accuracy desired and the computational effort we can afford, two given sizes of microscale domains should be chosen (see Figure \ref{fig:enlargingrve}). {\addedSimone The larger cell is denoted $\Omega_{\mu}^{H}$ and is used for the computation of reference solutions, so-called High-Fidelity (HF) solutions, by using a chosen classical BC. On the other hand,  $\Omega^R_{\mu}$ is the smaller MD (reduced) where the corrector problem is solved }{\added with an enriched, non-classical, BC. Instrumental for formulation of the admissible space that encodes such novel BC is the parametrised PDE setting, introduced in the following section. Such a structure is an essential step towards the use RB-ROM framework, {\addedSimone which justification is also provided along this section.}  
		
		\begin{figure}
			\centering
			\includegraphics[width=0.42\linewidth]{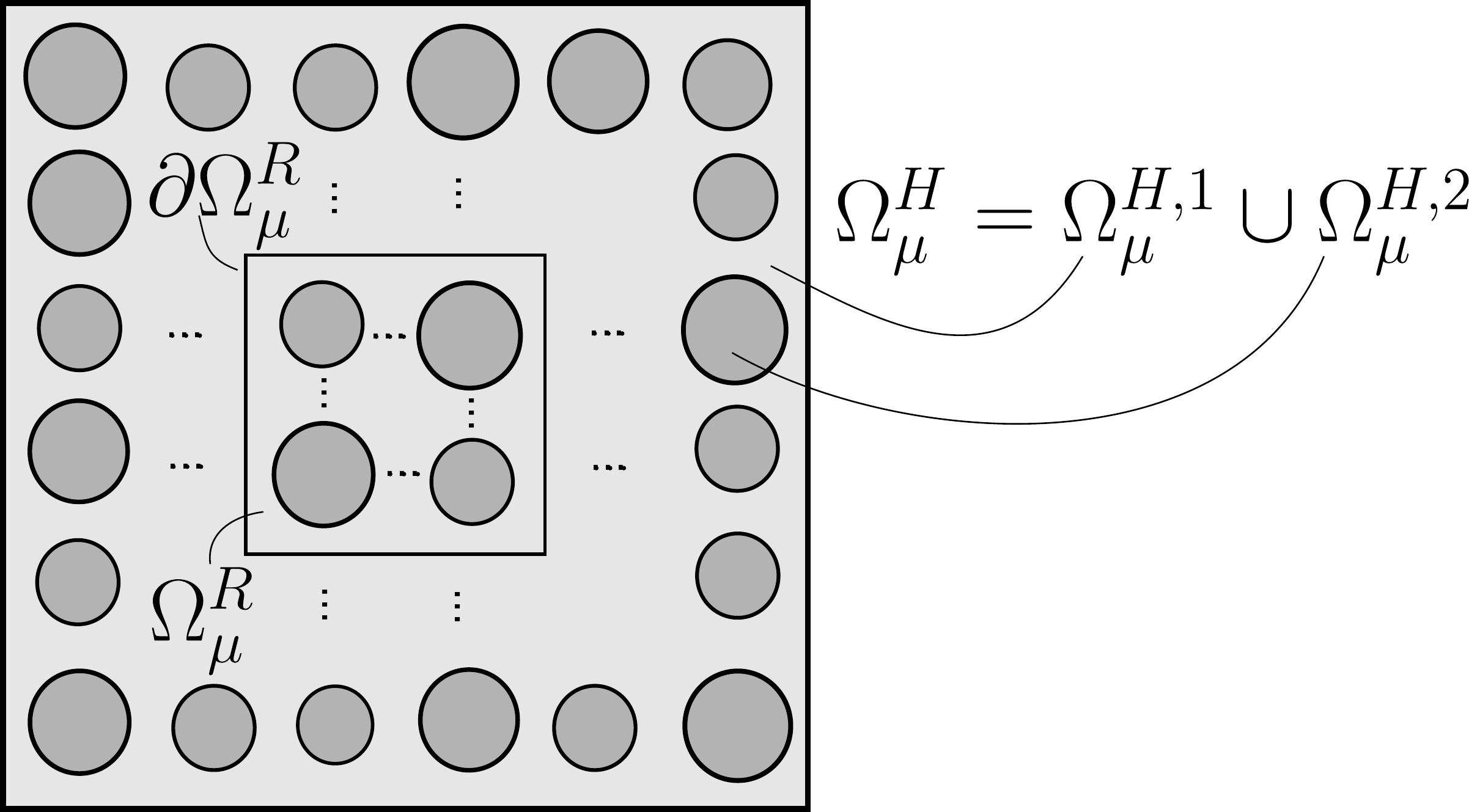}
			\caption{Larger and smaller microscale domains and boundary where the enhanced boundary condition is imposed.}
			\label{fig:enlargingrve}
		\end{figure}
		
	}
	
	\subsection{Reformulation of the multiscale problem in terms of parametrised PDE}
	\label{sec:multiscaleRB}
	
	In this section, we aim at providing a parametrised PDE structure to the corrector problem for solid mechanics. It is worth mentioning that the conversion of the general problem into its parametrised format encompasses necessarily a simplification. The practical effect is a new problem with a finite number of parameters that spans a controlled scenario of randomness, representative of the problem at hand. {\added However,} the general idea of the method remains true independently of the particular parametrisation chosen. 
	
	We are interested in approximating the fluctuation $\uf^{H}$, solution of Problem \ref{prob:correctorProblemSolidMechanics} in an enlarged domain $\Omega_{\mu}^{H}$, which depends on some parameters $\Vec{p} \in \mathcal{P} \subset \Real^{N_p}$, with $N_p$ the number of parameters. 
	Noticing that the triple $(\bs{\varepsilon}, \mathbb{C}_{\mu}, V_{\mu})$ fully defines the latter problem, the role of $\Vec{p}$ is to parametrise each of these objects, as respectively detailed in the following subsections. Importantly, we reserve the word \textit{parameter}, and notation $\Vec{p}$, to designate non-fixed properties, i.e., those that should be sampled. {\addedSimone For other numbers helping to define the problem, but that are fixed, we omit their explicit dependence}. 
	{\added Importantly, the space $V_{\mu}$ stands for different spaces depending if the problem is formulated on $\Omega_{\mu}^{H}$ or $\Omega^R_{\mu}$. For $\Omega_{\mu}^H$, due to lack of previous information, we fix $V_{\mu}(\Omega_{\mu}^{H}) = V_{\mu}^P(\Omega_{\mu}^{H})$, leading to Problem \ref{prob:correctorProblemSolidMechanics} 
		rewritten as follows.}
	\begin{problem}[High-fidelity parametrised corrector PDE] \label{prob:highFidelityParametricPDE}
		Given $\Vec{p} \in \mathcal{P}$, find $\uf^{H} \in V_{\mu}^P(\Omega_{\mu}^{H})$ such that
		\begin{align} \label{eq:parametricSolidMechanicsHF}
		a^{H}(\Vec{p} ; \uf^H,\Vec{v}) = b^{H}(\Vec{p};\Vec{v}) \quad \forall \Vec{v} \in V_{\mu}^P(\Omega_{\mu}^{H}),
		\end{align}
		where
		\begin{subequations}
			\begin{align}
			a^{H}(\Vec{p} ; \uf^H,\Vec{v}) = (\mathbb{C}_{\mu}(\Vec{p}) \nabla^s \uf^H, \nabla^s \Vec{v})_{L^2(\Omega_{\mu}^{H})}, \\ 
			b^{H}(\Vec{p} ; \Vec{v}) = -(\mathbb{C}_{\mu}(\Vec{p}) \bs{\varepsilon}(\Vec{p}), \nabla^s \Vec{v})_{L^2(\Omega_{\mu}^{H})}.
			\end{align}
		\end{subequations}
	\end{problem}

	To reduce the problem into $\Omega_{\mu}^R$, we choose $\Vec{v} \in V_{\mu}^P(\Omega_{\mu}^{H})$ such that $\Vec{v}|_{\Omega^R_{\mu}} \in V_{\mu}^L(\Omega_{\mu}^R)$ and $\Vec{v}|_{\Omega^{H}_{\mu} \backslash \Omega_{\mu}^R} = \Vec{0}$, {\added inducing the following trial set for the fluctuations 
		\begin{align} \label{eq:trialSpaceFluctuation}
		V_{\mu}(\Omega_{\mu}^R) = V_{\mu}(\Vec{p})(\Omega_{\mu}^R) := \left \{\bs{\eta} \in [H^1(\Omega_{\mu}^R)]^d ; \avg{\bs{\eta}}_{{\Omega}_{\mu}^R} = \avg{\uf^H(\Vec{p})}_{{\Omega}_{\mu}^R} ; T_{\partial \Omega_{\mu}^R} \bs{\eta} =  T_{\partial \Omega_{\mu}^R} \uf^H(\Vec{p}) \right \},  
		\end{align}
		where the dependence on $\Vec{p}$ is made explicitly and $T_{\partial \Omega_{\mu}^R} : [H^{1}(\Omega_{\mu}^{H})]^d \to [L^{2}(\partial \Omega_{\mu}^R)]^d$ is a trace operator. Using the previously defined spaces, the constrained version of Problem \ref{prob:highFidelityParametricPDE} on $\Omega^R_{\mu}$ read as:}

	\begin{problem}[Reduced parametrised corrector PDE] \label{prob:reducedParametricPDE}
		Given $\Vec{p} \in \mathcal{P}$, find $\uf^R \in V_{\mu}(\Omega_{\mu}^R)(\Vec{p})$ such that
		\begin{align} \label{eq:parametricSolidMechanics}
		a^R(\Vec{p} ; \uf^R ,\Vec{v}) = b^R(\Vec{p};\Vec{v}) \quad \forall \Vec{v} \in V_{\mu}^L(\Omega_{\mu}^R),
		\end{align}
		where 
		\begin{subequations}
			\begin{align} \label{eq:bilinearReduced}
			a^R(\Vec{p} ; \uf^R,\Vec{v}) = (\mathbb{C}_{\mu}(\Vec{p}) \nabla^s \uf^R, \nabla^s \Vec{v})_{L^2(\Omega_{\mu}^R)}, \\ 
			b^R(\Vec{p} ; \Vec{v}) = -(\mathbb{C}_{\mu}(\Vec{p}) \bs{\varepsilon}(\Vec{p}), \nabla^s \Vec{v})_{L^2(\Omega_{\mu}^R)}.
			\end{align}
		\end{subequations}
	\end{problem}
	
	{\added Importantly, we should notice that $\uf^H|_{\Omega_{\mu}^R} = \uf^R$ holds by construction and, as already commented in Remark \ref{rmk:translation}, the homogenised stress is invariant upon translations. Therefore, it is convenient to define
		\begin{align} \label{eq:trialSpaceFluctuation_redefined}
		V_{\mu}^*(\Omega_{\mu}^R) = V^*_{\mu}(\Vec{p})(\Omega_{\mu}^R) := \left \{\bs{\eta} \in V_{\mu}^{ZA}(\Omega_{\mu}^R) ; T_{\partial \Omega_{\mu}^R} \bs{\eta} =  \Vec{w} (\Vec{p}) \right \},
		\end{align}
		where $\Vec{w}(\Vec{p}) \in L^2(\partial \Omega_{\mu}^R)$ denotes
		\begin{align}
		\Vec{w}(\Vec{p}) = T_{\partial \Omega_{\mu}^R} \uf^{H}(\Vec{p}) - \avg{\uf^{H}(\Vec{p})}_{\Omega_{\mu}^R}. 
		\end{align}
		Notice that Problem \ref{prob:reducedParametricPDE} can be rephrased replacing  $V_{\mu}(\Vec{p})(\Omega_{\mu}^R)$ by $V^*_{\mu}(\Vec{p})(\Omega_{\mu}^R)$. As commented, the new solution is mechanically equivalent (leads to the same homogenised stress). As a result, the only information needed coming from the solution of Problem \ref{prob:highFidelityParametricPDE} is its trace to the internal boundary $\partial \Omega_{\mu}^R$ {\addedSimone and the necessary translation to obtain a zero-averaged fluctuation field in $\Omega_{\mu}^R$. For convenience, we rearranged this operation in the goal solution $\Vec{w}$, as above.} 
		
		The key idea of proposed method is the following: Instead of solving the larger Problem \ref{prob:highFidelityParametricPDE} to obtain $\Vec{w}$, we aim at obtaining an approximation $\Vec{w}^{\mathcal{N}}$ by a method that combines RB-ROM and DNN. This approximation yields to the definition of the novel set $V^{\mathcal{N}}_{\mu}(\Omega_{\mu}^R)(\Vec{p})$ (or just $V^{\mathcal{N}}_{\mu}$). Our final goal is to show that $V^{\mathcal{N}}_{\mu}$ outperforms the classical choices of BCs.

		Finally, before proceeding with the parametrisation of $\mathbb{C}_{\mu}$ and $\bs{\varepsilon}$ in Sections \ref{sec:constitutiveLaw_parametrisation} and \ref{sec:strain_parametrisation}, respectively, it is useful to assume the split $\Vec{p} = (\Vec{p}_c,\Vec{p}_{\bs{\varepsilon}})$, where $\mathbb{C}_{\mu}(\Vec{p}) = \mathbb{C}_{\mu}(\Vec{p}_c)$
		and $\bs{\varepsilon}(\Vec{p}) = \bs{\varepsilon}(\Vec{p}_{\bs{\varepsilon}})$. The explicitly definition of $V^{\mathcal{N}}_{\mu}(\Omega_{\mu}^R)(\Vec{p})$ is postponed to Section \ref{sec:fluctuationSpace_parametrisation}
	}
	
	%

	\subsubsection{Parametrisation of $\mathbb{C}_{\mu}$} \label{sec:constitutiveLaw_parametrisation} 
	In our context, parametrising the constitutive law is intrinsically linked with the parametrisation of microstructural geometrical features of material constituents. In this work, since the number of parameters have a direct impact on the quality of sampling, if we want to keep the number of samples limited, we aim at keeping the geometry and number of materials simple. 
	
	We admit a partition of two materials on $\Omega^{H}_{\mu}$ where those properties are constant by parts. Let be $\Omega_{\mu}^{H,1}$, $\Omega_{\mu}^{H,2}$ be such that $\overline{\Omega}_{\mu}^{H,1} \cup \overline{\Omega}_{\mu}^{H,2} = \overline{\Omega}_{\mu}^{H}$ and $\Omega_{\mu}^{H,1} \cap \Omega_{\mu}^{H,2} = \emptyset$, i.e., a partition of $\Omega^{H}_{\mu}$ (see Figure \ref{fig:enlargingrve}). 
	Each of these materials are linear isotropic and thus modelled by the Hookean law characterised by the two spatially varying Lam{\'e} coefficients $\lambda_{\mu}$ and $G_{\mu}$, or conversely the pair Young modulus and Poisson ratio. It reads
	\begin{align}
	\mathbb{C}_{\mu}(\y) = \lambda_{\mu}(\y) \ten{I} \otimes \ten{I} + 2 G_{\mu} (\y) \mathbb{I}^s \quad ,\forall \y \in \Omega_{\mu}^{H}.
	\end{align}
	
	Let us assume $\lambda_{\mu}$ and $G_{\mu}$ are linked with $\Omega_{\mu}^{H,1}$ and $\Omega_{\mu}^{H,2}$ via a constant $\gamma > 0$. It useful to define $\chi_{\gamma} : \Omega_{\mu}^{H} \to \Real$ as follows  
	\begin{align}
	\chi_{\gamma}(\y) =
	\begin{cases} 
	\gamma, \quad \text{if } \y \in \Omega_{\mu}^{H,2} \\
	1, \quad \text{otherwise} 
	\end{cases},
	\end{align}
	so $\lambda_{\mu}(\y) := \chi_{\gamma}(\y) \lambda_{\mu}^1$, $G_{\mu}(\y) = \chi_{\gamma}(\y) G_{\mu}^1$, i.e., $\mathbb{C}_{\mu}(\y) = \chi_{\gamma}(\y) \mathbb{C}^1_{\mu}$.
	
	For the sake of simplicity, let us consider {\addedSimone the two-dimensional setting with $\Omega^{H}_{\mu} = (0,1)^2$  and $\Omega_{\mu}^{H,2}$ composed by union of disjoint balls, distributed on a lattice}. Let $N_b$ the number of balls, $\y_i^c$ and $r_i$, $i = 1, 2,  \dots, N_b$ the centre positions and radii, respectively. Hence $\Omega_{\mu}^{H,2} = \bigcup_{i=1}^{N_b} B(\y_i^c,r_i)$, where $B(\y_0,r) := \{ \y \in \Real^{d} ; \|\y - \y_0 \| < r \}$ denotes the ball centred in $\y_0$ with radius $r$. Furthermore, we assume a regular grid of $\sqrt{N_b} \times \sqrt{N_b}$ balls, for $N_b$ a perfect square number, with positions as follows
	\begin{align} \label{eq:grid}
	\y_{(i-1) \sqrt{N_b} + j + 1}^c = \frac{2}{\sqrt{N_b}} ( 2 j - 1, 2 i - 1), \quad i, j = 1,2,\dots, \sqrt{N_b}.
	\end{align}  
	Regarding other properties used in the model, for the sake of simplicity, let us take $\lambda_{\mu}^1$, $G_{\mu}^1$, and $\gamma$ fixed. {\addedSimone Notice that each element of  $(\lambda_{\mu}^1, G_{\mu}^1, \gamma, \y_1^c, \dots, \y_1^{N_b})$ has fixed value and is not considered as a parameter of the problem} . In turn, the \textit{de facto} parameters ruling $\mathbb{C}_{\mu}$ are $\Vec{p}_c =  (r_1, r_2, \dots, r_{N_b})$.
	
	\subsubsection{Parametrisation of $\bs{\varepsilon}$}
	\label{sec:strain_parametrisation}
	
	Differently {\addedSimone as in the last section}, the parametrisation of macroscale strain is straightforward and we do not need to assume any simplification, besides the symmetry of $\bs{\varepsilon}$. Focusing on $d=2$, we introduce $\Vec{p}_{\bs{\varepsilon}} = (p_{\bs{\varepsilon},1}, p_{\bs{\varepsilon},2}, p_{\bs{\varepsilon},3})$  and the mapping $\Vec{p}_{\bs{\varepsilon}} \mapsto \bs{\varepsilon}(\Vec{p}_{\bs{\varepsilon}}) = \bmat p_{\bs{\varepsilon},1} & \frac{1}{2} p_{\bs{\varepsilon},3} \\ \frac{1}{2}p_{\bs{\varepsilon},3} &  p_{\bs{\varepsilon},2} \emat$. This choice of parametrisation is not unique, but it turns out that in this specific choice each component of $\Vec{p}_{\bs{\varepsilon}}$ is physically meaningful: i) $p_{\bs{\varepsilon},{1,2}}$ is the axial strain along the $1,2$ directions and ii) $p_{\bs{\varepsilon},3}$ is the shear strain. {\addedSimone Indeed, $\Vec{p}_{\bs{\varepsilon}}$ corresponds to the so-called Voigt's notation of $\bs{\varepsilon}$, here represented by the mapping $\mathrm{Voigt} : \Real^{2,2}_{\mathrm{sym}} \to \Real^3$, such that $ \bs{\varepsilon} \mapsto \Vec{p}_{\bs{\varepsilon}} = \mathrm{Voigt}(\bs{\varepsilon}) = (\varepsilon_{11}, \varepsilon_{22}, \varepsilon_{12} + \varepsilon_{21})$}. Finally, the extension for $d =3$ is immediate, but will not be used in this manuscript.
	
	Taking advantage of the problem structure in the right-hand side of \eqref{eq:parametricSolidMechanicsHF}, we can exploit its linearity with respect to $\Vec{p}_{\bs{\varepsilon}}$ by writing
	\begin{align} 
	b^{H}(\Vec{p};\Vec{v}) = \sum_{i=1}^3 p_{\bs{\varepsilon},i} b^{H}((\Vec{p}_c, \Vec{e}_i); \Vec{v}).
	\end{align}
	Accordingly, it is easy to see that linearity is also valid for $\Vec{w}$ as follows
	\begin{align} \label{eq:unitaryAuxSolution}
	\Vec{w}(\Vec{p}) &= \sum_{i=1}^3 p_{\bs{\varepsilon},i} \Vec{w}^{(i)}(\Vec{p}_c),
	\end{align}
	where $\Vec{w}^{(i)}(\Vec{p}_c) := T_{\partial \Omega_{\mu}} \uf^{H}((\Vec{p}_c,\Vec{e}_i)), \text{ for } i = 1, 2, 3$.  
	We exploit this decomposition in Section \ref{sec:RBfindingStrategy} to reduce the number of RB constructions. 
	
	\subsubsection{Parametrisation of {\added $V_{\mu}^{\mathcal{N}}$}}
	\label{sec:fluctuationSpace_parametrisation}
	
	By defining a learning-based set $V_{\mu}^{\mathcal{N}}$, the goal is to {\added characterise} the set where fluctuations live as accurate as possible (not necessarily as large as possible) based on the prior knowledge decoded in the RB-ROM and DNN model. 
	{\addedSimone  Indeed, we must recall that $V_{\mu}^{\mathcal{N}}$ aims to approximate $V_{\mu}(\Omega_{\mu}^R)$ in Problem \ref{prob:reducedParametricPDE}.} To this end, let us assume we dispose of an $N_{rb}$-dimensional RB $\mathcal{B}_{N_{rb}} = \{ \bs{\xi}_i; \bs{\xi}_i \in L^2(\partial \Omega_{\mu}^R)  ,i = 1,\cdots, N_{rb} \}$ and the mapping $\mathcal{N} : \Real^{N_p} \to \Real^{N_{rb}}$. We postulate 
	\begin{align} \label{eq:parametrisedVmu}
	V^{\mathcal{N}}_{\mu}(\Vec{p})&= \left \{ \bs{\eta} \in V_{\mu}^{ZA}; T_{\partial \Omega_{\mu}^R}\bs{\eta} = \sum_{i=1}^{N_{rb}}[ \mathcal{N}(\Vec{p}) ]_i \bs{\xi}_i \right \} = V_{\mu}^L + \uf^{\mathcal{N}}(\Vec{p}), 
	\end{align}
	where $\uf^{\mathcal{N}}(\Vec{p}) \in V_{\mu}^{ZA}$ is a continuous extension of $\sum_{i=1}^{N_{rb}}[ \mathcal{N}(\Vec{p}) ]_i \bs{\xi}_i \in L^2(\partial \Omega_{\mu}^R)$ to {\addedSimone $\Omega_{\mu}^R$}, whose construction is detailed next. 
	
	Solving Problem \ref{prob:reducedParametricPDE} with solutions in $V_{\mu}(\Vec{p})$ involves finding $\uf^0(\Vec{p}) \in V_{\mu}^L$ and the explicit form of $\uf^{\mathcal{N}}(\Vec{p}) \in V_{\mu}^{\mathcal{N}}$, such that $\uf^0(\Vec{p}) + \uf^{\mathcal{N}}(\Vec{p}) \in V_{\mu}(\Vec{p})$. Hence, Problem \ref{prob:reducedParametricPDE} encompasses the solution of two sub-problems: i) find $\uf^0(\Vec{p}) \in V_{\mu}^{\mathcal{N}}$ such that 
	\begin{align} \label{eq:subproblemForUf0}
	a^R(\Vec{p} ; \uf^0 ,\Vec{v}) = b^R(\Vec{p};\Vec{v}) \quad \forall \Vec{v} \in V_{\mu}^L,
	\end{align}
	and ii) find $\uf^{\mathcal{N}}(\Vec{p}) \in V_{\mu}^L$ such that 
	\begin{align} \label{eq:subproblemForUfN}
	a^R(\Vec{p} ; \uf^{\mathcal{N}}(\Vec{p}) ,\Vec{v}) = 0 \quad \forall \Vec{v} \in V_{\mu}^L.
	\end{align}
	This decoupling is possible thanks to the linearity of the problem. 
	\begin{remark} \label{rmk:not_a_subset}
		{\addedSimone It is worth mentioning that in general $V^{\mathcal{N}}_{\mu}$ is not a subset of $V_{\mu}^M$. In Appendix \ref{sec:admissibility} we prove this does not impose any limitation to the method and that, in fact, we can equivalently rewrite the original problem using a set $(V^{\mathcal{N}}_{\mu})^* \subset (V_{\mu}^M)$ just by incorporating a new term on the right-hand side of the variational formulation in Problem \ref{prob:reducedParametricPDE}.
		} 
	\end{remark}
	
	\begin{remark}
		In practice, it is more convenient to solve Problem \ref{prob:reducedParametricPDE} {\added(using $V_{\mu}^{\mathcal{N}}$)} at once rather \eqref{eq:subproblemForUf0} and \eqref{eq:subproblemForUfN} separately. However, this split motivates and justifies the structure  $V_{\mu}^{\mathcal{N}}$ postulated in \eqref{eq:parametrisedVmu}.
	\end{remark}

	{\moved 
		\subsection{Overview of the method} \label{sec:overviewMethod}
		
		Before digging into methodological details {\added of (Deep) Neural Networks (Section \ref{sec:designingNeuralNet})}, it is useful to start a short overview of the method {\added we propose. For easiness of notation, thereafter we denote our method as \DeepBC, which stands for Deep Boundary, as reference of a deep learning methodology to enrich multiscale problems with more accurate BCs. It is worth highlighting that} the basic requirements behind the development of \DeepBC and the strategies adopted to address each of them are the following:
		\begin{enumerate}
			\item Computational efficiency: coupled multiscale simulations are computationally demanding by their very nature. Since the local-problem solving is embarrassingly parallel, this step is the most critical regarding the optimisation of the method. Therefore, {\addedSimone such speed-up} is mandatory and our main strategy is to reduce the microscale domain. 
			\item Accuracy: To keep the same accuracy of more computationally demanding larger microscale domains, our strategy focuses on the choice of a more precise BCs for smaller microscale domains.
			\item Exploitation of \textit{a priori} knowledge: To be able to enrich the corrector problem with more suitable BCs, we exploit previous knowledge from pre-computed accurate solutions, also called high-fidelity (HF) solutions, which are necessarily more computationally demanding. This is performed by fully exploiting the potential of Reduced Order Modelling using Reduced Basis Method (RB-ROM) in combination with (Deep) Neural Networks (DNN). {\added The use of RB-ROM allows for a more compact representation of the BC, reducing the size of the needed DNN, by automatically ranking the basis in a meaningful order of importance, which will be also exploited to suitably calibrate the training process of the DNN in next section.}
			\item {\added Easiness of} implementation: The very mathematical {\addedSimone structure} of the local problem is unchanged, only the BCs are touched. Our method leads to non-homogeneous Dirichlet BCs, {\addedSimone not leading to any additional difficulties in terms of implementation.}
		\end{enumerate}

		{\added As usual in RB-ROM strategies, \DeepBC can be structured into \textit{offline} and \textit{online} phases}. Having chosen some meaningful microstructural parametrisation for the HF problem at hand, the main goal of the \textit{offline phase} is to find a mapping that translates parameters {\addedSimone from the}  high-fidelity problem into an accurate BC of the problem solved in the \textit{online phase}. It consists in the following pipeline: i) the dataset acquisition by simulating those models for some representative sampling of the parameter space, ii) the RB construction using the previous dataset, iii) extraction of relevant features of the dataset through the orthogonal projections onto this RB and, iv) the training of a DNN model aiming to predict these projections. Mathematically, this procedure results in a) $\mathcal{B}_{N_{rb}} = \{\bs{\xi}_i \in [L^2(\partial \Omega_{\mu})]^d\}_{i=1}^{N_{rb}}$, a RB with $N_{rb}$ functions and b) the mapping $\mathcal{N} : \Real^{N_p} \to \Real^{N_{rb}}$, where $N_p$ is the number of parameters of the high-fidelity model. {\added As for the \textit{online} phase, we can predict the Dirichlet-like BC encoding these two previous outputs into $V_{\mu}^{\mathcal{N}}$ for a given parameter of the HF problem $\Vec{p} \in \Real^{N_p}$. Finally, we proceed by solving the reduced corrector problem followed by some homogenisation post-processing. The general outline of the method is depicted in Algorithm \ref{box:workflow}}.

		\begin{algorithmbox}
			\begin{minipage}[t]{.45\textwidth}
				\textbf{Offline phase:}
				\begin{enumerate}
					\item Dataset generation
					\item Reduced Basis Construction
					\item Features extractions (Input-Output)
					\item DNN training
				\end{enumerate}
			\end{minipage}
			\hfill
			\noindent
			\begin{minipage}[t]{.45\textwidth}
				\textbf{Online phase:}	
				\begin{enumerate}
					\item Given a microstructure and a micro-strain: Prediction of BC using the DNN
					\item {\coupling Solving the microscale problem (Problem \ref{prob:reducedParametricPDE})} using the BC of step 1.
					\item  Homogenisation of the obtained solution
				\end{enumerate}
			\end{minipage}
			\caption{General workflow of the \DeepBC method detailing Offline and Online phases.}
			\label{box:workflow}
		\end{algorithmbox}

		It is worth mentioning that the RB-ROM and DNN concentrate most of the computational efforts in the \textit{offline phase} rather than in the \textit{online phase}.  To generate the dataset, some parametrisation of the microstructure should be chosen, together with the range of these parameters, their probability distributions, and sample strategy. This is a problem-dependent choice and we discuss in Section \ref{sec:dataset} the examples treated in the manuscript. In turn, the final cost added in comparison to a standard coupled multiscale simulation is the evaluation of a DNN model to determine the BC that feeds the corrector problem solver. The DNN evaluation step (prediction) is many orders of magnitude less computationally demanding than the corrector problem {\addedSimone numerical solution}, by using FEM for example.
		
		\section{Designing the Neural Network} \label{sec:designingNeuralNet}
		
		In this section, we introduce the DNN model to predict the BCs. As mentioned before, we aim at characterising the nonlinear mapping from the parameter space into RB coefficients $\mathcal{N} : \Real^{N_p} \to \Real^{N_{rb}}$. In terms of DNN architecture, our model is based on two instances of the standard fully connected Multilayer Perceptron (MLP) model, {\added which is revisited in Section \ref{sec:notationNN},} combined seamlessly to exploit: i) the splitting between macro-strain and microscale constitutive parameters and ii) the problem symmetries. {\addedSimone We also explore such principles for dataset and RB generation (Section \ref{sec:RBfindingStrategy}), which, in fact, naturally leads to a specific format for combining MLP submodels, as explained in Section \ref{sec:combiningNNmodels}.} Regarding the loss function, we also adopt a non-standard setting motivated from the RB-ROM (see Section \ref{sec:lossFunction}).

		\subsection{Basic notation for Deep Neural Networks } \label{sec:notationNN}
		
		\begin{figure}
			\centering
			\includegraphics[width=0.5\linewidth]{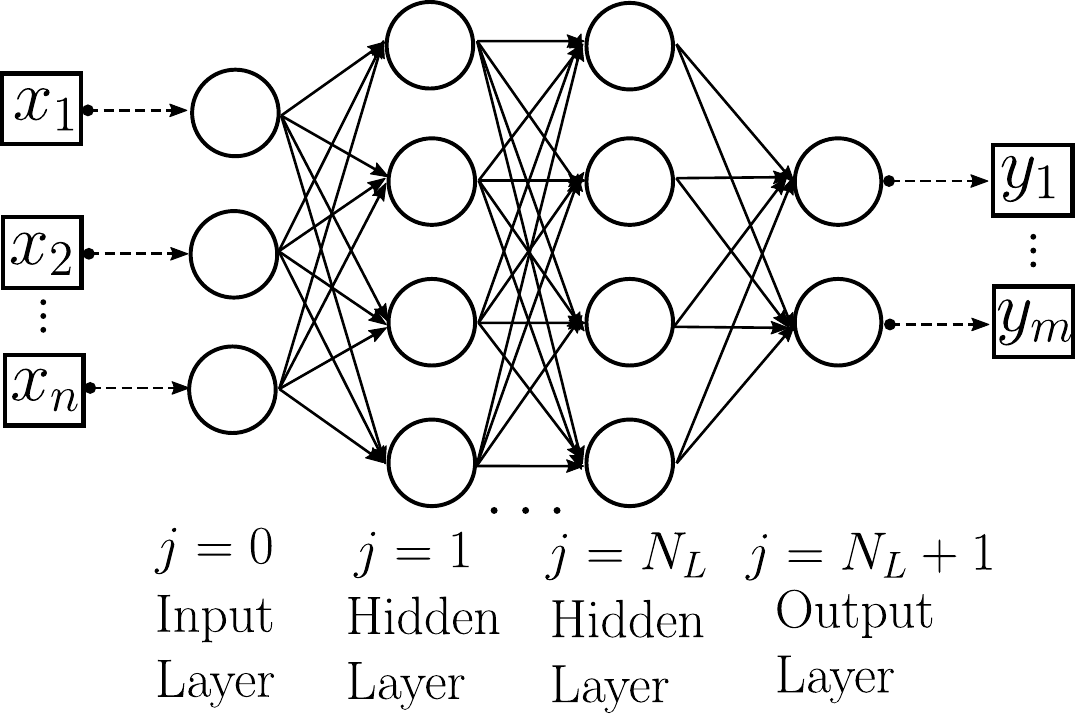}
			\caption{Multilayer Perceptron Neural Network {\added with $N_L$ hidden layers}.}
			\label{fig:mlp}
		\end{figure}

		Here we review some basic concepts and notation for supervised learning using DNN. For the interested reader, we refer to \citep{Higham2019,Goodfellow2016} for a deeper presentation. 
		
		Let $D = \{(\bs{\alpha}^{(i)},\bs{\beta}^{(i)}) \in \Real^m \times \Real^n; i = 1,\dots, N_s\}$ be a \textit{dataset} with $N_s$ samples (or snapshots). We aim at finding the mapping {\addedSimone $\mathcal{N}(\bs{\theta};\cdot) : \Real^m \to \Real^n$} such that it minimises the regularised empirical risk
		\begin{align} \label{eq:minimisationNN}
		\min_{\bs{\theta}} \frac{1}{N_s} \sum_{i=1}^{N_s} 
		\mathcal{L}(\mathcal{N}({\bs{\theta}}; \bs{\alpha}^{(i)}), \bs{\beta}^{(i)}) \, + \lambda \mathcal{R}(\bs{\theta}), 
		\end{align}
		where $\mathcal{L} : \Real^m \times \Real^m \to \Real^+$ is the so-called cost or loss function, and $\mathcal{R} : \Real^p \times \Real^m \to \Real^+$ is a regularising function controlled by a intensity factor $\lambda \in \Real^+$. The vector $\bs{\theta}$ should be understood as a collection of parameters that defines the mapping {\addedSimone $\mathcal{N}(\bs{\theta};\cdot)$, which is taken accordingly with the usual definition of the MLP architecture for DNN \citep{Goodfellow2016}. The size of $\bs{\theta}$, and thus the complexity of the DNN, depends on the pre-fixed number of hidden layers $N_L$ and the quantity of neurons in each layer (see Figure \ref{fig:mlp} to help with the interpretation).}
		
		%
		%
		%

		Moreover, concerning the practical implementation of the training, i.e., the optimisation procedure to solve \eqref{eq:minimisationNN}, we use the ADAM version of the stochastic gradient with mini-batches, and, as regularisation techniques, we make use of dropout \citep{dropout}, $l^2$ regularisation, and early stop based on the validation error \citep{Goodfellow2016}. The exact numerical values used in our examples are specified in the appropriate section.  
		
		\FloatBarrier 
		
		{\moved 
			\subsection{{\added Dataset and RB construction strategies}} \label{sec:RBfindingStrategy}
			
			As seen in Section \ref{sec:strain_parametrisation}, the trace $\Vec{w}$ of the fluctuations on $\partial \Omega_{\mu}^R$, is the linear combination of auxiliary solutions obtained to each load direction ($\Vec{w}^{(I)}$), for $I= 1, 2, 3$. Obtaining them involve fewer parameters, since the loads are fixed in unitary directions, which lead to a reduced parameter space to sample. In other words, exploiting this property in the dataset generation means that for a fixed number of high-fidelity snapshots, we can sample better a smaller parameter set $\mathcal{P}_r \subset \Real^{N_b}$ instead of $\mathcal{P} \subset \Real^{N_p}$, where $N_p = N_b + 3$ (in 2D). In this regard, it is also important to maximise the covering of $\mathcal{P}_r$ by using a smart sampling strategy, for example one of the Quasi-Monte Carlo Family such as the Latin Hypercube Sampling (LHS) used in this work. 
			
			Furthermore, in terms of mechanical load, solutions  $\Vec{w}^{(1)}$ and $\Vec{w}^{(2)}$ as in \eqref{eq:unitaryAuxSolution} are equivalent in the sense one can be obtained by rotation of another, as depicted in Figure \ref{fig:schemaVerticalEvaluation}. {\added Hence, choosing $N_{rb}'$ basis functions to satisfy some accuracy criteria, we obtain $\mathcal{B}^{(I)}_{N_{rb}'} =  \left\{\bs{\xi}^{(I)}_i \right\}_{i=1}^{N_{rb}'}$ associated to $\Vec{w}^{(I)}$ for $I = 1,3$, by applying the POD procedure of Algorithm \ref{box:PODprocedure}. The missing basis $\mathcal{B}^{(2)}_{N_{rb}'}$ can be obtained by $\mathcal{B}^{(2)}_{N_{rb}'} = \ten{Q}_{\frac{\pi}{2}} \mathcal{B}^{(1)}_{N_{rb}'}$, where the anti-clockwise rotation of $\pi/2 \mathrm{rad}$  is applied for each function in the basis by means of the orthonormal matrix $\ten{Q}_{\frac{\pi}{2}}$. On the other hand, for a given $\Vec{p}_c \in \Real^{N_c}$, the solution $\Vec{w}^{(2)}(\Vec{p}_c)$ transforms according to 	
				\begin{align} \label{eq:rotation90}
				\Vec{w}^{(2)}(\Vec{p}_c) = \Vec{w}((\Vec{p}_c, \Vec{e}_2))
				= \ten{Q}_{\frac{\pi}{2}} \Vec{w}((\bs{\pi}^{-1}_{\frac{\pi}{2}}(\Vec{p}_c), \Vec{e}_1))= \ten{Q}_{\frac{\pi}{2}}  \Vec{w}^{(1)}(\bs{\pi}^{-1}_{\frac{\pi}{2}}(\Vec{p}_c)), 
				\end{align}
				where $\bs{\pi}^{-1}_{\frac{\pi}{2}}$ denotes a permutation on the parameter vector entries to reflect a clockwise rotation of $\pi/2 \, \mathrm{rad}$ on the microstructure pattern (see Figure \ref{fig:schemaVerticalEvaluation}). The actual format for the permutation is problem- and implementation-dependent. In the next section, we retake these ideas to construct the DNN model that predicts the BCs.}
			
			\begin{figure}[h]
				\centering
				\includegraphics[width=0.5\linewidth]{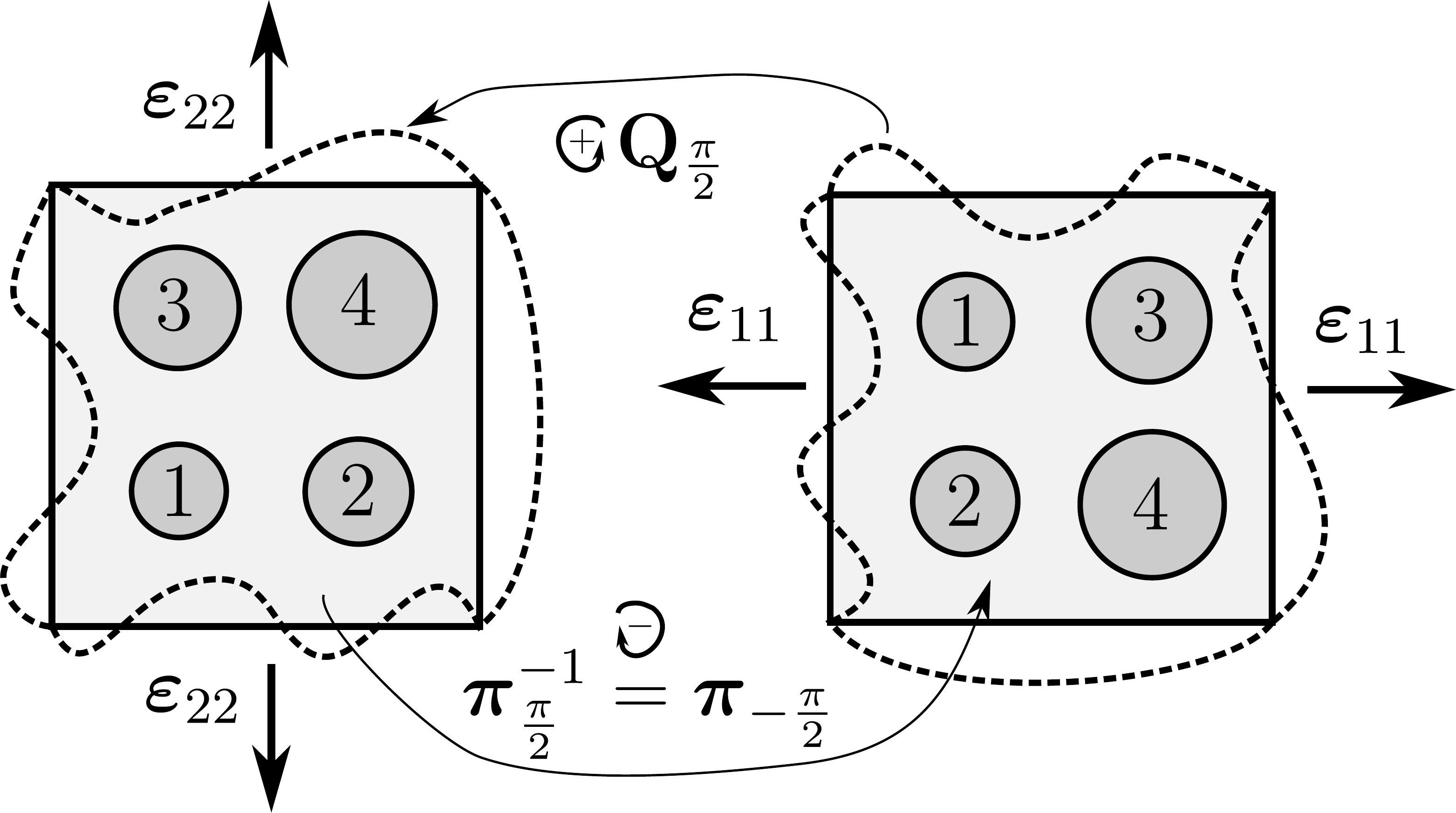}
				\caption{Symmetry in load-geometry relation. Here $\bs{\pi}^{-1}_{\frac{\pi}{2}}((1,2,3,4)) = \bs{\pi}_{-\frac{\pi}{2}}((1,2,3,4)) = (2,4,1,3).$ }
				\label{fig:schemaVerticalEvaluation} 
			\end{figure}

		}

		\subsection{\DeepBCNN model: definition and combination of the axial and shear DNN submodels} \label{sec:combiningNNmodels}
		
		\begin{figure}[!h]
			\centering
			\includegraphics[width=0.75\linewidth]{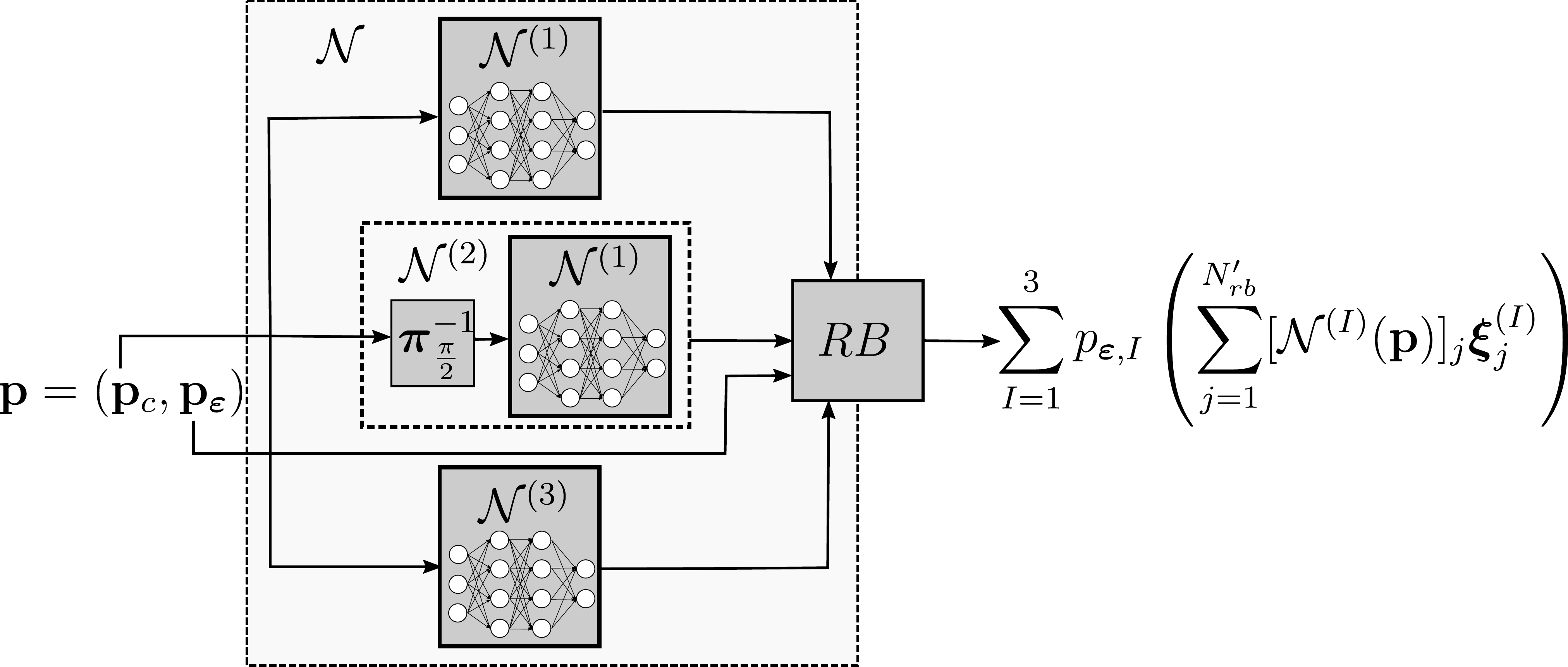}
			\caption{\DeepBCNN architecture: {\addedSimone exploiting the physical structure of the problem using a combination of MLP submodels and RB functions.} }
			\label{fig:NNarchitecture}
		\end{figure}
		
		As already commented, before defining the DNN model, we need to obtain a DNN submodel for the axial unitary macro-strain (without loss of generality, chosen here as the horizontal one)  and another DNN submodel for shear unitary macro-strain. Following the previous conventions, we denote them as $\mathcal{N}^{(1)}, \mathcal{N}^{(3)}: \Real^{N_b} \to \Real^{N_{rb}'}$, {\added for the axially horizontal and shear loading, respectively. In turn, for the axially vertical loading, we invoke \eqref{eq:rotation90} to define $\mathcal{N}^{(2)}(\Vec{p}) = \mathcal{N}^{(1)}\left(\Vec{\bs{\pi}_{\frac{\pi}{2}}^{-1}(p)}\right)$.}
		
		
		{\added 
			Finally, the BC prediction for some $\Vec{p}\in \Real^{N_p}$ reads as
			\begin{align}
			\uf^{\mathcal{N}}(\Vec{p})|_{\partial \Omega_{\mu}^R}  = \sum_{I=1}^3 p_{\bs{\varepsilon},I} \left( \sum_{j=1}^{N'_{rb}} [\mathcal{N}^{(I)}(\Vec{p_c})]_j \bs{\xi}_{j}^{(I)} \right) \in L^2(\partial \Omega_{\mu}^R).
			\end{align}
		}
		
		{\moved
			It is important to mention that by taking into account the physical structure of the problem into the DNN architecture ({\added reviewed in Figure \ref{fig:NNarchitecture}}), we are able to: i) minimise the size of datasets needed since the size of parameter space (input) is smaller for each submodel; ii) use fewer number of neurons and layers (smaller complexity of the DNN architecture) since the prediction of directional load-wise solutions are also less complex; iii) profit from the full advantage of the {\addedSimone RB} representation {\addedSimone due to their specialisation according to the load}; iv) use fewer RB modes for the same accuracy; v) naturally take into account symmetries concerning horizontal and vertical loads and vi) training of 2 instead of 3 submodels. All points mentioned here are intrinsically related and overall render better accuracy and smaller computational burden in the dataset generation and training. 
		}
		
		\FloatBarrier

		\subsection{DNN-POD error analysis and Loss Function choice} \label{sec:lossFunction}
		
		In this section, we precisely define the total error committed in the BC prediction, which accounts for the intrinsic truncation error when a RB is chosen (POD error) and also the error due to the DNN prediction. The latter naturally yields to the natural choice for the loss-function $\mathcal{L}$ used in \eqref{eq:minimisationNN}. All notation here should be understood for one specific load, axial or shear, with indices omitted for convenience.
		
		First, let us consider some $\mathcal{B}_{N_{rb}}$ and $\mathbb{S}$ given. {\added For some $i=1,\dots,N_s$ and target HF solution $\Vec{w}^{(i)} \in \mathbb{S} \subset L^2(\partial \Omega_{\mu}^R)$, let us denote its projection upon $\Span{\mathcal{B}_{N_{rb}}}$ as $\Pi_{N_{rb}} \Vec{w}^{(i)} := \sum_{j=1}^{N_{rb}} \beta_j^{(i)} \bs{\xi}_j$, where $\beta^{(i)}_j$ represents its projection along $\bs{\xi}_j$. Analogously, we denote the DNN prediction of $\Vec{w}^{(i)}$ as $ \hat{\Vec{w}}^{(i)}  := \sum_{i=1}^{N_{rb}} \hat{\beta}^{(i)}_j \bs{\xi}_j \in \Span{\mathcal{B}_{N_{rb}}}$, where $\hat{\beta}^{(i)}_j$ follows the same interpretation of $\beta^{(i)}_j$}. 
		
		A natural choice to measure the expected error is to consider the empirical averaged error over all elements of $\mathbb{S}$ namely
		\begin{align} \label{eq:totalError}
		\mathcal{E}^2_{T}(N_{rb}) = \frac{1}{N_s} \sum_{i=1}^{N_s} \|\hat{\Vec{w}}^{(i)} - \Vec{w}^{(i)} \|^2_{L^2(\partial \Omega^R_{\mu})}, 
		\end{align}
		where the sub-index $(\cdot)_{T}$ stands for total error, as opposed to the contributions to the POD and DNN errors defined next. To this aim, using the fact that $\Vec{w}^{(i)} - \Pi_{N_{rb}} (\Vec{w}^{(i)}) \in (\Span \mathcal{B}_{N_{rb}})^\perp$ and that $\mathcal{B}_{N_{rb}}$ is orthonormal, it yields    
		\begin{align}
		\|\hat{\Vec{w}}^{(i)} - \Vec{w}^{(i)} \|^2_{L^2(\partial \Omega_{\mu}^R)} &= 
		\|\hat{\Vec{w}}^{(i)} -  \Vec{w}^{(i)} + \Pi_{N_{rb}}(\Vec{w}^{(i)}) - \Pi_{N_{rb}}(\Vec{w}^{(i)}) \|^2_{L^2(\partial \Omega_{\mu}^R)} \nonumber \\
		&=\|\hat{\Vec{w}}^{(i)} - \Pi_{N_{rb}}(\Vec{w}^{(i)}) \|^2_{L^2(\partial \Omega_{\mu}^R)} +  \|\Vec{w}^{(i)} - \Pi_{N_{rb}}(\Vec{w}^{(i)}) \|^2_{L^2(\partial \Omega_{\mu}^R)} \nonumber \\ 
		&=\|\hat{\bs{\beta}}^{(i)} - \bs{\beta}^{(i)} \|^2_2 +  \|\Vec{w}^{(i)} - \Pi_{N_{rb}}(\Vec{w}^{(i)}) \|^2_{L^2(\partial \Omega_{\mu}^R)}, \quad {\added \text{for}\, i = 1, \dots, N_s.}
		\end{align}
		{\added The boldface notation for $\hat{\bs{\beta}}^{(i)}$ and $\bs{\beta}^{(i)}$ simply mean the component-wise  collection for their respective projections for a snapshot $i$}. Note that the second term of the above expression when averaged over the $N_s$ snapshots corresponds to the POD mean squared error ($\mathcal{E}_{POD}^2$) already defined in \eqref{eq:errorPOD}. On the other hand, the first term accounts for the DNN error, in which the mean squared version is defined below
		\begin{align}
		\mathcal{E}^2_{DNN}(N_{rb}) &:= \frac{1}{N_s} \sum_{i=1}^{N_s} \|\hat{\bs{\beta}}^{(i)} - \bs{\beta}^{(i)} \|^2_2.
		\end{align}
		The expression above provides a physically-based formula for the loss function to be considered in \eqref{eq:minimisationNN}. More precisely,
		\begin{align}
		\mathcal{L}(\hat{\bs{\beta}}^{(i)}, \bs{\beta}^{(i)}) := \frac{1}{N_{rb}} \sum_{j=1}^{N_{rb}} \omega_j (\hat{\bar{\beta}}_j^{(i)} - \bar{\beta}_j^{(i)} )^2, 
		\end{align}
		where we have used the min-max scaling $\textstyle \overline{(\cdot)} = \frac{(\cdot) - \beta_j^{min}}{\beta_j^{max} - \beta_j^{min}}$ (for $\hat{\beta}^{(i)}_j$ and $\beta^{(i)}_j$), $\beta_j^{min} = \min_{i} \beta_j^{(i)}, \beta_j^{max} = \max_{i} \beta_j^{(i)}, \omega_j = (\beta_j^{max} - \beta_j^{max})^2$, for $j=1,\dots,N_{rb}$.
		
		Finally, we can rewrite \eqref{eq:totalError} as
		\begin{align} \label{eq:totalError_split}
		\mathcal{E}^2_{T}(N_{rb}) =  
		\mathcal{E}^2_{DNN}(N_{rb}) + \mathcal{E}^2_{POD}(N_{rb}),
		\end{align}
		where we should notice that while $\mathcal{E}_{POD}^2$ has a monotonous decreasing (w.r.t. $N_{rb}$) estimation given by \eqref{eq:errorPOD}, the $\mathcal{E}^2_{DNN}$ depends on the DNN architecture, initialisation, optimisation algorithm, etc. Therefore, to assess an optimal $N_{rb}$ for a given DNN architecture  we should study the trade-off between these two errors. This is explored in Section \ref{sec:dataset}.

		{\added It is worth highlighting that} the RB-based strategy adopted also induces a specific choice of weighted mean-squared error, where weights are given by eigenvalues found in the POD procedure, which are ranked by magnitude. This choice also leads to a natural total estimator in the {\added $L^2(\partial \Omega_{\mu}^R)$-norm}, which is the intrinsic choice to measure the BC discrepancies, that links the loss function reached in the training and the error committed by POD truncation itself. Importantly, each submodel can be trained separately and be put together only for prediction purposes, which is interesting from the computational point of view.

		\section{DNN model construction: datasets, RB, and training} \label{sec:dataset}
		
		In this section, we specify the choice of numerical parameters adopted to assess our method. We choose to constrain the more general setting of Sections \ref{sec:multiscaleRB} and \ref{sec:designingNeuralNet}. As already discussed, the DNN model construction encompasses the RB extraction and the training of a DNN for a given dataset, denoted here as \textit{training/RB} dataset, with size $N^{trb}_s = 51200$. Moreover, from the point of view of DNN training, we also need to simulate auxiliary datasets for \textit{validation} and \textit{testing} reasons, in which the size is assumed to be respectively $20\%$ ($N^{val}_s = 10240$) and $10\%$ ($N^{tes}_s = 5120$) of the \textit{training/RB} dataset size. Importantly, each dataset is considered for the axial and shear loads separately, as summarised in Table \ref{tab:datasets}. 
		The fixed coefficients used in all simulations are presented in Table \ref{tab:fixedCoefficients}, while the random sampling of the inclusion sizes, which are the random parameters of the problem, are {\addedSimone shown} next.

		Concerning the inclusions, they are placed in a $6\times6$ grid (see \eqref{eq:grid}), so $N_b = 36$, with radii computed as follows 
		\begin{align} \label{eq:radiusSampling}
		r^{(i)}_j = \exp(a + \theta^{(i)}_j b) \in [r_{min}, r_{max}], \qquad  j= 1, \dots, N_b, \; i = 1, \dots, N_s,
		\end{align}
		where $a = \log(r_{max} r_{min})$, $b = \log(\tfrac{r_{max}}{r_{min}})$, $r_{max} =  0.4 \tfrac{L_{\mu}^{HF}}{\sqrt{N_b}} = \tfrac{1}{15}$ and $r_{min} = 0.1 \tfrac{L_{\mu}^{HF}}{\sqrt{N_b}} = \tfrac{1}{60}$, and $\theta^{(i)}_j \in [-1,1]$ is a random variable detailed next. The strategy \eqref{eq:radiusSampling} is used by \citep{Steinmann2019} in a similar scenario and has the advantages of being empirically inspired and also a simple manner to constrain the value range for the radii. Regarding the sampling, we use the LHS with an underlying uniform distribution to {\addedSimone sample $\theta^{(i)}_j$ for $j= 1, \dots, N_b, \; i = 1, \dots, N_s$. Notice that the LHS requires \textit{a priori} knowledge of $N_s$ and $N_b$, since it aims at achieving an optimal coverage of the high-dimensional parameter while keeping the number of samples small. Indeed, each realisation is not independent from each other.} We use the default LHS implementation provided \texttt{scikit-opt} library \citep{scikitopt}, with three different seeds for random generation and sample sizes according to Table \ref{tab:datasets}, yielding the samples $\mathbb{P}^{(1)}$, $\mathbb{P}^{(2)}$, and $\mathbb{P}^{(3)}$. Examples of some realisations of microstructures are found in Figure \ref{fig:microstrustures}.
		
		It is widely known in the multiscale solid mechanics literature 
		that the volume fraction $\phi = \tfrac{\sum_{i=1}^{N_b} \pi r_i^2}{(L_{\mu}^{HF})^2}$ is a random variable that plays a crucial role in the final solution. Consequently it is worth looking at how this variable is indirectly sampled by the assumed random generation algorithm, which focuses only on individual radii. In Figure \ref{fig:histogramsVolumeFraction}, checking the histograms of $\phi$ for $\mathbb{P}^{(1)}$, $\mathbb{P}^{(2)}$, and $\mathbb{P}^{(3)}$, we can observe that, although not implicitly enforced, the resulting underlying probability distribution seems to be consistent {\addedSimone along the several datasets. We can also }perfectly use another sample strategy, provided that the resulting volume fractions remain close to the range that the model was designed for. 
		
		{\added 
			For obtaining the HF solutions, we use the FEM with quadratic polynomials in triangular meshes. For the mesh generation using \textit{gmsh} \citep{gmsh}, we have imposed $20$ equally spaced divisions for each internal boundary segment (bottom, right, top and left of $\partial \Omega_{\mu}^R$), making a total of $160$ degrees of freedom on the boundary, for all snapshots. The overall unstructured mesh on $\Omega_{\mu}^H$ is then generated using the corresponding characteristic length size of the internal boundary. Accordingly, the number of functions in the RB will be at most the total number of degree of freedom on the boundary, and the associated POD-error should vanish close to this value. Indeed, the largest RB considered here is $140$, which actually gives practically zero POD-error. 
		}
		
		\begin{table}
			\begin{center}
				\begin{tabular}{|c|c|c|c|c|}
					\hline 
					Index & Load & Use & Size ($N_s$) & Sample \\ 
					\hline 
					$1$ & Axial &  Training/RB & $N_s^{trb} = 51200$ & $\mathbb{P}^{(1)}$ (LHS($N_s^{trb}$, \texttt{seed} $=1$))\\ 
					\hline 
					$2$ & Shear &  Training/RB & $N_s^{trb} = 51200$ & $\mathbb{P}^{(1)}$ (LHS($N_s^{trb}$, \texttt{seed} $=1$)) \\ 
					\hline 
					$3$ & Axial &  Validation  & $N_s^{val} = 10240$ & $\mathbb{P}^{(2)}$ (LHS($N_s^{val}$, \texttt{seed} $=2$)) \\ 
					\hline 
					$4$ & Shear & Validation & $N_s^{val} = 10240$ & $\mathbb{P}^{(2)}$ (LHS($N_s^{val}$, \texttt{seed} $=2$)) \\ 
					\hline 
					$5$ & Axial &  Testing & $N_s^{tes} = 5120$ & $\mathbb{P}^{(3)}$ (LHS($N_s^{tes}$, \texttt{seed} $=3$)) \\ 
					\hline 
					$6$ & Shear & Testing & $N_s^{tes} = 5120$ & $\mathbb{P}^{(3)}$ (LHS($N_s^{tes}$, \texttt{seed} $=3$)) \\ 
					\hline 
				\end{tabular} 
			\end{center}
			\caption{Summary datasets for the different combinations of ''Load/Use'' scenarios. Each ''Use'' is associated with one LHS run for a given \texttt{seed} and number of samples.}
			\label{tab:datasets}
		\end{table}

		\begin{table}
			\begin{center}
				\begin{tabular}{|c|c|c|}
					\hline 
					$\lambda_{\mu}^1$ & $G_{\mu}^1$ & $\gamma$ \\ 
					\hline 
					$0.576923$ & $0.384615$ & $10$\\ 
					\hline 
				\end{tabular} 
			\end{center}
			\caption{From left to right: the two Lam{\'e} coefficients ($\lambda_{\mu}^1$ and $G_{\mu}^1$) for the matrix and the {\addedSimone multiplicative constant factor ($\gamma$) in the inclusion applied for both previous coefficients}. The values for $\lambda_{\mu}^1$ and $G_{\mu}^1$ corresponds to $E = 1.0$ (units disregarded), the Young modulus, and $\nu = 0.3$, the Poisson ratio, using standard transformation expressions for Hookean isotropic media.}
			\label{tab:fixedCoefficients}
		\end{table}

		\begin{figure}[htp]
			\centering
			\subfigure[$\Vec{p}_c^{(63)} \in \mathbb{P}^{(1)}$. ]{
				\includegraphics[width=0.25\textwidth]{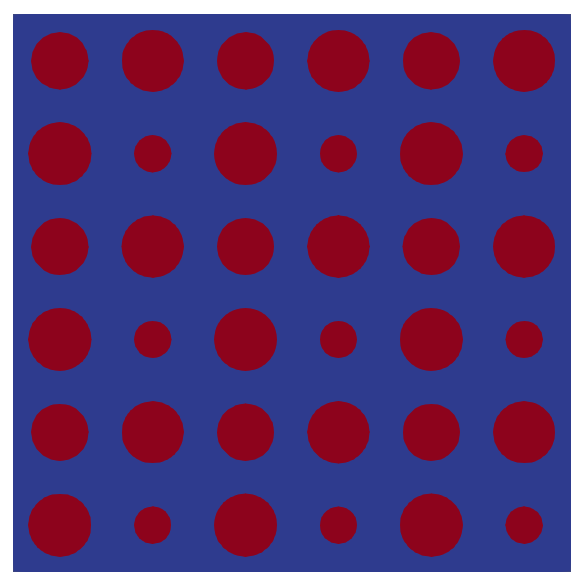}
			}
			\subfigure[$\Vec{p}_c^{(24567)} \in \mathbb{P}^{(1)}$.]{
				\includegraphics[width=0.25\textwidth]{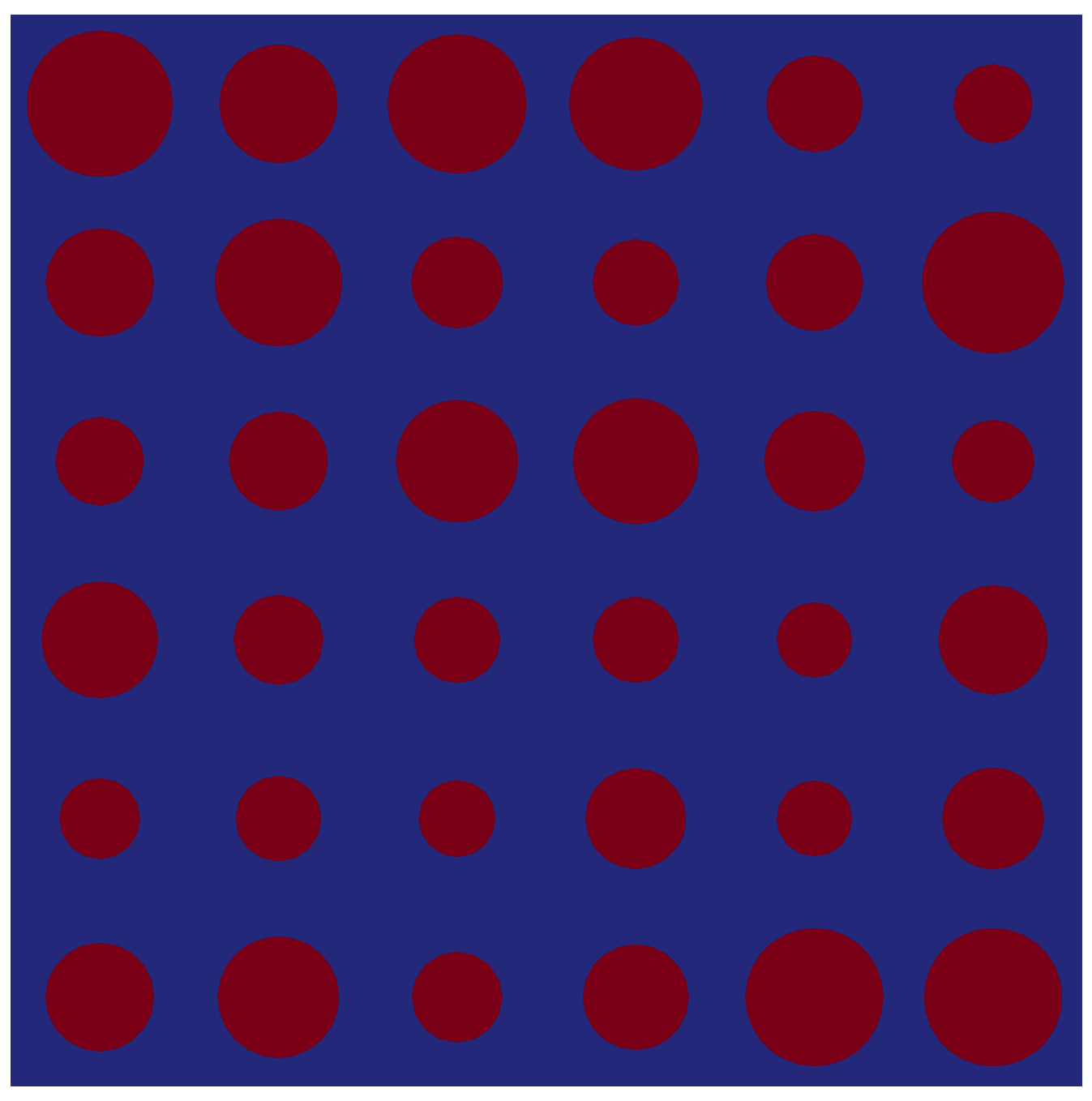}
			}
			\subfigure[$\Vec{p}_c^{(4231)} \in \mathbb{P}^{(2)}$.]{
				\includegraphics[width=0.25\textwidth]{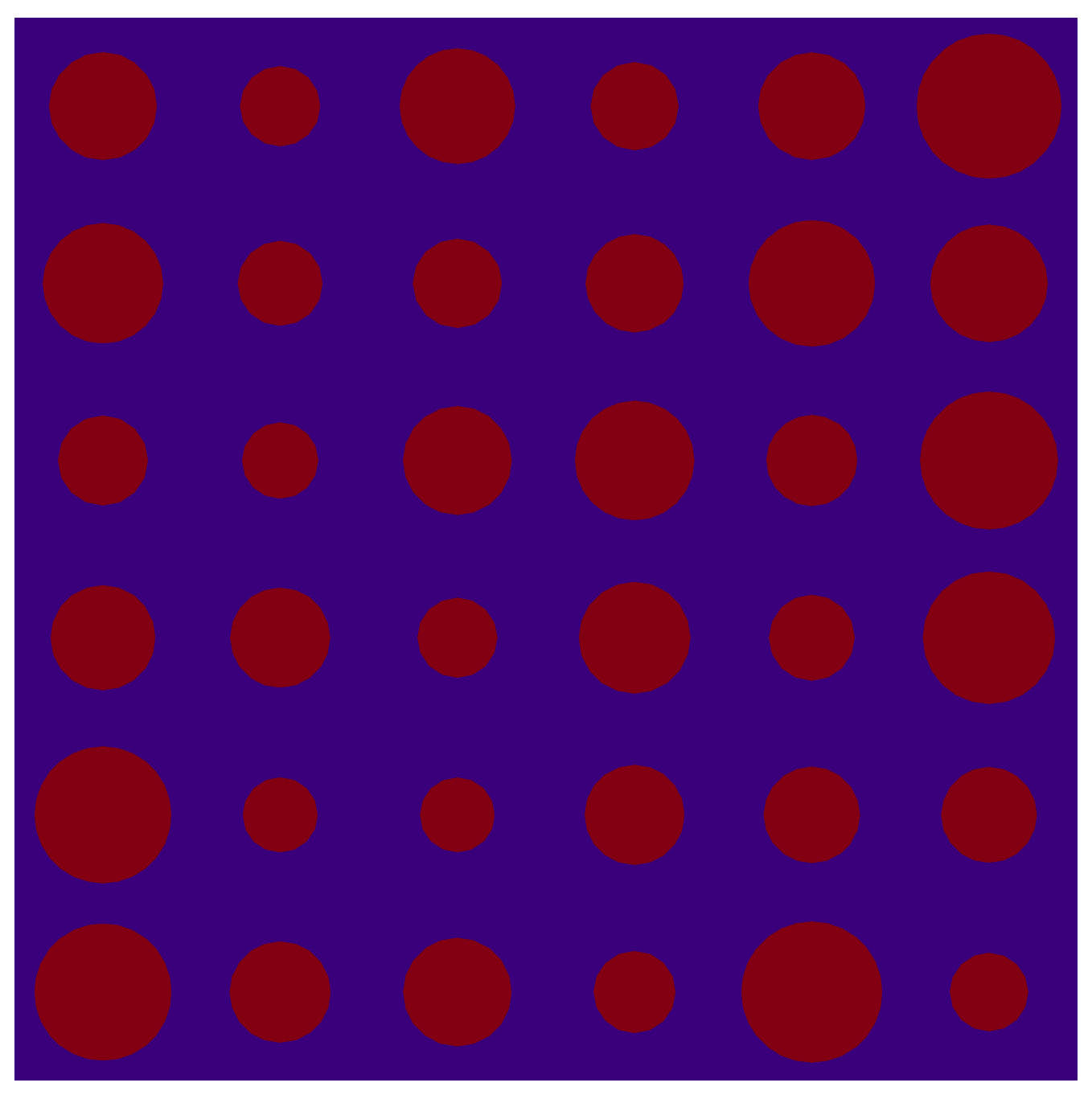}
			}
			\caption{Some realisations of microstructures for the high-fidelity model showing the different phases of $\Omega_{\mu}^{H}$: $\Omega_{\mu}^{H,1}$ in blue and $\Omega_{\mu}^{H,2}$ in red.}
			\label{fig:microstrustures}
		\end{figure}

		\begin{figure}
			\centering
			\subfigure{
				\includegraphics[width=0.4\textwidth]{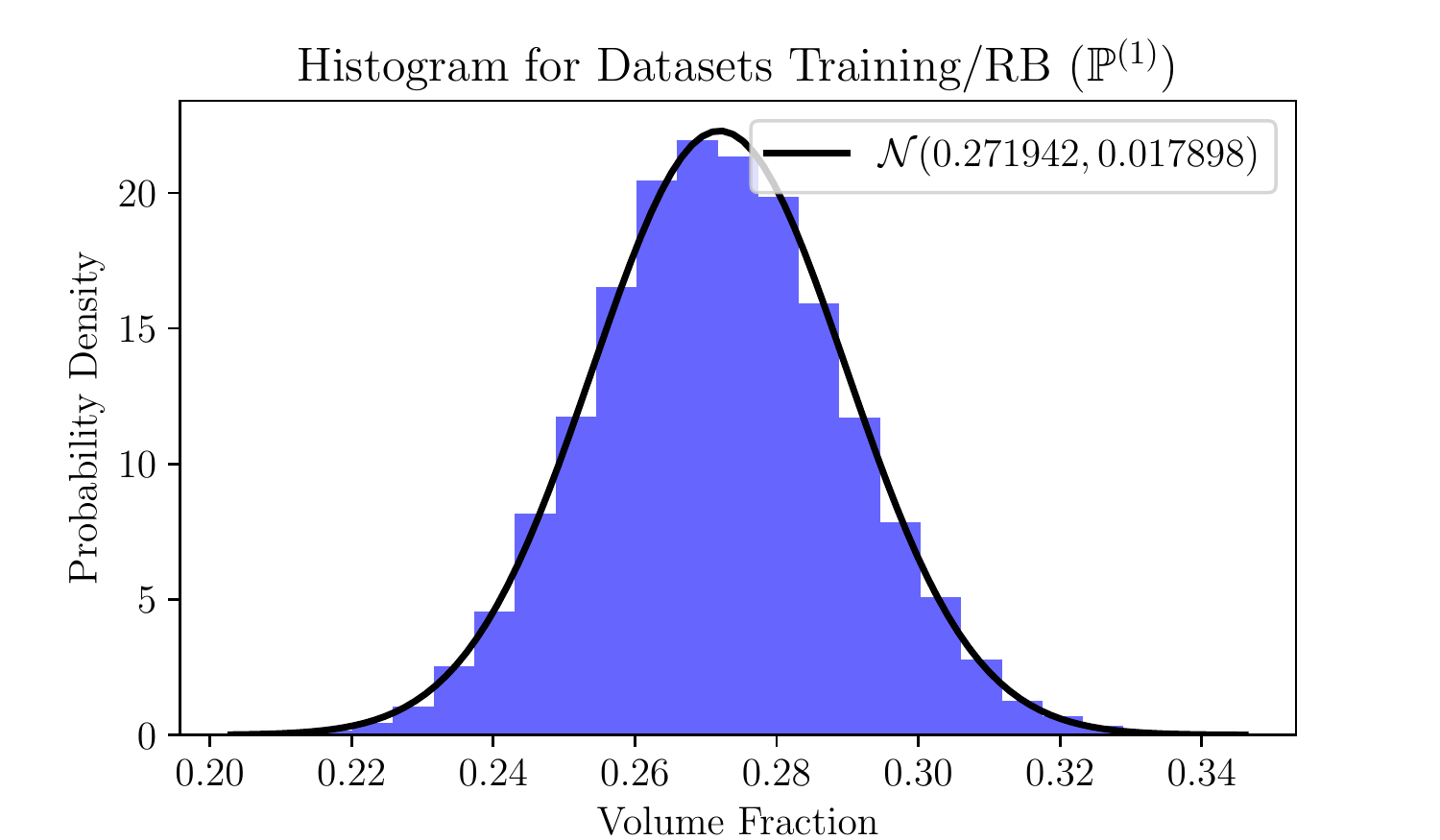}
			}
			\subfigure{
				\includegraphics[width=0.4\textwidth]{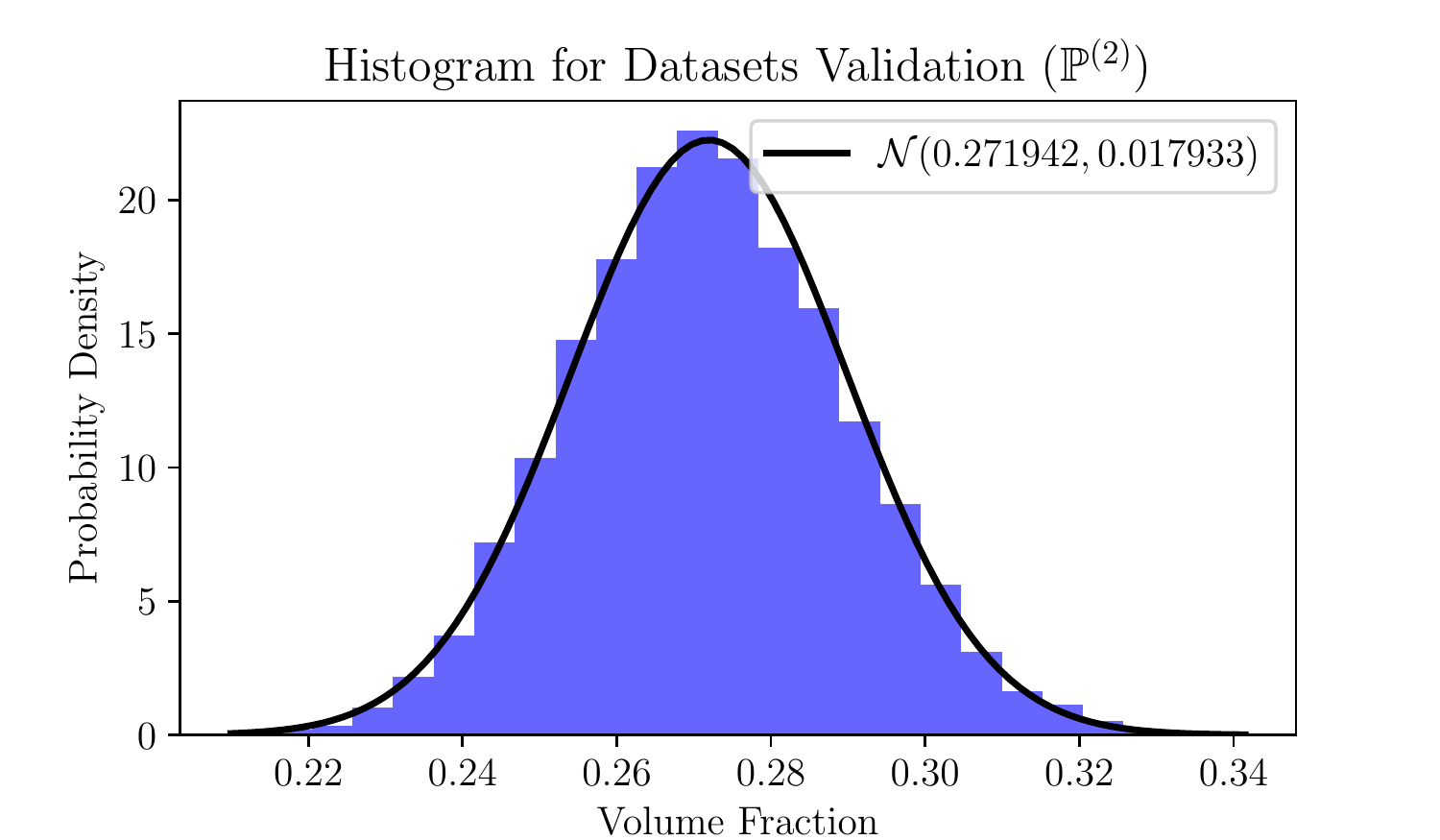}
			}
			\subfigure{
				\includegraphics[width=0.4\textwidth]{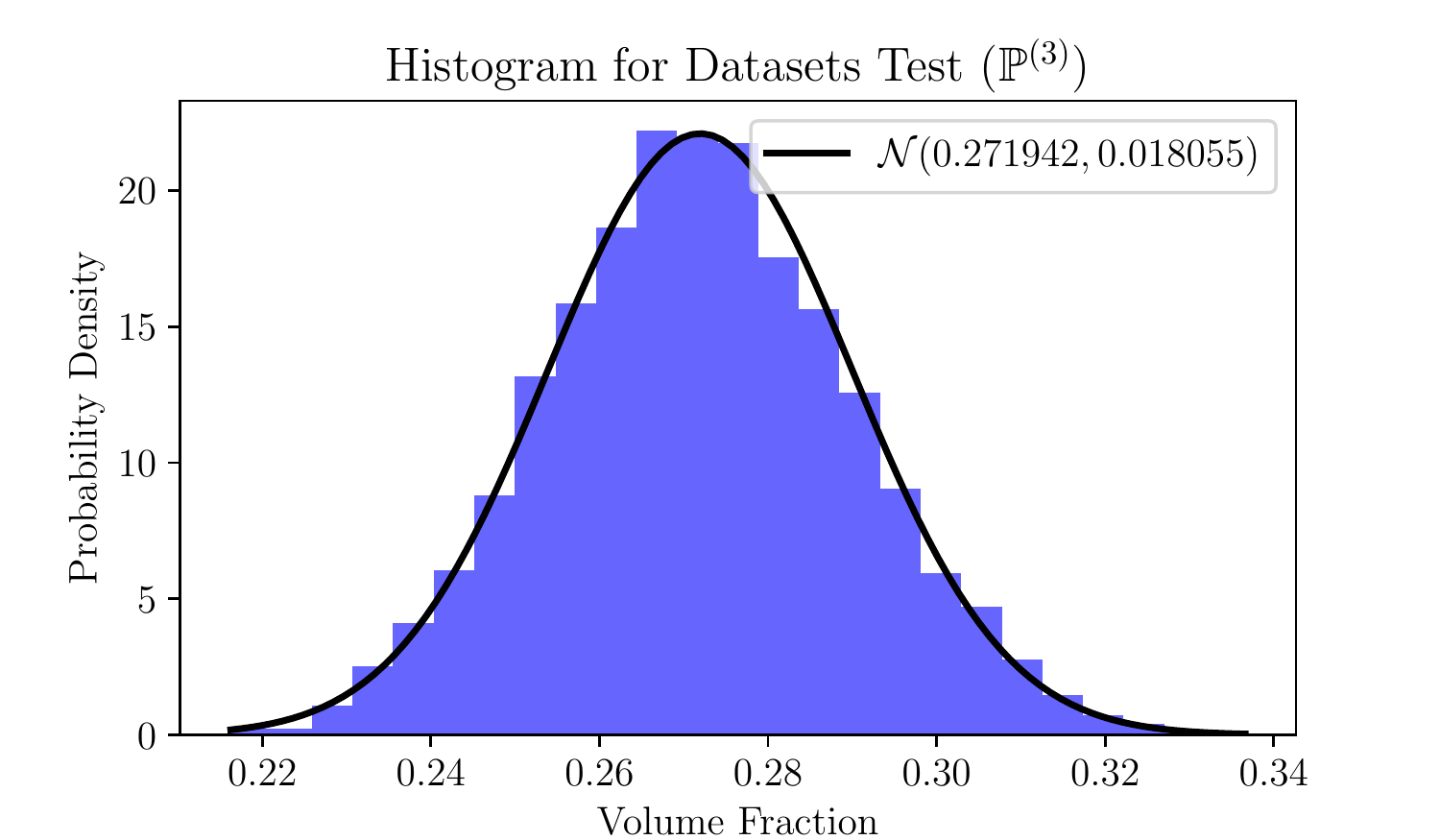}
			}
			\caption{Volume fraction histograms of three different types of datasets. The notation $\mathcal{N}(\cdot,\cdot)$ represents the fitted Normal distribution (do not mix with a DNN).}
			\label{fig:histogramsVolumeFraction}
		\end{figure}
		
		{\added
			Concerning the computational time, all simulations were performed in a cluster with 80 (160 threads) Intel(R) Xeon(R) Gold 6148 CPU @ 2.40GHz cores and the results of solving each HF simulation are resumed in Table \ref{tab:computationalTimeDataset}. The \textit{Total Time} displayed is an estimate based on the statistics computed over all snapshots, i.e., the number of snapshots times the mean time spent per snapshot. The Wall Time is computed supposing the use of $32$ cores and such that the Total Time is split equally. These estimates make sense only in embarrassingly parallel scenarios, as it is the case of the problem at hand. Notwithstanding, the ultimate time expended depend on many other factors, which are out of the scope of our analysis, such as the implementation, number of cores/threads, memory available, etc. Furthermore, the time expended in the RB-ROM steps and the mesh generation is negligible compared to the Wall Time for the \textit{training/RB} dataset, and for this reason, we do not report them. 
		}
		
		\begin{table}[h!]
			\begin{center}
				\begin{tabular}{|c|c|c|c|}
					\hline 
					& Training/RB & Validation & Test \\ 
					\hline 
					Time per snapshot(s) & $43.16 \pm 5.34$ & $45.11 \pm 10.60 $ & $47.91 \pm 12.96 $ \\
					Total Time (h) & $613.8$ & $128.3$ & $68.14$ \\
					Wall Time (h) ($32$ Cores) & $19.18$ & $4.009$ & $2.129$ \\
					\hline 
				\end{tabular} 
			\end{center}
			\caption{Computational cost of the dataset generation.}
			\label{tab:computationalTimeDataset}
		\end{table}
		
		{\added
			Finally, by applying POD procedure of Algorithm \ref{box:PODprocedure} to the datasets $1$ and $2$, we can evaluate the committed error in the sense  of \eqref{eq:errorPOD}, as shown in Figure \ref{fig:poderror}. This error represents the best expected approximation of the BC in the sense of $L^2(\partial \Omega_{\mu}^R)$, i.e., with no error due to the DNN approximation. We can verify these two errors are comparable, with a slight advantage to the shear over the axial model, for a given number of RB functions. As result, we have adopted the same DNN architecture and number of RB in both models, since they have comparable complexity. 
			
			\begin{figure}
				\centering
				\includegraphics[width=0.48\linewidth]{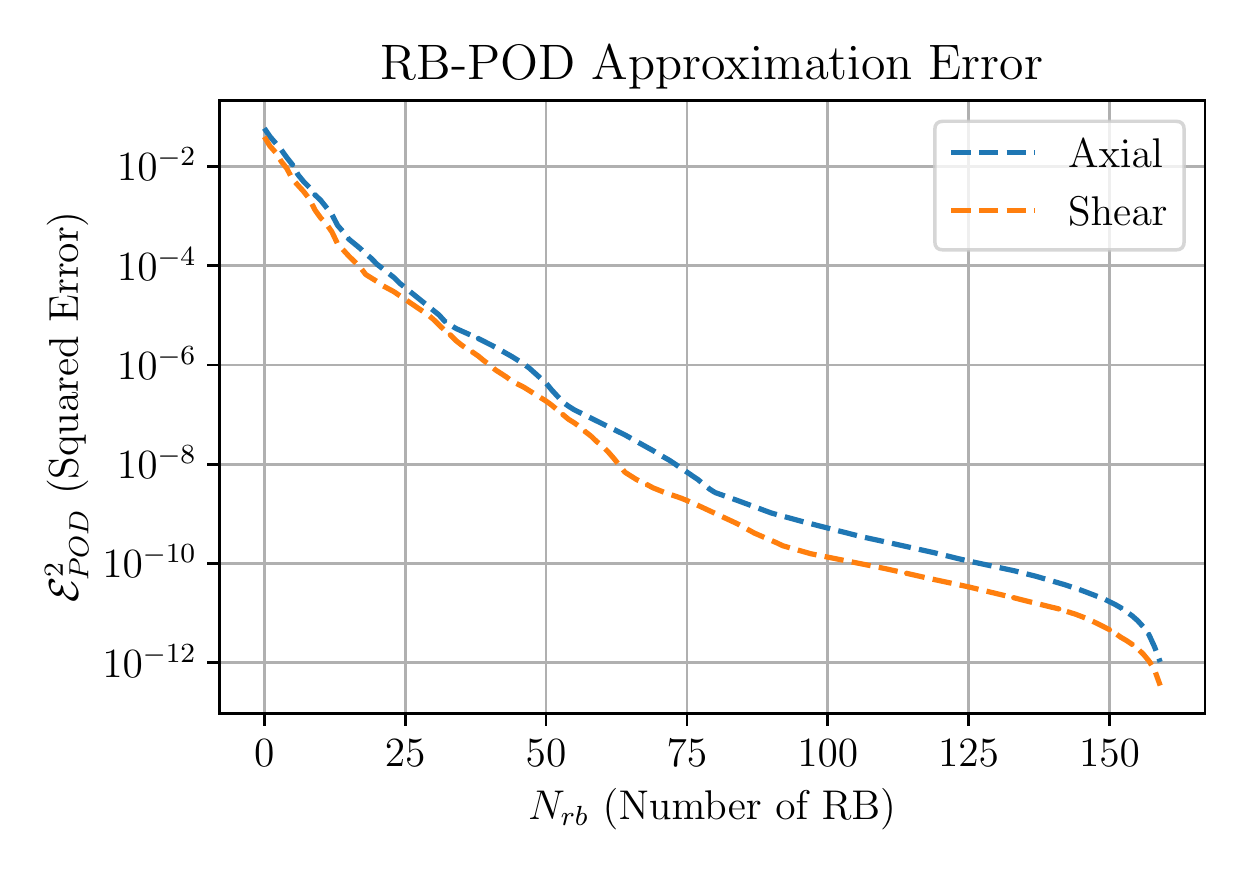}
				\caption{Error committed in the RB-ROM approximation.}
				\label{fig:poderror}
			\end{figure}
			
			To define the models $\mathcal{N}^{(1)}$ and $\mathcal{N}^{(3)}$, it is worth studying the sensibility of the total error $\mathcal{E}_T$ with the number of RB (see \eqref{eq:totalError}). It is reasonable to expect that at some point the decreasing of $\mathcal{E}_{POD}$ may not compensate the increasing of $\mathcal{E}_{DNN}$ as $N'_{rb}$ gets larger. Therefore, we have considered $N_{rb}' \in \{5,10,20,40,80,140\}$ and have trained $6$ different DNN architectures (one for each $N_{rb}'$), which do not differ regarding the hidden layers, but only the number of entries on the output. Also, we have considered a DNN architecture of $3$ hidden layers with $300$ neurons each. As activation functions, we have used the Swish \citep{swish} in the hidden layers followed by a linear activation in the output layer. To regularise the model, we have used the dropout with probability of retention $p = 99,5\%$ in all hidden layers and $\lambda = 10^{-8}$ for $l^2$ regularisation for all weights and biases. The DNN has been trained using the ADAM optimiser with default parameters and mini-batches of $32$ samples. An adaptive learning rate has been employed, decreasing linearly from $5.0 \times 10^{-4}$ until $5.0 \times 10^{-5}$ in epoch $5000$ (last). It is worth mentioning that we reached these values by systematically testing combinations of smaller architectures, regularisation values, and activation functions. Therefore, such values should be understood as guidance for problems with similar complexity and dataset sizes, although there might be still room for small tunings. Finally, we used \texttt{Tensorflow} \citep{tensorflow} for the implementation, and the final model is assumed to be the one yielding to the least loss function value in the validation dataset. 
			
			In Figures \ref{fig:trainingAxial_ny5} and \ref{fig:trainingAxial_ny140} we can see the loss function history ($\mathcal{E}_{DNN}^2$) for the axial DNN model using $N'_{rb} = 5$ and $N'_{rb} = 140$, respectively. Loss function values are naturally larger for $N'_{rb} = 140$, since the output has more entries. The loss function evaluated on training and validation differ considerably as result of regularisation, which is not used in validation, the latter turning out to be smaller. For the last epochs, the validation reaches a plateau while the training error continues to decrease slightly. This indicates that the model reached its full capacity while avoiding overfitting. The very same comments apply to Figures \ref{fig:trainingShear_ny5} and \ref{fig:trainingShear_ny140}, for the shear DNN model. For the results with different $N'_{rb}$, check Appendix \ref{sec:additionalPlots}.
			
			\begin{figure}
				\centering
				\subfigure[For $N'_{rb} = 5$.]{
					\includegraphics[width=0.48\linewidth]{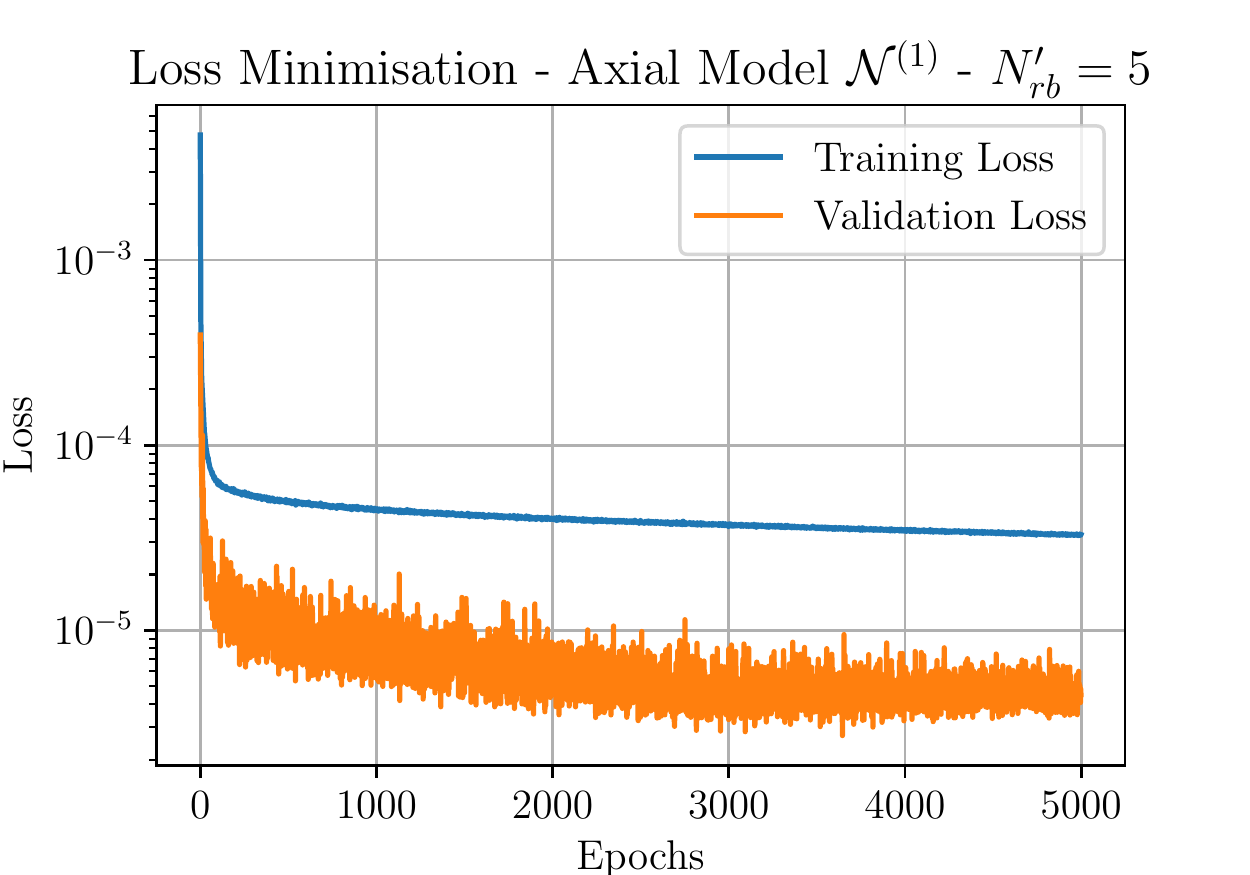}
					\label{fig:trainingAxial_ny5}
				}
				\subfigure[For $N'_{rb} = 140$.]{
					\includegraphics[width=0.48\textwidth]{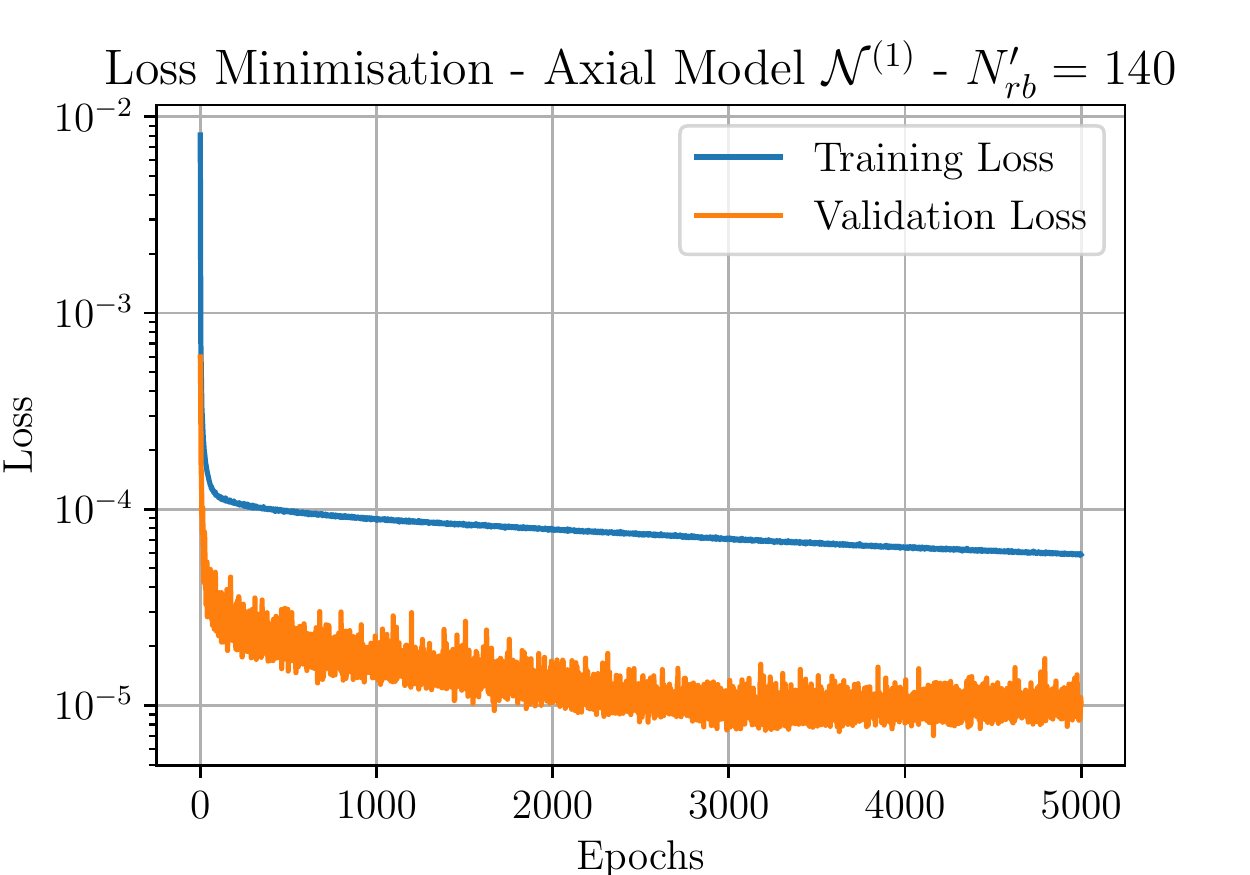}
					\label{fig:trainingAxial_ny140}
				}
				\caption{{\added Historic for the loss function minimisation in Axial model $\mathcal{N}^{(1)}$. Training loss evaluated with regularisation (dropout and $l^2$), while validation loss in prediction mode (regularisation disabled).}}
			\end{figure}
			
			\begin{figure}
				\centering
				\subfigure[For $N'_{rb}=5$.]{
					\includegraphics[width=0.48\linewidth]{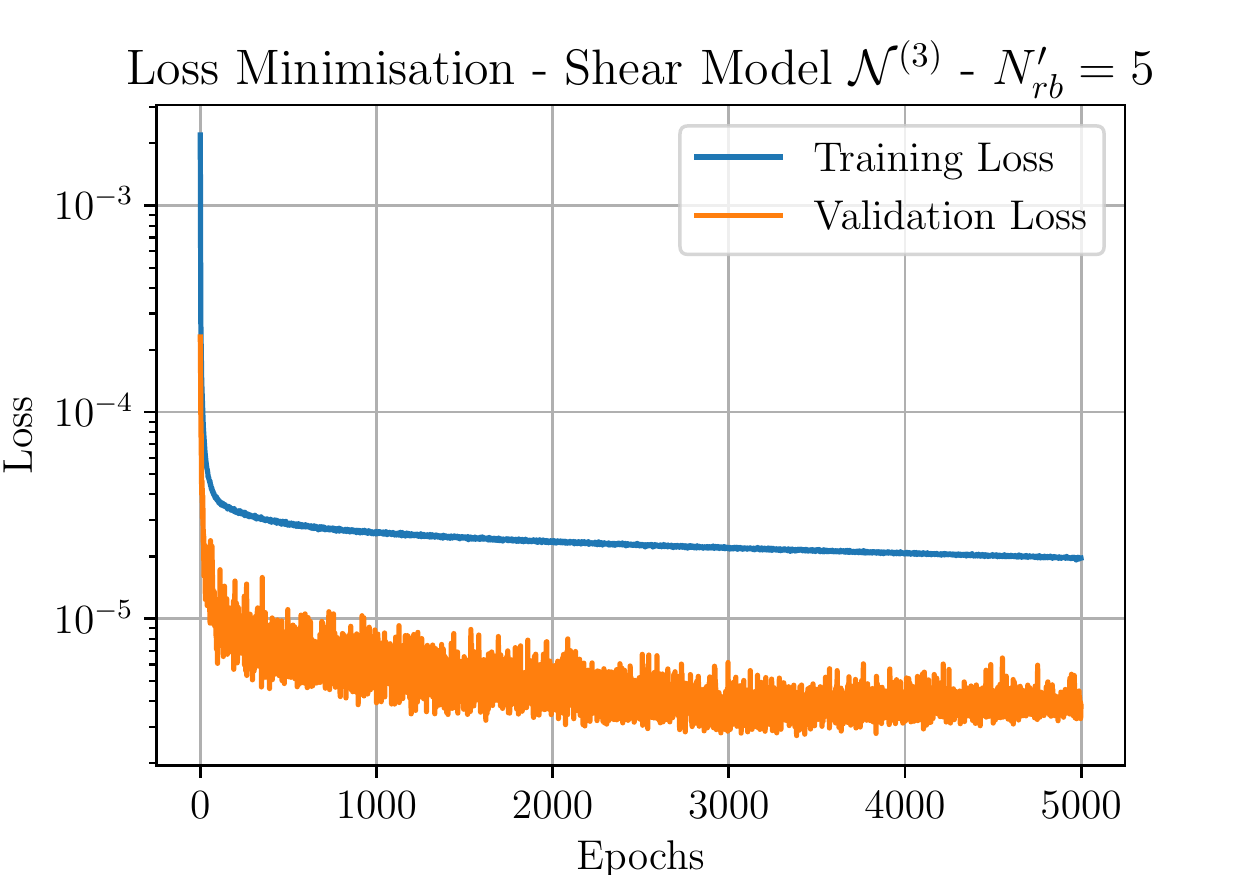}
					\label{fig:trainingShear_ny5}
				}
				\subfigure[For $N'_{rb}=140$.]{
					\includegraphics[width=0.48\textwidth]{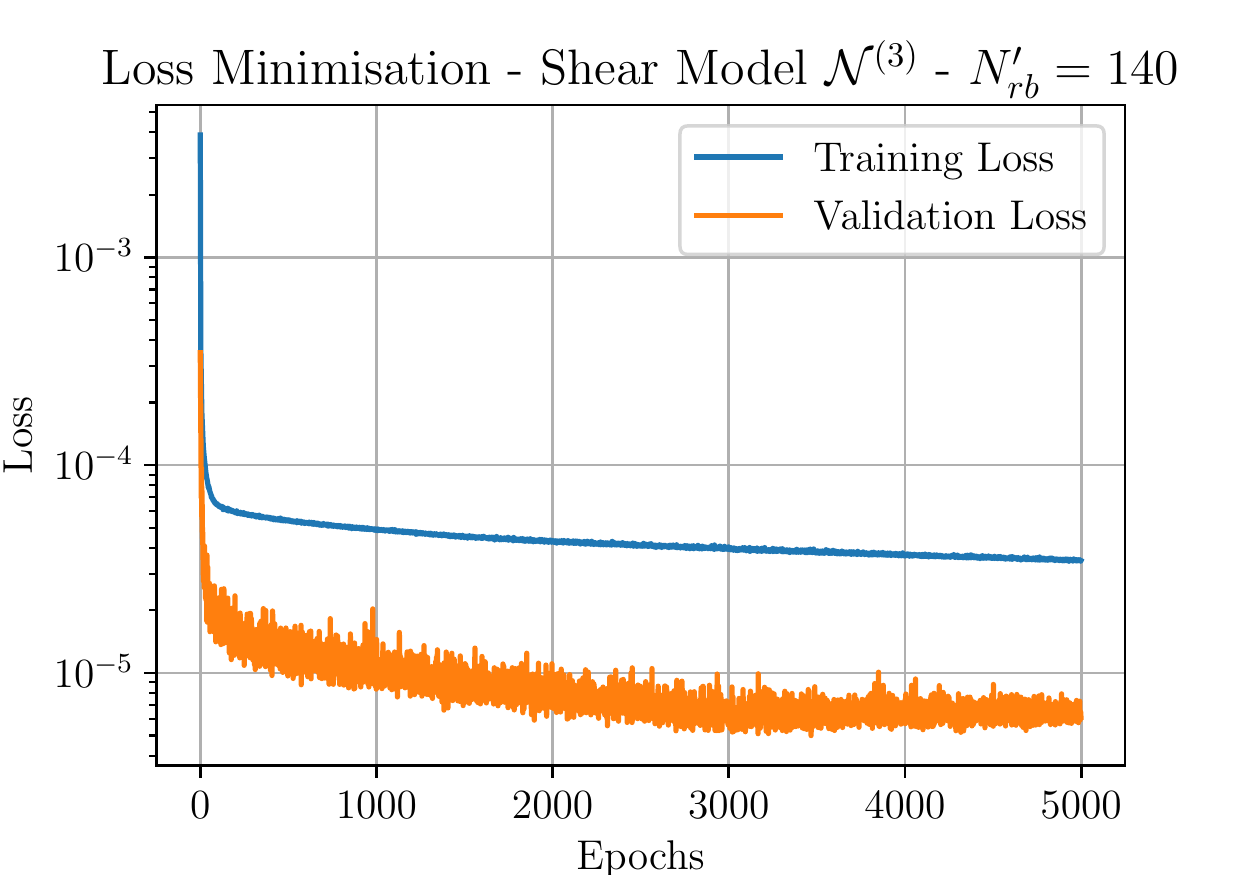}
					\label{fig:trainingShear_ny140}
				}
				\caption{{\added Historic for the loss function minimisation in Shear model $\mathcal{N}^{(3)}$. Training loss evaluated with regularisation (dropout and $l^2$), while validation loss in prediction mode (regularisation disabled).}}
			\end{figure}

			Lastly, the final trained DNN models are compared in terms of the DNN, POD, and total errors (in test datasets) in Figures \ref{fig:testAxial_error} and \ref{fig:testShear_error} for the axial and shear cases, respectively. For the axial case we observe the least total error for $N'_{rb}=140$ while in the shear case both $N'_{rb}=80$ and $N'_{rb}=80$ give equivalent quantities. Therefore, we choose to keep $N'_{rb}=140$ to the numerical examples of next section.

			\begin{figure}[h!]
				\centering
				\subfigure[Axial model ($\mathcal{N}^{(1)}$)]{
					\includegraphics[width=0.48\linewidth]{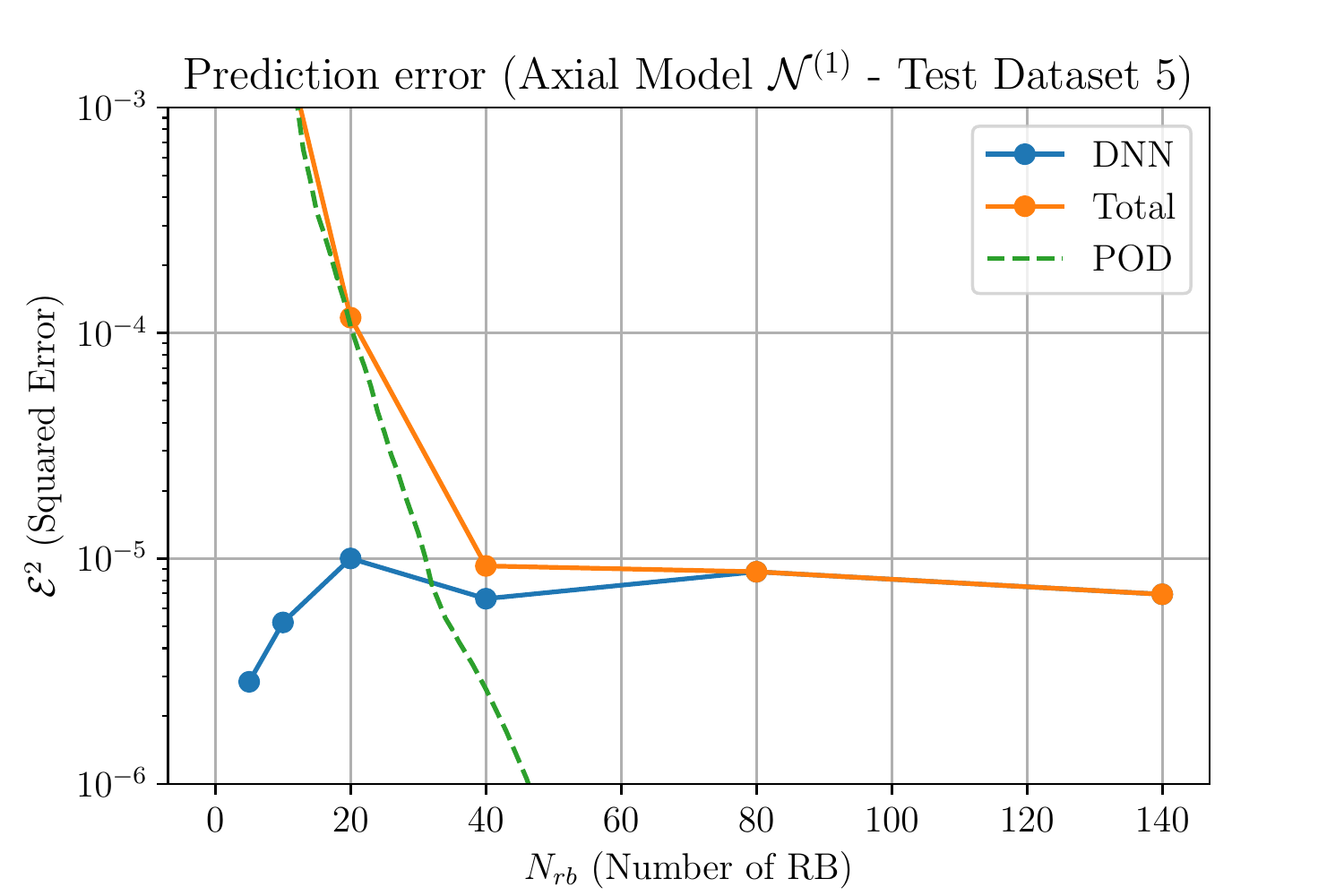}	
					\label{fig:testAxial_error}
				}
				\subfigure[Shear model ($\mathcal{N}^{(3)}$)]{
					\includegraphics[width=0.48\textwidth]{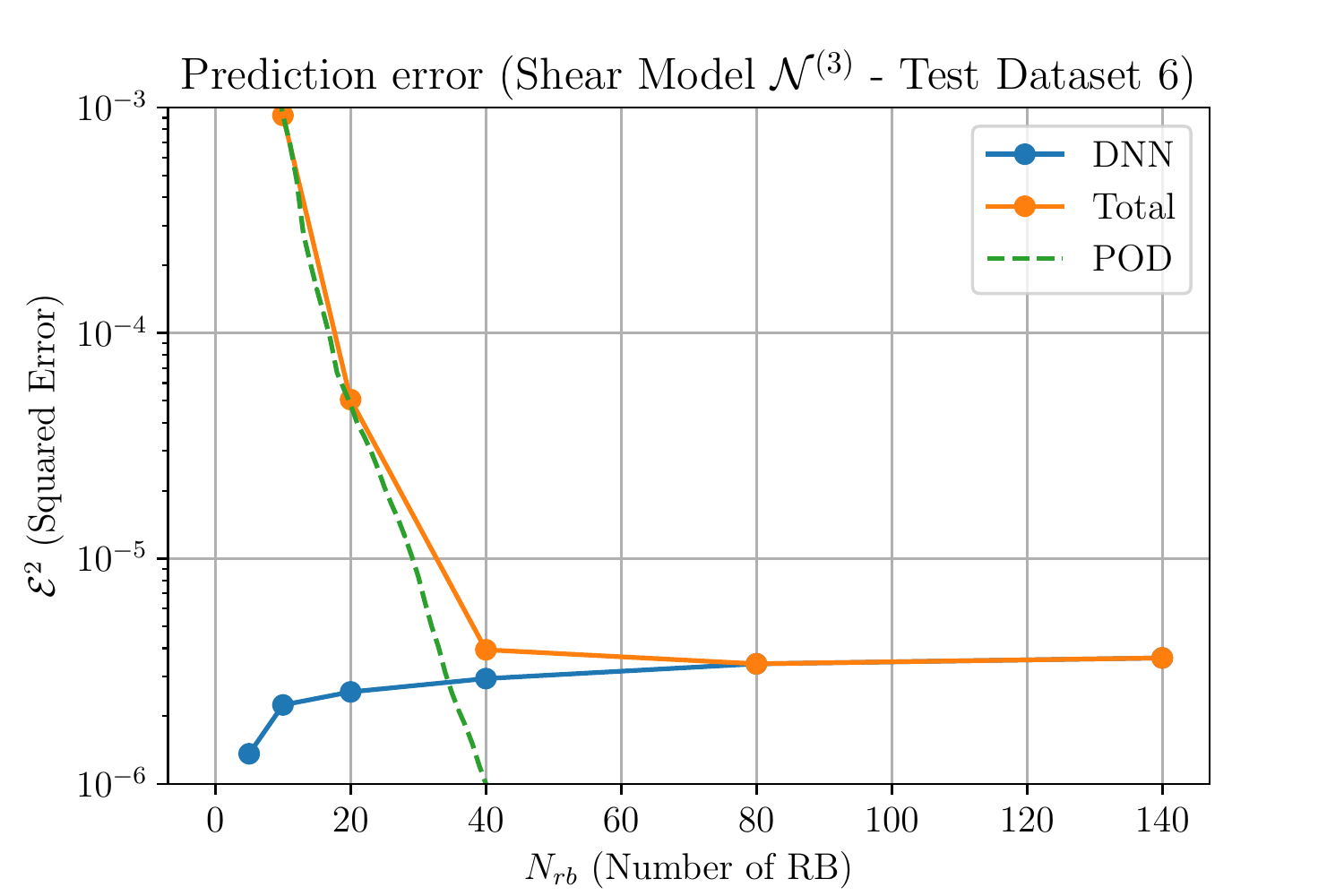}
					\label{fig:testShear_error}
				}
				\caption{{\added Mean Squared Prediction error (test datasets) for the $L^2(\partial \Omega_{\mu})$: DNN and POD contributions yields to the Total error.}}
			\end{figure}

		}

		\FloatBarrier

		\section{Numerical Examples}  \label{sec:NumericalExamples} 
		
		In this section we aim at assessing the \DeepBC method in practice, i.e., in coupled two-scale simulations, also known as \fe2 \citep{Feyel2003}. Two kinds of problems are considered ({\added recall} the notation of Figure \ref{fig:randommedia}): i) the first with a strong scale separation, i.e., the size of heterogeneities is several orders of magnitude smaller than the dimension of the body such that $L_{\mu}\ll L$; and ii) the second with $L_{\mu}< L$ with about one or two orders of magnitude of difference. As already commented, the computational homogenisation approach is best suited for the i)-kind-problem, however due practical reasons ii)-kind-problem is also interesting since only in this situation direct numerical simulations (DNS) are feasible. The two kinds of problems differ substantially in terms of how the microstructures are generated at integration points. We postpone the presentation of further details to Section \ref{sec:cookMembrane} and Section \ref{sec:DNS}.
		
		Before digging into the numerical results, it is worth mentioning that we aim at answering the following questions:
		\begin{enumerate}
			\item How does \DeepBC compare in terms of accuracy with the periodic condition using the high-fidelity solution as reference?
			\item How does \DeepBC compare in terms of computational cost the high-fidelity model?
			\item How does \DeepBC compare in terms of accuracy with the periodic and high-fidelity approaches using the DNS as the reference solution? 
		\end{enumerate}
		
		Here, it is worth recalling that the \DeepBC and periodic models refer to local problems posed in $\Omega_{\mu}^R$, thus with the BC imposed on $\partial \Omega_{\mu}^R$, while high-fidelity stands for local problems $\Omega_{\mu}^H$, with periodic conditions imposed on $\partial \Omega_{\mu}^H$. We decided to choose periodic BCs for comparisons, since it is widely accepted that this model usually leads to more accurate results in contrast with classical models for BCs.  As a natural choice, we answer questions 1 and 2 above using a i)-kind-problem (Section \ref{sec:cookMembrane}), and question 3 using a ii)-kind-problem (Section \ref{sec:DNS}). 
		
		Concerning the computational framework, in all numerical examples we use the \textit{micmacFenics} library \citep{micmacsFenics}, an open-source Fenics-based library for the implementation classical BCs in the \fe2 setting, which can be alternatively found as part of \textit{multiphenics} \citep{multiphenics}.
		
		\FloatBarrier
		
		\subsection{Cook Membrane} \label{sec:cookMembrane}
		
		\begin{figure}
			\centering
			\subfigure[Schema showing dimensions and points of interest.]{
				\includegraphics[width=0.35\linewidth]{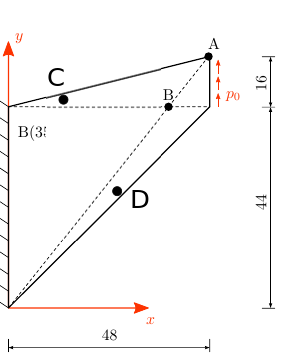}
				\label{fig:cook_schema}
			}
			\subfigure[Coaser mesh considered, containing $5$ elements on the vertical direction.]{
				\includegraphics[width=0.3\textwidth]{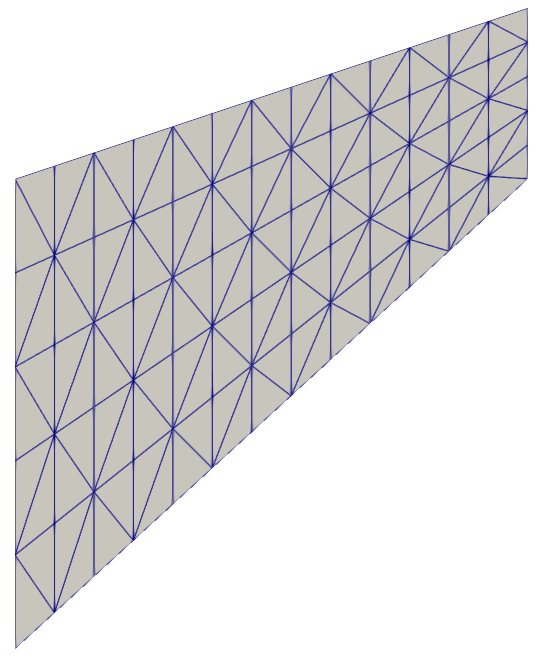}
				\label{fig:cook_mesh}
			}
			\caption{Cook membrane.}
		\end{figure}

		In this example for i)-kind-problem we apply the \DeepBC to the so-called Cook Membrane, a classical benchmark in solid mechanics normally used to test numerical methods for incompressible material \citep{Schroder2021}, but also used in computational homogenisation \citep{Steinmann2016} and to test ML-based surrogate models \citep{Rabczuk2020}. It consists in a solid structure with dimensions as in Figure \ref{fig:cook_schema} subjected to an homogeneous distributed vertical load $\Vec{t}$ on the rightmost face. The magnitude of $\Vec{t}$ and the material model itself varies among the literature considered. Therefore, for the sake of simplicity, we decided to consider the very same material as used in the dataset construction of Section \ref{sec:dataset} and $\Vec{t} = (0,0.05)$, which gives approximately the same vertical displacement of the tip $A$ found in \citep{Schroder2021}.
		
		We aim at comparing the results obtained using the boundary conditions i) \DeepBC and ii) Periodic against our reference solution, obtained from high-fidelity microstructures. For the sake of simplicity, we choose the MDs randomly from the \textit{validation} and \textit{test} datasets from Section \ref{sec:dataset}, already at hand. They are sorted among the elements without repetition, which is possible since there are fewer integration points in the finer mesh than the pre-simulated structures microstructures available. To obtain more reliable results, we always report results considering the average over 10 realisations, where each realisation should be understood as a different sorting {\added of MDs along the integration points}. Importantly, also notice that due the scale separation, the MDs choices keeps no neighbouring relation.
		
		As first study, we consider the convergence of the vertical displacement of the tip A (Figure \ref{fig:cook_schema}) with respect to the mesh divisions in the vertical direction. The coarser mesh looks like the one depicted in Figure \ref{fig:cook_mesh}, with $5$ regular divisions along the vertical direction and a the horizontal number of divisions such as we keep the same characteristic length of the rightmost side. We then build the other finer meshes with $10, 20, 40,$ and $80$ divisions. We can see the results in Figure \ref{fig:dispA}, where a clear advantage of \DeepBC against the periodic can be observed when compared with the high-fidelity results. In numerical terms, for the finest mesh we gain between one and two orders of magnitude in the relative error norm (see Table \ref{tab:relativeErrorsFinerMesh}). 
		
		\begin{figure}
			\centering
			\includegraphics[width=0.5\linewidth]{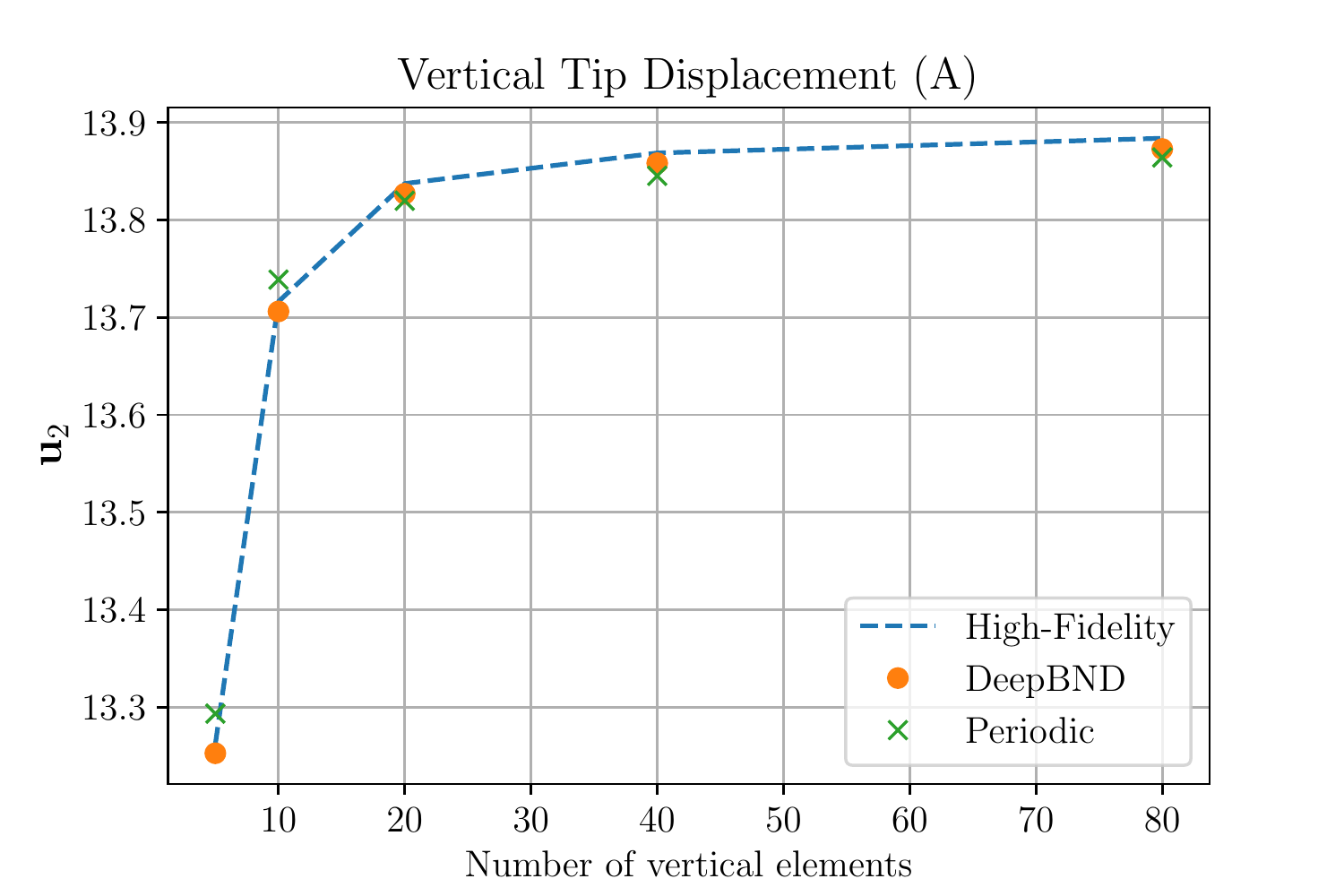}
			\caption{Convergence of the vertical displacement at the tip (A) in averaged values for 10 realisations of the MDs sorting.}
			\label{fig:dispA}
		\end{figure}

		\begin{table}
			\begin{center}
				\begin{tabular}{|c|c|c|c|}
					\hline 
					& \DeepBC & Periodic \\ 
					\hline 
					Vertical Displacement A & $7.923\times 10^{-4} \pm 1.091\times 10^{-4}$ & $1.417\times 10^{-3} \pm 5.496\times 10^{-4} $ \\ \hline
					Von Mises B & $2.185\times 10^{-4} \pm 1.758\times 10^{-4}$ & $3.126\times 10^{-2} \pm 1.722\times 10^{-2}$ \\ \hline
					Von Mises C & $2.444\times 10^{-4} \pm 1.799\times 10^{-4}$ & $1.096\times 10^{-2} \pm 1.060\times 10^{-2} $	 \\ \hline
					Von Mises D  & $1.288\times 10^{-4} \pm 9.262\times 10^{-5}$ & $6.096\times 10^{-3} \pm 4.362\times 10^{-3}$ \\ \hline 
				\end{tabular} 
			\end{center}
			\caption{Relative errors in the finest mesh.}
			\label{tab:relativeErrorsFinerMesh}
		\end{table}
		
		\begin{table}
			\begin{center}
				\begin{tabular}{|c|c|c|c|}
					\hline 
					& \DeepBC & {\added HF MD (Linear BC)}\\ 
					\hline 
					Time per snapshot(s) & $0.5148 \pm 0.09548$ & $30.56 \pm 7.238 $ \\ \hline
					Total Time Finer Mesh (h) & $1.831$ & $108.6$ \\ \hline
					Wall Time Finer Mesh (h) ($32$ Cores) & $0.05722$ & $3.394$	 \\
					\hline
					Speed-up (ratio time per snapshot) & $59.36 \times$ &  $1\times$ \\ \hline 
				\end{tabular} 
			\end{center}
			\caption{Computational cost for the simulations using \DeepBC and High-Fidelity cases. To a fairer comparison in terms of computational times, note that the values reported for the High-Fidelity corresponds to the case that is implemented via homogeneous Dirichlet BCs {\added (Linear BC space)}, thus differing from those from Table \ref{tab:computationalTimeDataset}, implemented with periodic BCs.}
			\label{tab:computationalTimeReduced}
		\end{table}
		
		As a second comparison, Figure \ref{fig:vonMisesStress} shows the von Mises equivalent stress field for the different models, where we can notice that the results for the periodic case (rightmost plot) are considerably noisier, which means worse, compared to the \DeepBC simulation (central plot) that visually presents identical results to High-fidelity solution (leftmost plot). Such discrepancies become clearer by plotting errors of the Von Mises stresses, as depicted Figure \ref{fig:vonMisesStressError}, in which \DeepBC performs about 3 orders of magnitude better than the periodic case. Similar conclusions for the stress are also retrieved component-wise in Figure \ref{fig:stressCook}. 
		Analysing Figure \ref{fig:vonMisesStress}, \ref{fig:stressCook}, and \ref{fig:cook_schema}, we can see that points B, C and D were placed where extreme values for stresses are retrieved. Comparing the Von Mises stress along all realisations, \DeepBC gains between one and two orders of magnitude in accuracy with respect to the periodic case, depending on the point considered, as seen in Table \ref{tab:relativeErrorsFinerMesh}.
		
		\begin{figure}
			\centering
			\subfigure[$\sigma^{HF}_{vm}$.]{
				\includegraphics[width=0.26\textwidth]{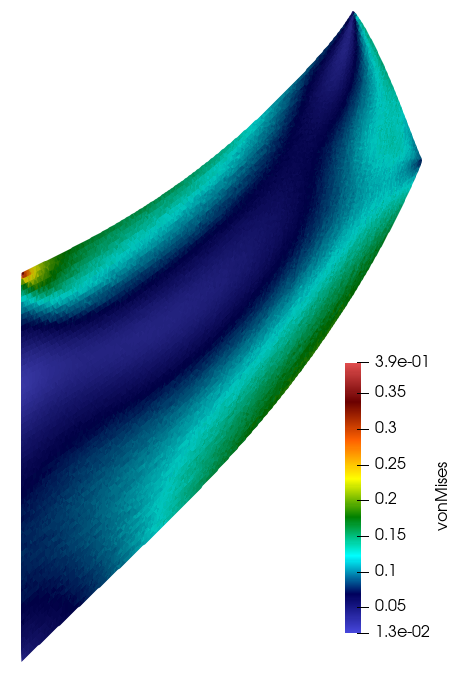}
			}
			\subfigure[$\sigma^{dnn}_{vm}$.]{
				\includegraphics[width=0.26\textwidth]{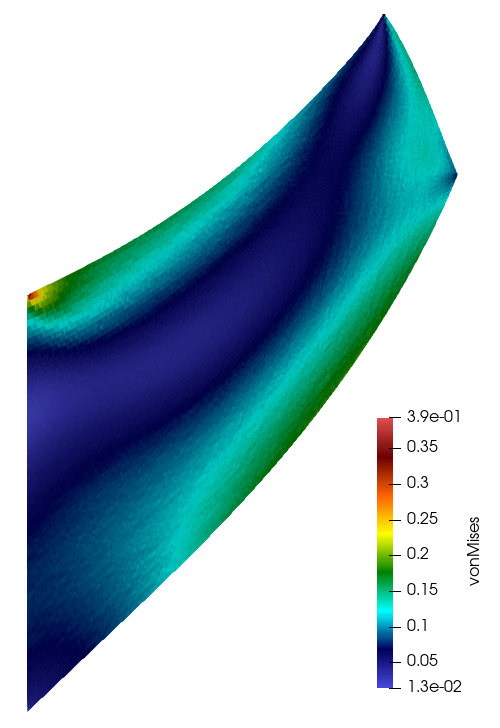}
			}
			\subfigure[$\sigma^{per}_{vm}$.]{
				\includegraphics[width=0.26\textwidth]{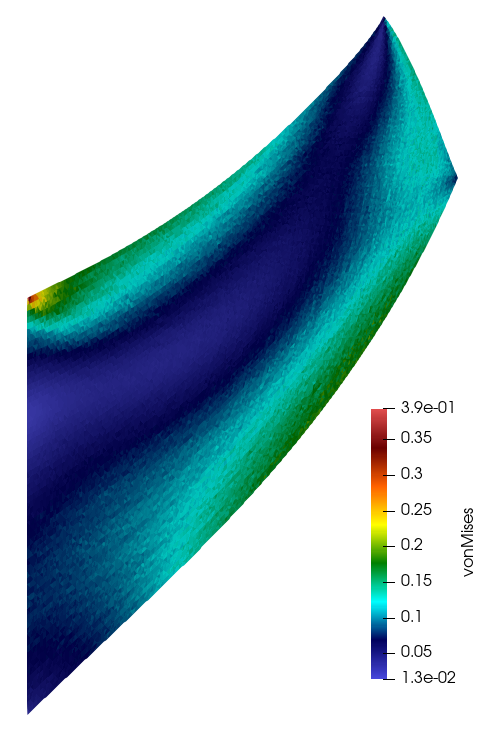}
			}
			\caption{Comparison for the Von Mises Stress ($\sigma_{vm}$) in the Cook membrane.}
			\label{fig:vonMisesStress}
		\end{figure}

		\begin{figure}
			\centering
			\subfigure[$|\sigma^{dnn}_{vm} - \sigma^{ref}_{vm}|$.]{
				\includegraphics[width=0.26\textwidth]{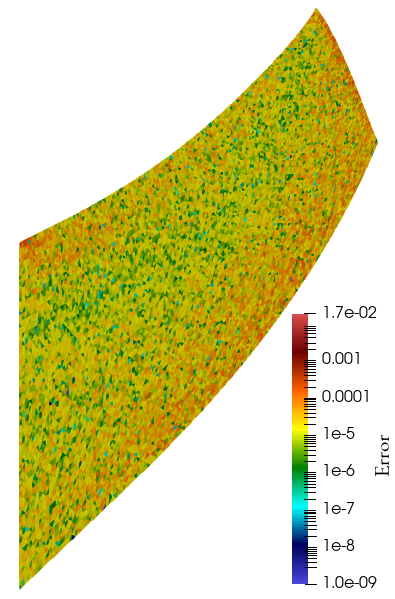}
			}
			\subfigure[$|\sigma^{per}_{vm} - \sigma^{ref}_{vm}|$.]{
				\includegraphics[width=0.26\textwidth]{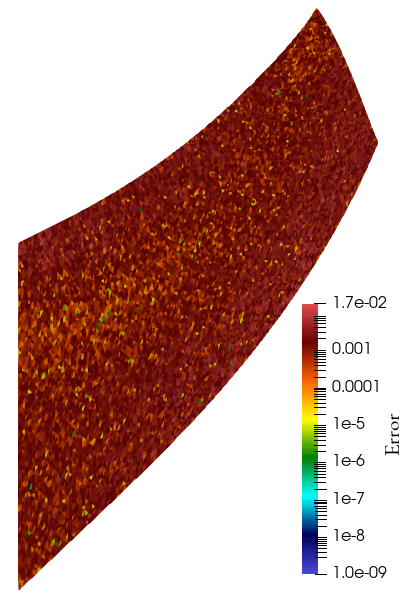}
			}
			
			\subfigure[$\frac{|\sigma^{dnn}_{vm} - \sigma^{ref}_{vm}|}{\sigma^{ref}_{vm}}$.]{
				\includegraphics[width=0.26\textwidth]{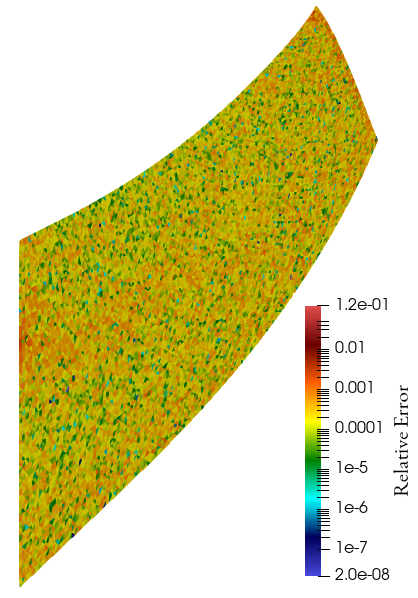}
			}
			\subfigure[$\frac{|\sigma^{per}_{vm} - \sigma^{ref}_{vm}|}{\sigma^{ref}_{vm}}$.]{
				\includegraphics[width=0.26\textwidth]{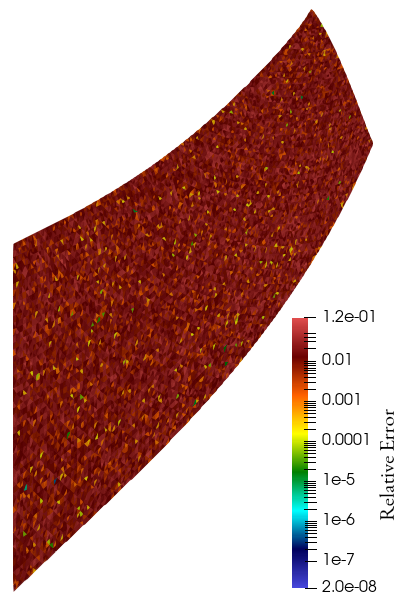}
			}
			\caption{Comparing \DeepBC and periodic simulations in terms of absolute and relative errors of the Von Mises Stress ($\sigma_{vm}$) for the Cook membrane.}
			\label{fig:vonMisesStressError}
		\end{figure}
		
	}

	In terms of computational times, we can see in Table \ref{tab:computationalTimeReduced} that the computational speed-up reached by \DeepBC {\added with respect to HF with linear BCs} was about $59\times$. This number refers only to the embarrassingly parallel part of the problem, which in this case is the assembly of FEM matrices for the macroscale problem, specifically the solution of the local problems. The remaining fixed time for solving the macroscale problem is the same, hence no reason of including it here. Also, the mesh generation is not taken into account, but if taken, the speed-up would be even higher. The \DeepBC was not compared directly with periodic BCs times, since due to the more efficient native Fenics/Multiphenics implementation of enforcing Dirichlet BCs {\added with \DeepBCn}, instead of the periodicity which is imposed by Lagrange Multipliers with additional pure Python routines.

	\begin{figure}
		\centering
		\subfigure[$\sigma^{HF}_{11}$.]{
			\includegraphics[width=0.26\textwidth]{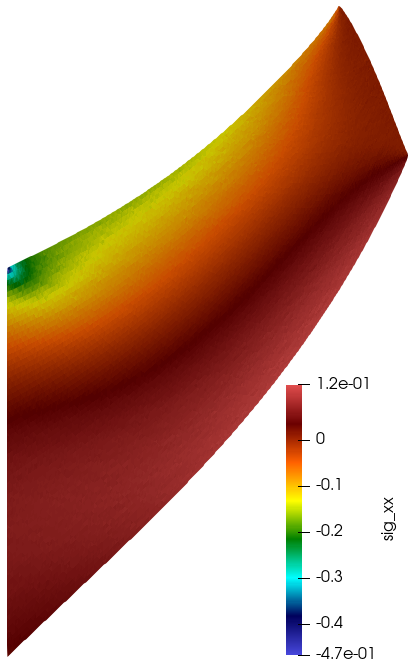}
		}
		\subfigure[$\sigma^{dnn}_{11}$.]{
			\includegraphics[width=0.26\textwidth]{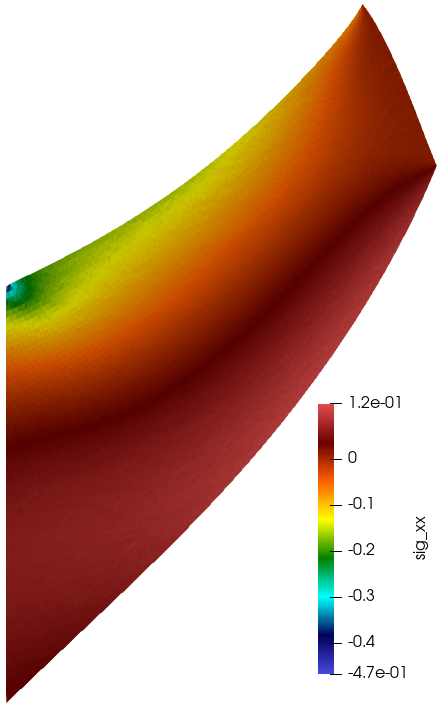}
		}
		\subfigure[$\sigma^{per}_{11}$.]{
			\includegraphics[width=0.26\textwidth]{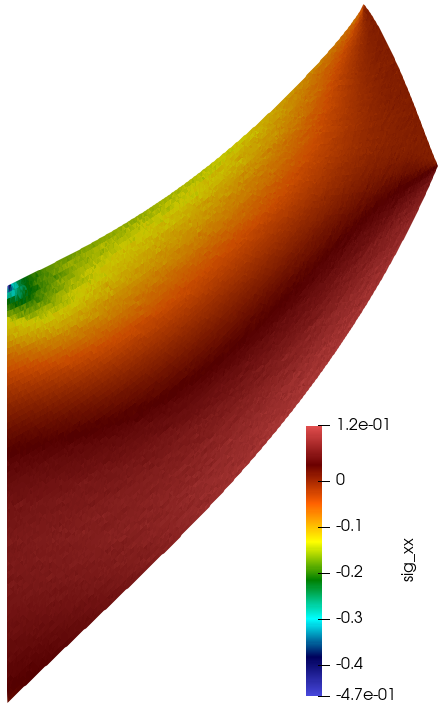}
		}
		\subfigure[$\sigma^{HF}_{22}$.]{
			\includegraphics[width=0.26\textwidth]{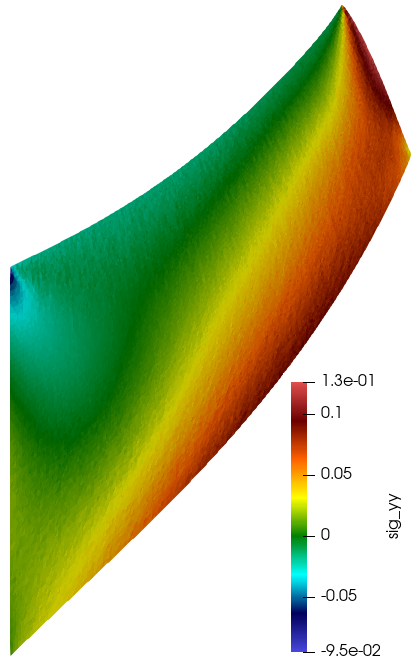}
		}
		\subfigure[$\sigma^{dnn}_{22}$.]{
			\includegraphics[width=0.26\textwidth]{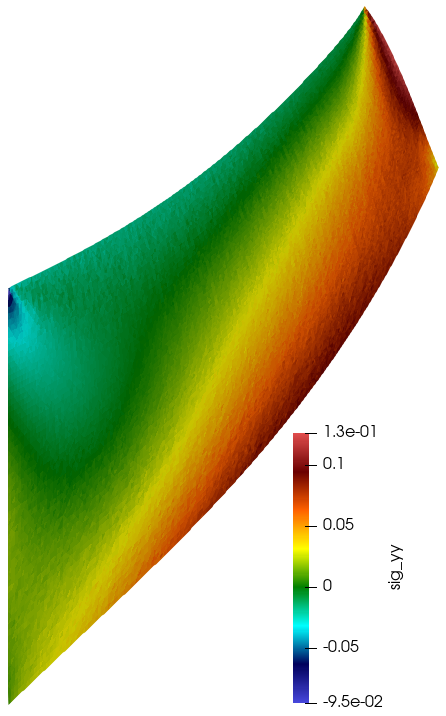}
		}
		\subfigure[$\sigma^{per}_{22}$.]{
			\includegraphics[width=0.26\textwidth]{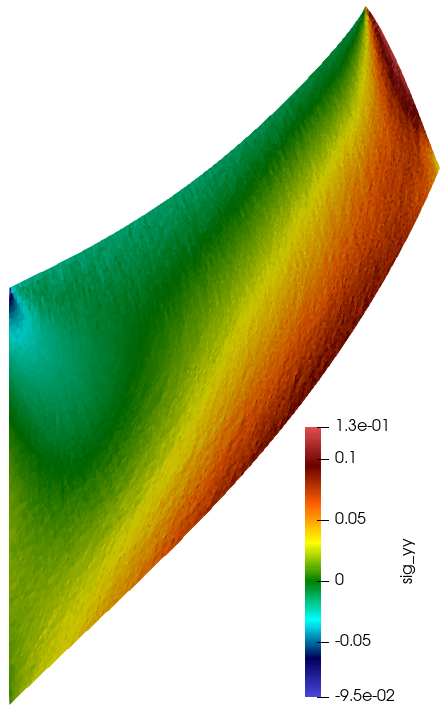}
		}
		\subfigure[$\sigma^{HF}_{12}$.]{
			\includegraphics[width=0.26\textwidth]{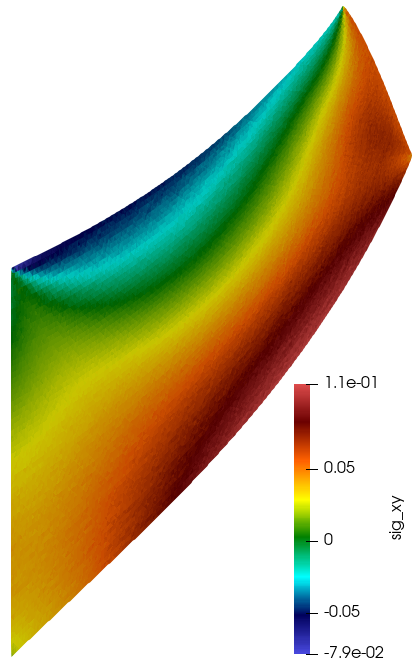}
		}
		\subfigure[$\sigma^{dnn}_{12}$.]{
			\includegraphics[width=0.26\textwidth]{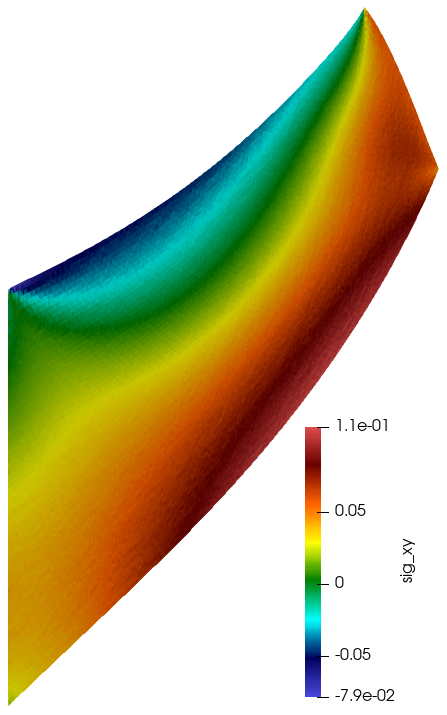}
		}
		\subfigure[$\sigma^{per}_{12}$.]{
			\includegraphics[width=0.26\textwidth]{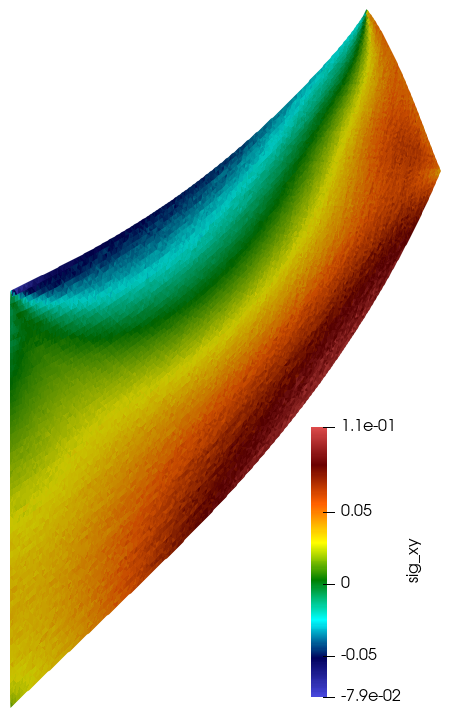}
		}
		\caption{Comparison Cauchy stress components in the Cook membrane.}
		\label{fig:stressCook}
	\end{figure}

	\FloatBarrier

	%
	%


	\subsection{Clamped bar using DNS}\label{sec:DNS}
	The last test of the new multiscale methodology aims at assessing its performance in terms of accuracy and computational cost {\added with respect to} other multiscale methods in comparison to direct numerical simulations (DNS). In this second numerical example, for easiness of microstructures constructions, we adopt a rectangular bar $\Omega = [0,L_x = 4]\times[0, L_y = 1]$ as macroscopic domain such that it is clamped on the left and loaded with a uniform shear on the right $\Vec{t} = (0,-0.2)$, 
	as in Figure \ref{fig:barClamped}. The material is organised by a $N_x^{DNS} \times N_y^{DNS}$ grid of $H^{DNS} = L_x/N_x^{DNS} = L_y/N_y^{DNS}$ squared blocks in which each block contains a random sized circular inclusion.
	The mesh for the DNS simulations is generated with the characteristic length $h \approx H^{DNS}/15$, which is similar to the one used in the dataset generation. 
	We consider $N_y^{DNS} \in \{24,72\}$, to study two different approximations with respect to the DNS, which is expected to be more suitable for larger $N_y^{DNS}$. Just for visualisation purposes, we show in Figure \ref{fig:dispDNS} one realisation of a structure with $N_y^{DNS}=7$. Moreover, it is worth remarking that the mesh generation for larger $N_y^{DNS}$ becomes cumbersome and computationally costly.
	Finally, it is also worth mentioning that the bar geometric dimensions were inspired by  \citep{Rabczuk2020}, however as in that reference the results presented were merely visual, we decided to compare \DeepBC against HF and periodic models only.
	
	Concerning the multiscale simulations, we use a regular triangular mesh with $100$ divisions vertically and $400$ horizontally. Notice that this mesh is several orders of magnitude coarser than the DNS, which makes the \fe2 simulations still appealing, despite its also elevated cost. As usual, one MD should be associated at each integration point. Differently from the previous case (Cook's Membrane, Section \ref{sec:cookMembrane}), now the microstructures have to keep some relation with its neighbourhood. Given an integration point, we choose the associated MD as the 6x6 window having the closest centre point, and then we proceed as before by identifying the reduced 2x2 window and predicting the BCs. This strategy is depicted in Figure \ref{fig:barClamped}, where the dashed windows are associated with the respective centre (circles) of the same colour. Notice that in total the number of different MDs are $(N_x^{DNS} - 5) \times (N_y^{DNS} - 5)$, which makes $1729$ and $18961$ for $N_y^{DNS}=24$ and $N_y^{DNS}=72$, respectively.
	
	\begin{figure}
		\centering
		\includegraphics[width=0.9\linewidth]{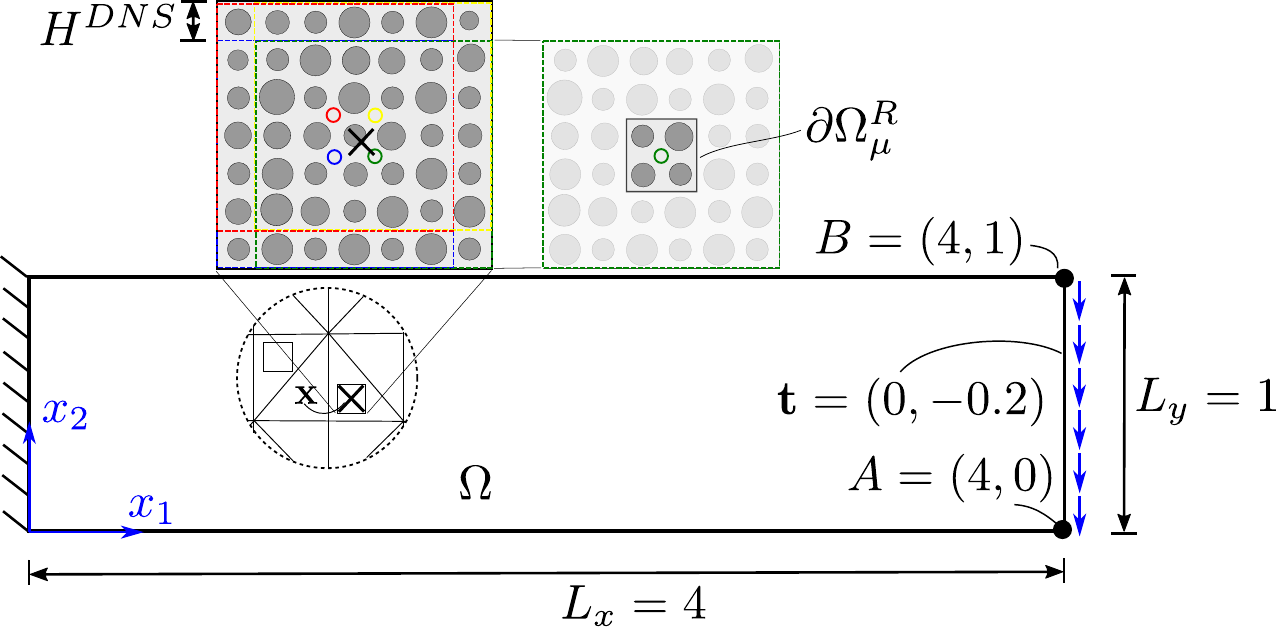}
		\caption{Clamped bar scheme showing geometric dimensions and loads. The integration point is represented by ($\times$) and underlying HF cell centres by circles. Colours indicate the correspondence between the centres and respective window.}
		\label{fig:barClamped}
	\end{figure}
	
	\begin{figure}
		\centering
		\includegraphics[width=0.9\linewidth]{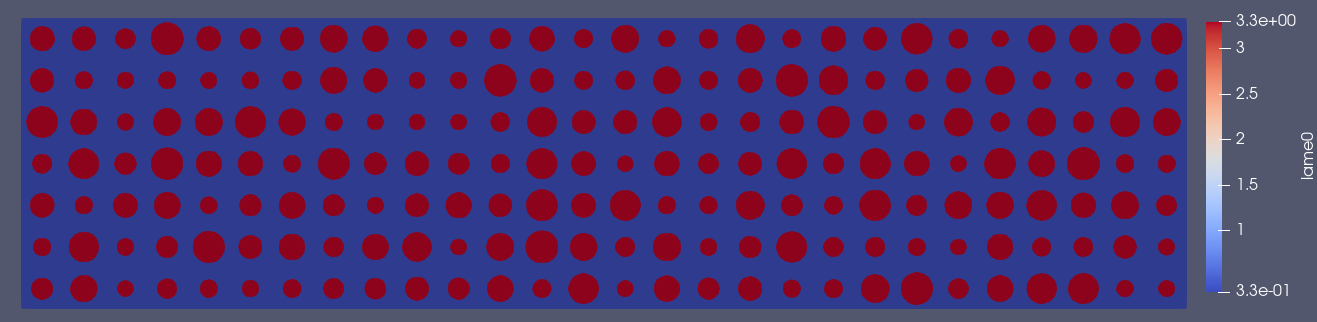}
		\caption{Realisation of bar with $N_y^{DNS} = 7$ and $N_x^{DNS} = 28$.}
		\label{fig:dispDNS}
	\end{figure}
	
	Concerning the results, each structure has been solved using quadratic polynomial elements for the DNS mesh (DNS solution) and using the multiscale approach, with the MDs selected as aforementioned, for the \DeepBCn, periodic, and the HF cases. The latter cases have been compared with the DNS solution with respect to the $L^2$ norm, and the Euclidean norm for displacements at points $A$ and $B$ (see Figure \ref{fig:barClamped}). We can see in Figure \ref{fig:error_rel_DNS_ny24} that relative errors for the $N_y^{DNS} = 24$ case are of order $1\%$, with slight advantages for the HF and \DeepBC contrasted with the periodic case. On the other hand, for the $N_y^{DNS} = 72$ case, \DeepBC model reaches $0.05\%$ in relative error for the $L^2(\Omega)$-norm while the periodic model attains $0.2\%$ in the same norm, i.e., a factor $4$ of improvement. Similar gains were found for the point-wise displacement norms. Unexpectedly, comparing with the HF case, \DeepBC performs even better, which shows that our approach does not only interpolate the HF targets, in the sense of averaging results taking close inputs, but it is also capable of learning more general relations associated to the whole assemble of data.    
	

	\begin{figure}
		\subfigure[]{
			\includegraphics[width=0.45\linewidth]{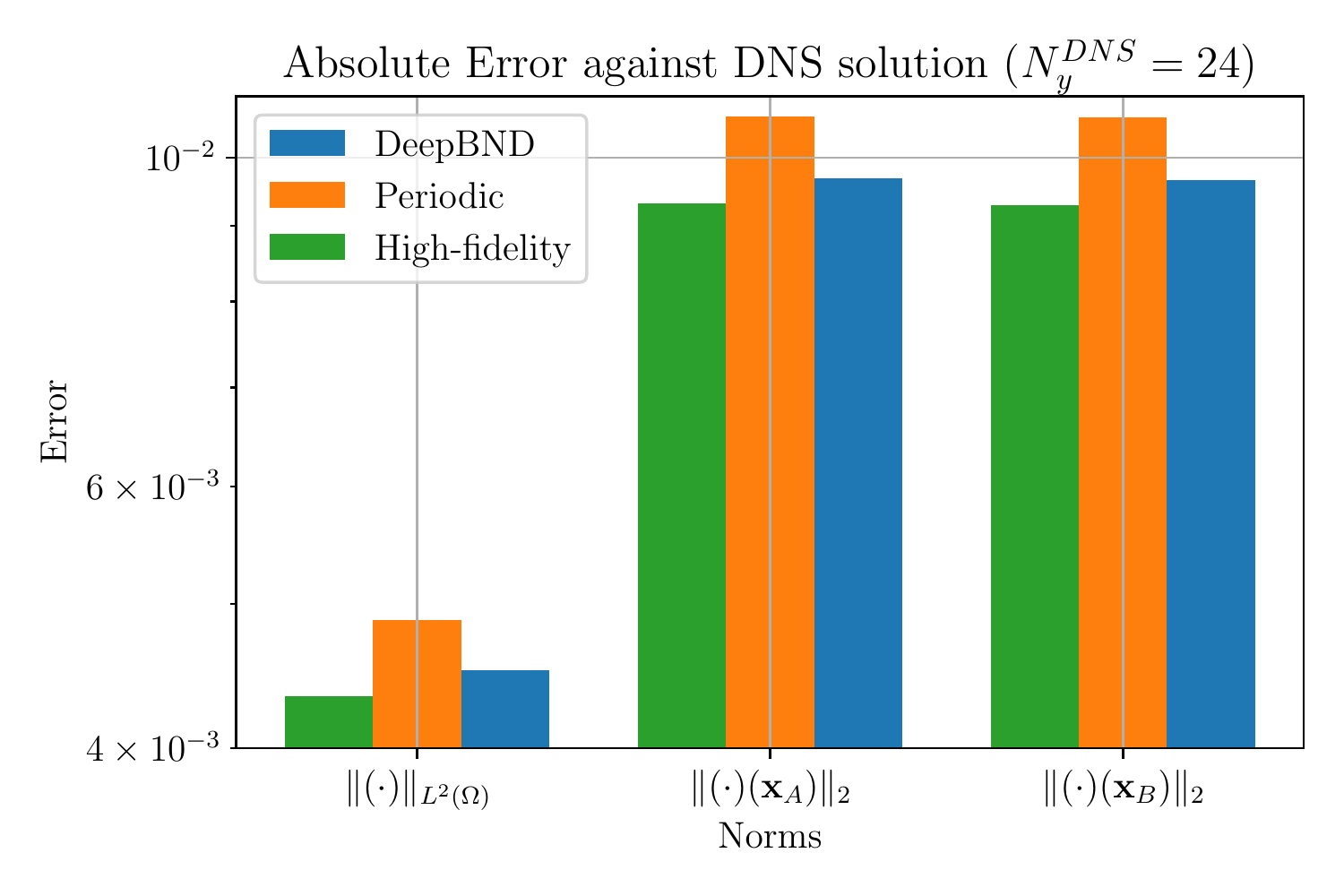}
			\label{fig:errorDNS_ny24}
		}
		\centering
		\subfigure[]{
			\includegraphics[width=0.45\linewidth]{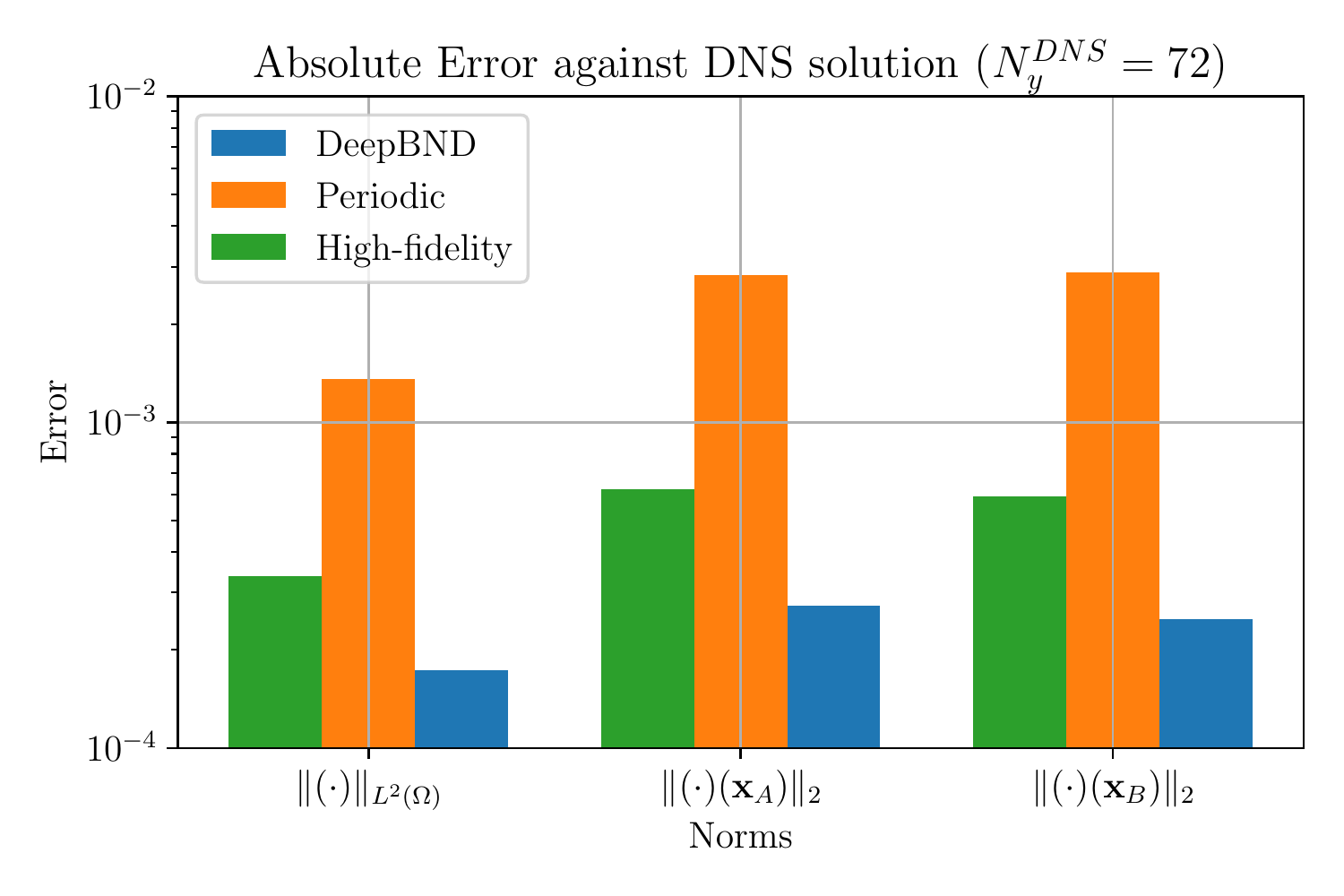}
			\label{fig:errorDNS_ny72}
		}
		\subfigure[]{
			\includegraphics[width=0.45\linewidth]{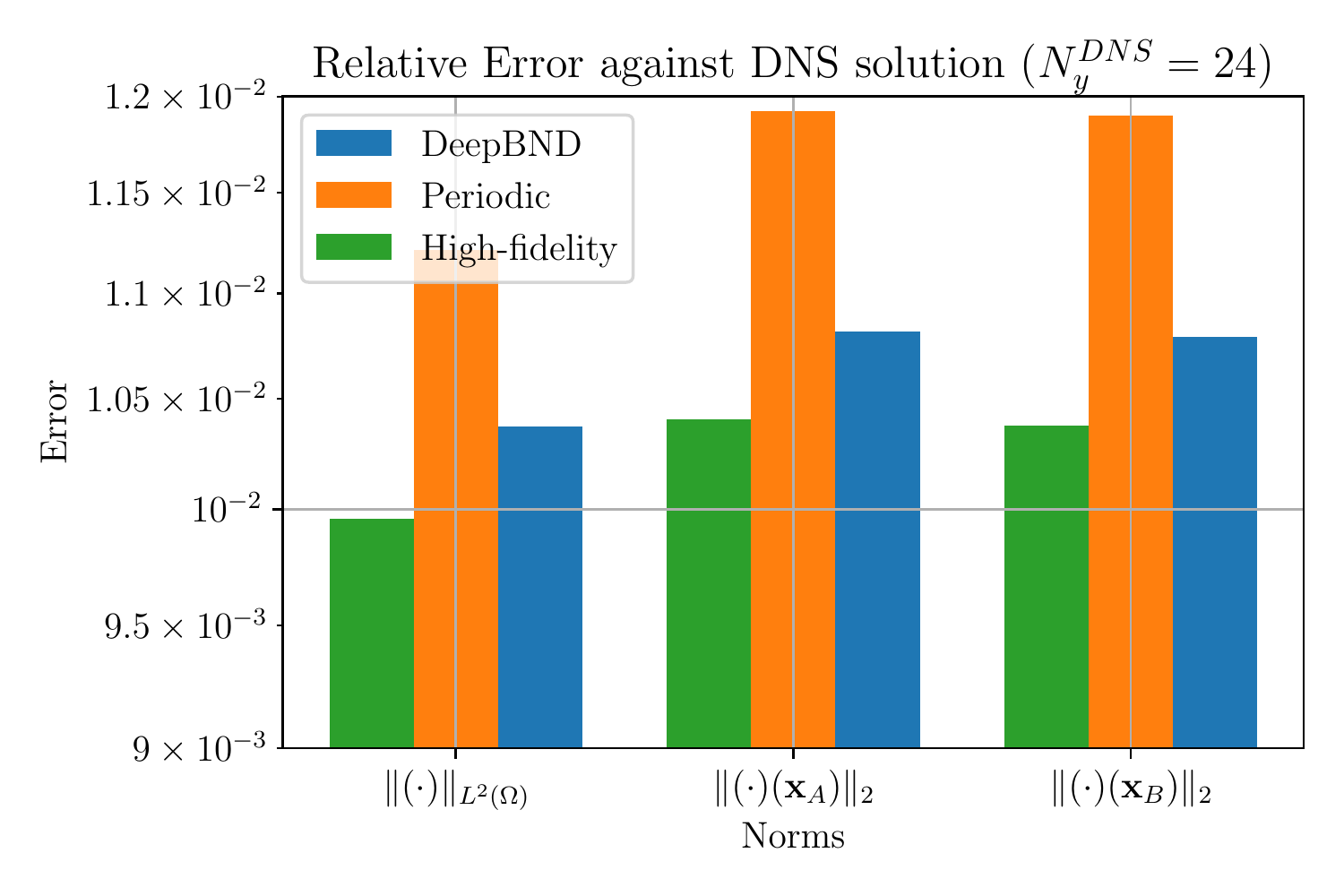}
			\label{fig:error_rel_DNS_ny24}
		}
		\centering
		\subfigure[]{
			\includegraphics[width=0.45\linewidth]{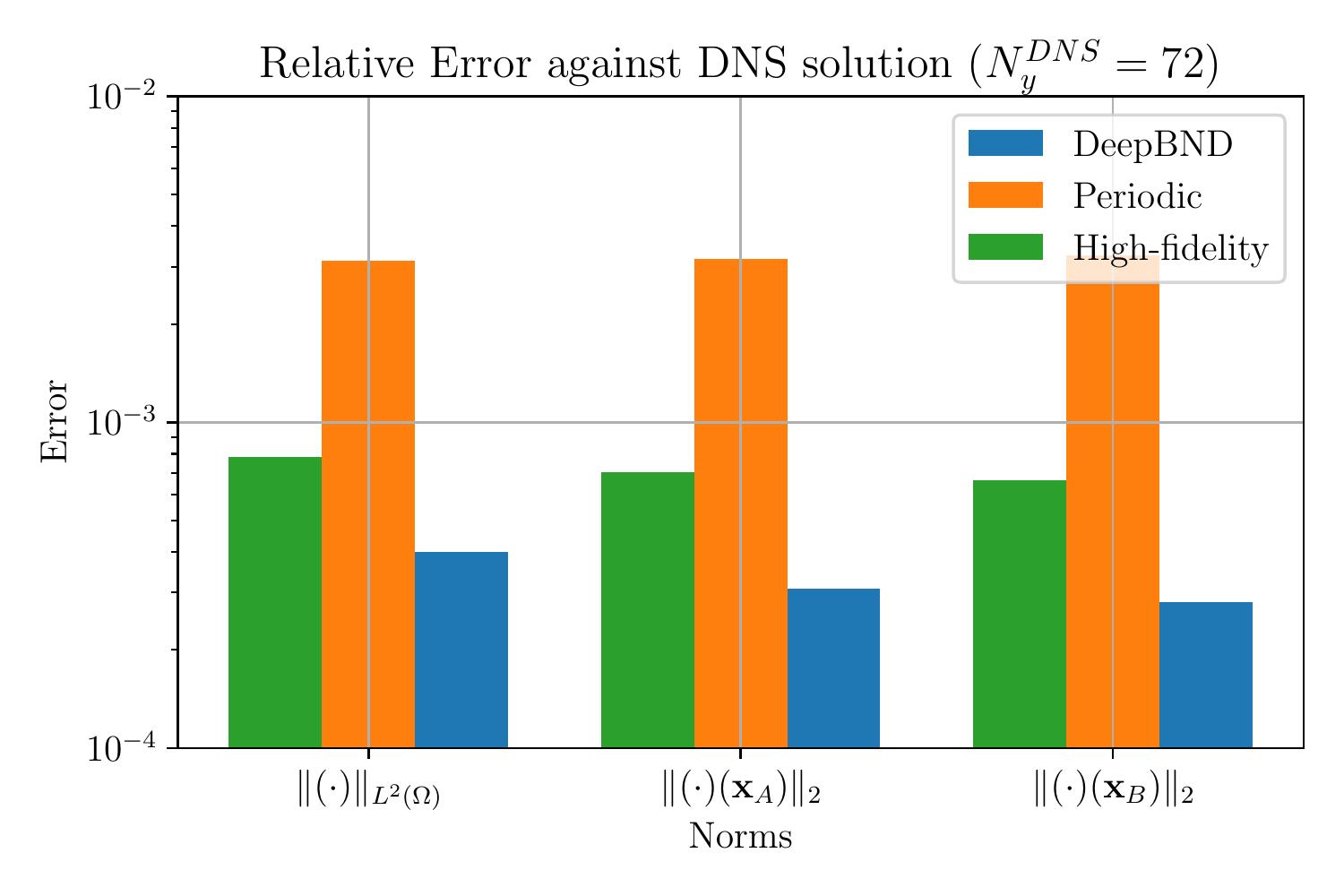}
			\label{fig:error_rel_DNS_ny72}
		}
		\caption{Comparison errors for DNS solution as reference for the clamped bar problem.}
		\label{fig:errorsDNS}
	\end{figure}

	{\added
		
		
		%
		%

		%
		%
		
	}
	
	\section{Concluding Remarks}
	
	In this work, we have introduced a novel ML-based technique, so-called \DeepBCn, to predict more suitable BCs for multiscale problems. Exploiting the combination between ROM-RB and DNN, our method delivers physically admissible BCs and can learn from previous HF simulations. As result, numerical examples have shown that \DeepBC reaches improvements in accuracy up to two orders of magnitude, keeping the same computational cost, compared to classical methods.  While the construction of the \DeepBC model might be computationally expensive (offline phase), we reached speed-up gains of approximately $60$ times in the real-use \fe2 case (online phase) when compared to HF simulations, which justifies the fixed additional effort. Indeed, as seen in Algorithm \ref{box:workflow}, \DeepBC method only adds one new step to the online phase, namely the DNN model evaluation, which is negligible in terms of additional computational complexity concerning the standard methodologies. Consequently, \DeepBC can be straightforwardly implemented in already existing \fe2 codes \citep{micmacsFenics}.
	
	It is worth mentioning that our method presents clear advantages concerning \fe2 approaches with classical BCs and also those in which surrogate models replace the constitutive law. Concerning the former approach, as already discussed, the main asset is the computational speed-up by keeping similar accuracy. {\added As for those using surrogate models, our method is naturally more precise since we still solve the physical problem}. Indeed, \DeepBC can be understood as an {\added hybrid methodology}, in the sense that it lies between the standard computational homogenisation with relatively large MDs and its complete replacement by surrogate models. Up to the best of authors' knowledge, our approach is unique in the literature.
	
	
	Finally, it is worth remarking that although we have particularly focused on problems arising in multiscale solid mechanics, our method is sufficiently generic and remains essentially the same for {\added other boundary value problems at the microscale}. Other possible example of applications are porous media with random porosity fields or fluid mechanics with random obstacles. \DeepBC  can also tackle MDs with different kinds of heterogeneities, e.g., elliptical inclusions not necessarily distributed in grid, continuous varying properties and so forth, with the only constraint that larger datasets are needed as the number of parameters increases.   
	
	\section{Acknowledgements}
	This work was supported in part by the Swiss National Science Foundation under projects FNS-200021\_197021 and 200021E-168311.
	
	%
	%
	%

}

\FloatBarrier

\bibliographystyle{plainnat}       
\bibliography{references}   



\appendix 

\section{Multiscale modelling (additional details)} \label{sec:PMVP}

{\added Here, for convenience, we provide some additional details concerning of the multiscale model of Section \ref{sec:microscaleModel}. For the interested reader,} our presentation follows closely the Principle of Multiscale Virtual Power (PMVP) introduced in \citep{SouzaNeto2015,Blanco2016a}, which can be seen as a generalised version of the Hill-Mandel Principle \citep{Mandel1972}.

To fully characterise kinematically admissible fluctuations, {\added yielding to definition of the Minimally Constrained Space for fluctuations in \eqref{eq:VM}}, we need to consider two kinematical compatibility hypotheses below:
\begin{enumerate}[i)]
	\item Compatibility of microscale displacements:
	\begin{align}
	\u = \avg{\umu}_{\Omega_{\mu}}.
	\end{align}
	{\added 
		The equivalent restriction on the fluctuation follows from \eqref{eq:insertion} straightforwardly as
		\begin{align*}
		\avg{\umu}_{\Omega_{\mu}} &= \avg{\u}_{\Omega_{\mu}} + \bs{\varepsilon} \avg{\y}_{\Omega_{\mu}} + \avg{\uf}_{\Omega_{\mu}}
		\\ &= \u + \avg{\uf}_{\Omega_{\mu}},
		\end{align*}
		which yields
		\begin{align*}
		\avg{\uf}_{\Omega_{\mu}} = \Vec{0}.
		\end{align*}
	}
	\item Compatibility of microscale strains:
	\begin{align} \label{eq:compatibilityDisplacements}
	\bs{\varepsilon} = \avg{\bs{\varepsilon}_{\mu}(\umu)}_{\Omega_{\mu}}.
	\end{align}
	{\added 
		The equivalent restriction on the fluctuation follows from \eqref{eq:insertion} and using standard vector calculus identities
		\begin{align*}
		\avg{\bs{\varepsilon}_{\mu}(\umu)}_{\Omega_{\mu}} = \avg{\bs{\varepsilon}}_{\Omega_{\mu}} + \avg{\bs{\varepsilon}_{\mu}(\uf)}_{\Omega_{\mu}}
		= \bs{\varepsilon} + \avg{\uf \otimes^s \Vec{n}}_{\partial \Omega_{\mu}},
		\end{align*}
		which yields\footnote{\added This constraint is valid if and only if the MD has no voids crossing $\partial \Omega_{\mu}$. On the contrary, a generalisation for this condition only accounting for integrals on the solid part of the boundary is necessary \citep{Rocha2018}.}}
	\begin{align*}
	\avg{\uf \otimes^s \Vec{n}}_{\partial \Omega_{\mu}} = \Vec{0},
	\end{align*}
	where $\Vec{n}$ is the outward unit vector along $\partial \Omega_{\mu}$.
\end{enumerate}

{\added
	Finally, to derive Problem \ref{prob:correctorProblemSolidMechanics}, 
	we enunciate the PMVP which states the virtual power balance between scales, as follows.
}
\begin{problem}[PMVP] 
	For a given $\bs{\varepsilon} \in \Real^{d,d}_{sym}$, the pair $(\bs{\sigma}, \uf) \in \Real^{d,d}_{sym} \times V_{\mu}$ is such that
	\begin{align} \label{eq:PMVP}
	\bs{\sigma} \cdot \hat{\bs{\varepsilon}} = \avg{ \bs{\sigma}_{\mu} \cdot \hat{\bs{\varepsilon}}_{\mu}}_{\Omega_{\mu}},  \quad \forall  \hat{\bs{\varepsilon}} \in \Real^{d\times d}_{sym}, \tilde{\Vec{v}}_{\mu} \in  V_{\mu}.
	\end{align}
	where $\hat{\bs{\varepsilon}}_{\mu} = \hat{\bs{\varepsilon}} + \bs{\varepsilon}_{\mu}(\tilde{\Vec{v}}_{\mu})$ and $\bs{\sigma}_{\mu} = \mathbb{C}_{\mu}(\bs{\varepsilon} + \bs{\varepsilon}_{\mu}(\uf) )$. The tensors  $\bs{\varepsilon}$ and $\bs{\sigma}$ denote the macroscale strain and stress respectively.
\end{problem}

{\added
	As a result of the principle above we have the following variational consequences that constitute the microscale problem
	
	\begin{enumerate}
		\item Corrector problem: Taking $\hat{\bs{\varepsilon}} = \Vec{0}$ in \eqref{eq:PMVP}, it yields 
		\begin{align*}
		\avg{ \bs{\sigma}_{\mu} \cdot \hat{\bs{\varepsilon}}_{\mu}(\vf)}_{\Omega_{\mu}} = 0, \quad \forall \vf \in V_{\mu},
		\end{align*}  
		which corresponds to \eqref{eq:micro_problem}.
		\item Stress homogenisation: Taking $\vf = \Vec{0}$ in \eqref{eq:PMVP}, it yields 
		\begin{align*}
		&\bs{\sigma} \cdot \hat{\bs{\varepsilon}} = \avg{ \bs{\sigma}_{\mu} \cdot \hat{\bs{\varepsilon}}}_{\Omega_{\mu}} = \avg{ \bs{\sigma}_{\mu}}_{\Omega_{\mu}} \cdot \hat{\bs{\varepsilon}}, \quad \forall  \hat{\bs{\varepsilon}} \in \Real^{d\times d}_{sym}, \\
		&\Rightarrow \bs{\sigma} = \avg{ \bs{\sigma}_{\mu}}, 
		\end{align*}  
		which corresponds to \eqref{eq:homogenisation}.
	\end{enumerate}
}

{\moved
	\section{{\addedSimone Remark on the admissibility of $\mathcal{V}_{\mu}^{\mathcal{N}}$}} \label{sec:admissibility} 
	
	{\added Let us retake the discussion initiated in Remark \ref{rmk:not_a_subset} 
		. Note that the functions of a given RB are generally such that $\avg{ \bs{\xi}_i \otimes \Vec{n}}_{\partial \Omega_{\mu}^R} \neq \Vec{0}$, for all $i=1,\dots,N_{rb}$. Recalling the definition $V_{\mu}^{\mathcal{N}}$  in  \eqref{eq:parametrisedVmu}, from the former observation we have that $V_{\mu}^{\mathcal{N}}\not \subset V_{\mu}^M$, which renders Problem \ref{prob:reducedParametricPDE}, just in principle, not consistent with \eqref{eq:micro_problem}.} However, Problem \ref{prob:reducedParametricPDE} is still valid, since it derives from Problem \ref{prob:highFidelityParametricPDE} using well grounded variational arguments. To transform $V^{\mathcal{N}}_{\mu}(\Vec{p})$ into a subset of $V_{\mu}^M$, we should consider rearranging $\uf^{\mathcal{N}}(\Vec{p})$ as follows
	\begin{align}
	\uf^{\mathcal{N}}(\Vec{p})(\y) = (\uf^{\mathcal{N}}(\Vec{p})(\y) -   \avg{\uf^{\mathcal{N}} \otimes \Vec{n}}_{\partial \Omega_{\mu}^R} \y )   + \avg{\uf^{\mathcal{N}} \otimes \Vec{n}}_{\partial \Omega_{\mu}^R}, \quad \text{for } \y \in \Omega_{\mu}^R.
	\end{align}
	Defining $\tilde{\bs{\varepsilon}}(\Vec{p}) := \avg{\uf^{\mathcal{N}}(\Vec{p}) \otimes \Vec{n}}_{\partial \Omega_{\mu}^R}$ and $(\uf^{\mathcal{N}})'(\Vec{p}):= \uf^{\mathcal{N}}(\Vec{p})(\y) -   \tilde{\bs{\varepsilon}}(\Vec{p}) \y $,  we can see that $(\uf^{\mathcal{N}})' \in V_{\mu}^{ZA}$ and satisfies  $\avg{(\uf^{\mathcal{N}})'(\Vec{p}) \otimes \Vec{n}}_{\partial \Omega_{\mu}^R} = \Vec{0}$, so $(\uf^{\mathcal{N}})' \in V_{\mu}^{M}$  by construction. Applying the very same idea to each of RB $\bs{\xi}_i$, $i=1,\dots,N_{rb},$ we define the auxiliary basis $\mathcal{B}_{N_{rb}}'= \{ \bs{\xi}'_i := \bs{\xi}_i - \avg{ \bs{\xi}_i \otimes \Vec{n}}_{\partial \Omega_{\mu}^R} \y \in L^2(\partial \Omega_{\mu}^R)\}$, yielding to 
	\begin{align} \label{eq:parametrisedVmuStar}
	(V_{\mu}^{\mathcal{N}})^*(\Vec{p})&= \left \{ \bs{\eta} \in V_{\mu}^{M}; T_{\partial \Omega_{\mu}^R}\bs{\eta} = \sum_{i=1}^{N_{rb}}[ \mathcal{N}(\Vec{p}) ]_i \bs{\xi}'_i \right \} = V_{\mu}^L + (\uf^{\mathcal{N}})', 
	\end{align}
	such that $(V_{\mu}^{\mathcal{N}})^*(\Vec{p}) \subset V_{\mu}^M$. It is also useful to rewrite \eqref{eq:insertion} as follows
	\begin{align} \label{eq:insertionModified}
	\umu(\y) = \u + \bs{\varepsilon}\y + \uf^0 + \uf^{\mathcal{N}}   = \u + ( \bs{\varepsilon} + \tilde{\bs{\varepsilon}} ) \y + \uf^0   + (\uf^{\mathcal{N}})',
	\end{align}  
	where the dependence on $\Vec{p}$ has been omitted for the sake of simplicity.  Note that the final product in the process of converting $V^{\mathcal{N}}_{\mu}(\Vec{p}) \not \subset V_{\mu}^M$ to $(V^{\mathcal{N}}_{\mu}(\Vec{p}))^*  \subset V_{\mu}^M$ is an additional term added into the homogenised macroscale strain. The variational formulation of Problem \ref{prob:reducedParametricPDE} should be modified accordingly, resulting in the problem:  find $\uf \in (V_{\mu}(\Vec{p}))^*$ such that
	\begin{align} \label{eq:parametricSolidMechanics_bilinear}
	a'(\Vec{p} ; \uf ,\Vec{v}) = b'(\Vec{p};\Vec{v}) \quad \forall \Vec{v} \in V_{\mu}^L,
	\end{align}
	where
	\begin{subequations}
		\begin{align} \label{eq:bilinearReduced_2}
		a'(\Vec{p} ; \uf,\Vec{v}) &= (\mathbb{C}_{\mu}(\Vec{p}) \nabla^s \uf, \nabla^s \Vec{v})_{L^2(\Omega_{\mu}^R)} = a^R(\Vec{p} ; \uf,\Vec{v}), \\ 
		b'(\Vec{p} ; \Vec{v}) &= -(\mathbb{C}_{\mu}(\Vec{p}) (\bs{\varepsilon} + \tilde{\bs{\varepsilon}}) (\Vec{p}), \nabla^s \Vec{v})_{L^2(\Omega_{\mu}^R)} = b^R(\Vec{p} ; \Vec{v})  -(\mathbb{C}_{\mu}(\Vec{p})  \tilde{\bs{\varepsilon}} (\Vec{p}), \nabla^s \Vec{v})_{L^2(\Omega_{\mu}^R)}.
		\end{align}
	\end{subequations}
	In practical terms, Problem \ref{prob:reducedParametricPDE} can be still solved with no modifications. The homogenised stress also remains unchanged since the final microscale displacements in both formulations are the same by comparing both sides of \eqref{eq:insertionModified}.
	
}
%

{\added
	\section{Additional plots (training)} \label{sec:additionalPlots}
	Here we provide, for the sake of completeness, additional plots concerning the training of the different architectures of Section \ref{sec:dataset}.
}
\begin{figure}[h!]
	\centering
	\subfigure[For $N'_{rb} = 10$.]{
		\includegraphics[width=0.48\linewidth]{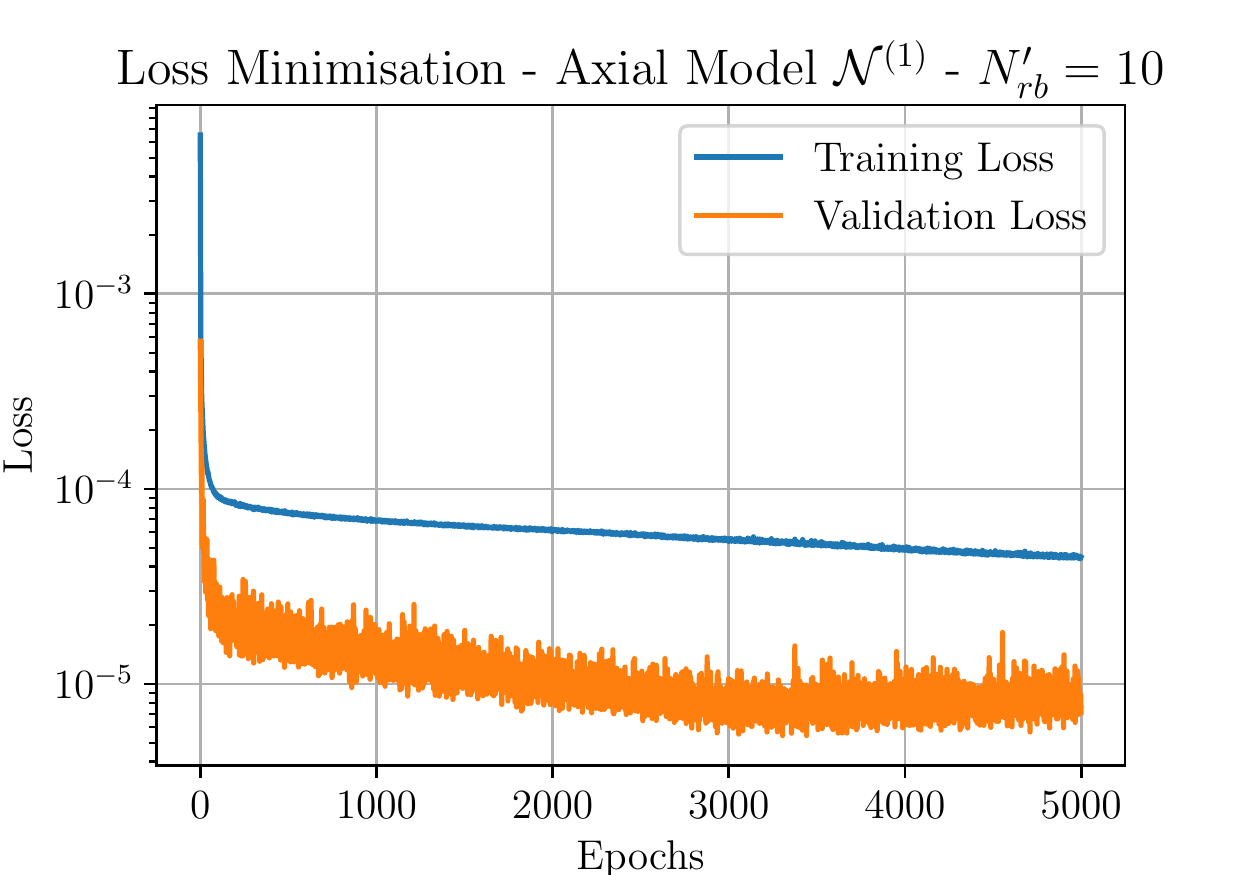}
		\label{fig:trainingAxial_ny10}
	}
	\subfigure[For $N'_{rb} = 20$.]{
		\includegraphics[width=0.48\textwidth]{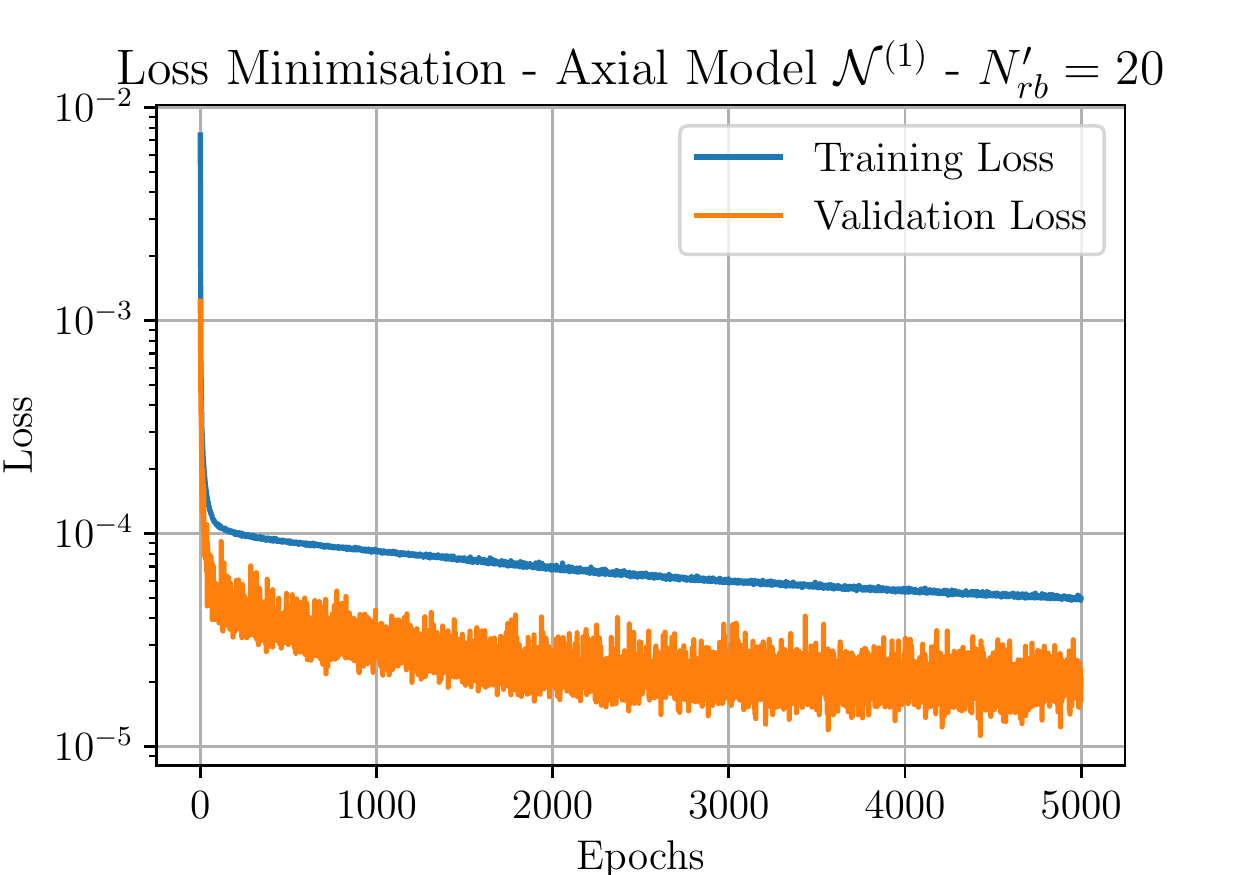}
		\label{fig:trainingAxial_ny20}
	}
	\subfigure[For $N'_{rb} = 40$.]{
		\includegraphics[width=0.48\linewidth]{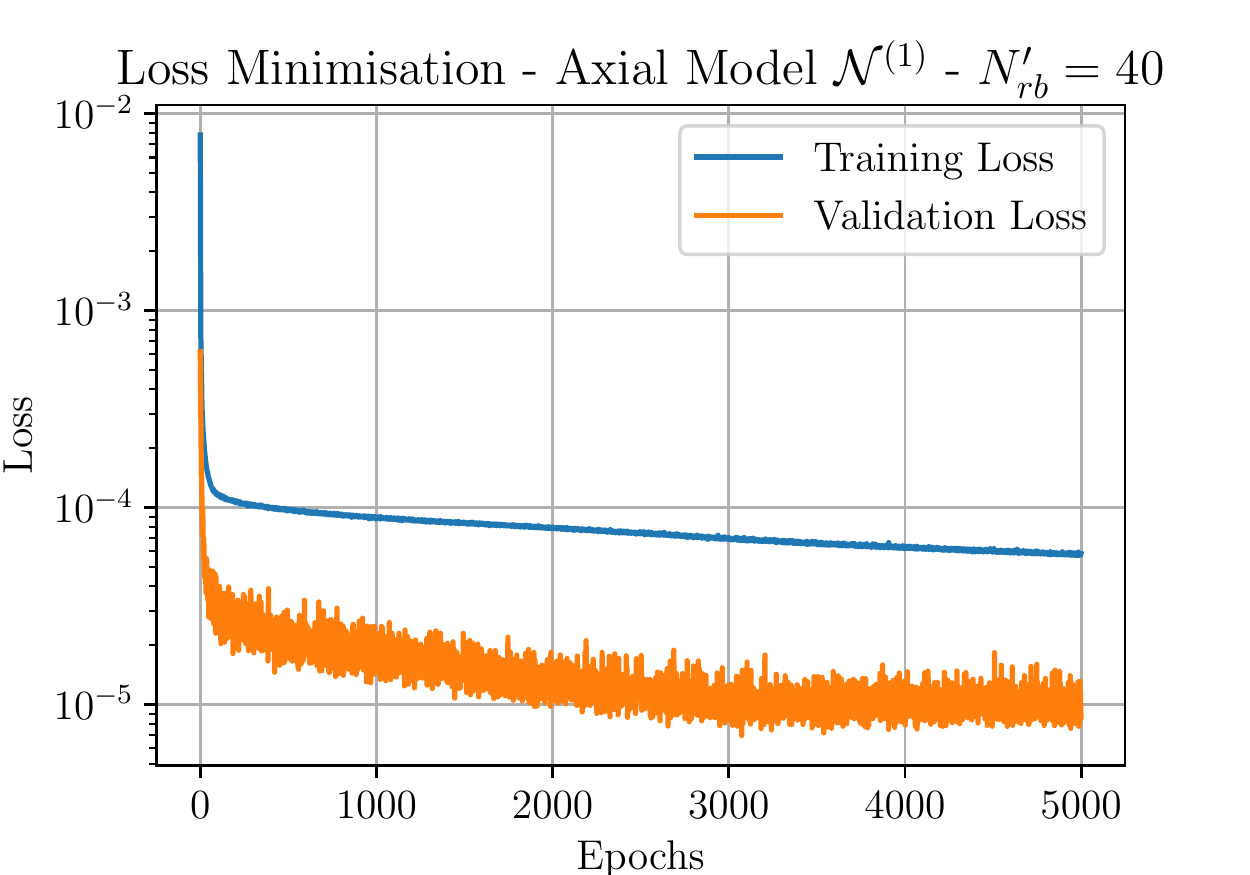}
		\label{fig:trainingAxial_ny40}
	}
	\subfigure[For $N'_{rb} = 80$.]{
		\includegraphics[width=0.48\textwidth]{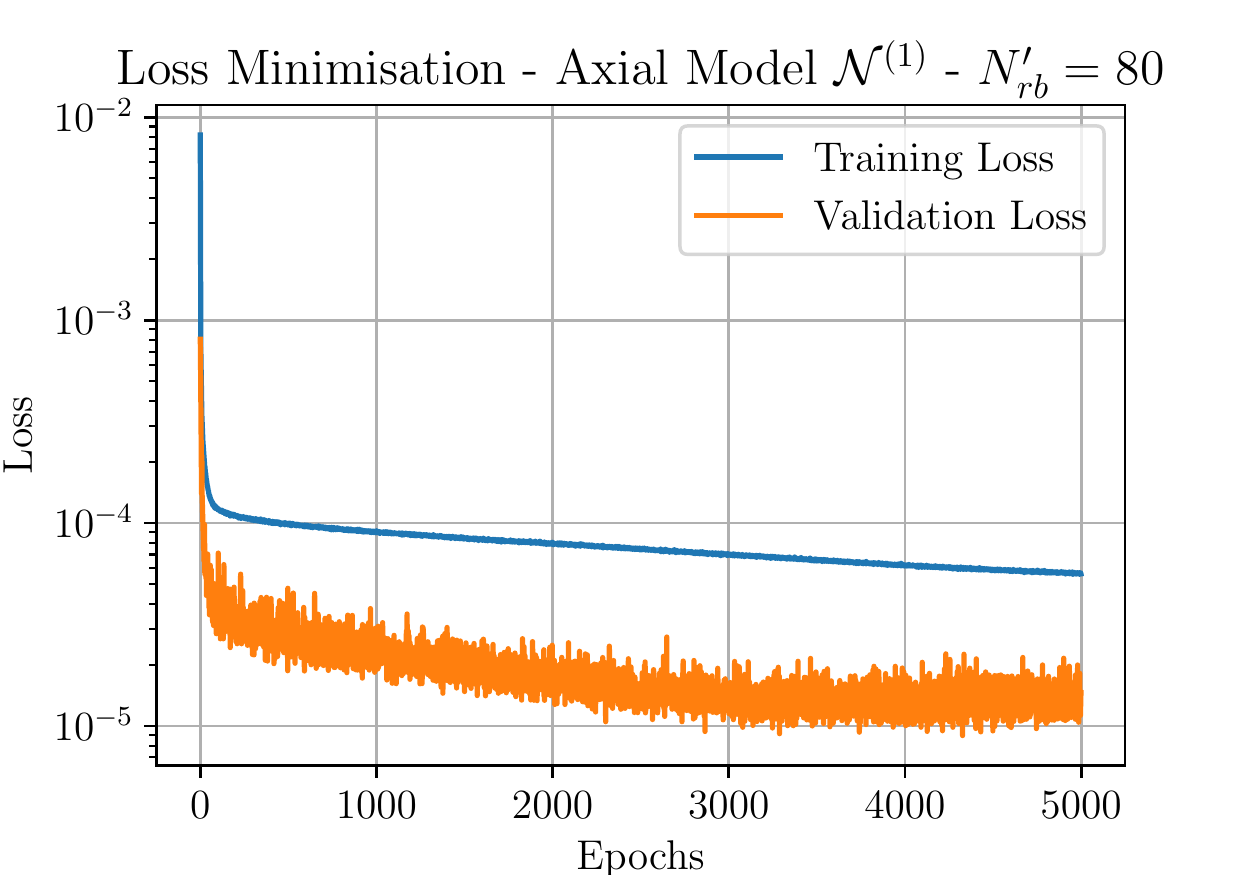}
		\label{fig:trainingAxial_ny80}
	}
	\caption{{\added Historic for the loss function minimisation in Axial model $\mathcal{N}^{(1)}$. Training loss evaluated with regularisation (dropout and $l^2$), while validation loss in prediction mode (regularisation disabled).}}
\end{figure}

\begin{figure}[h!]
	\centering
	\subfigure[For $N'_{rb} = 10$.]{
		\includegraphics[width=0.48\linewidth]{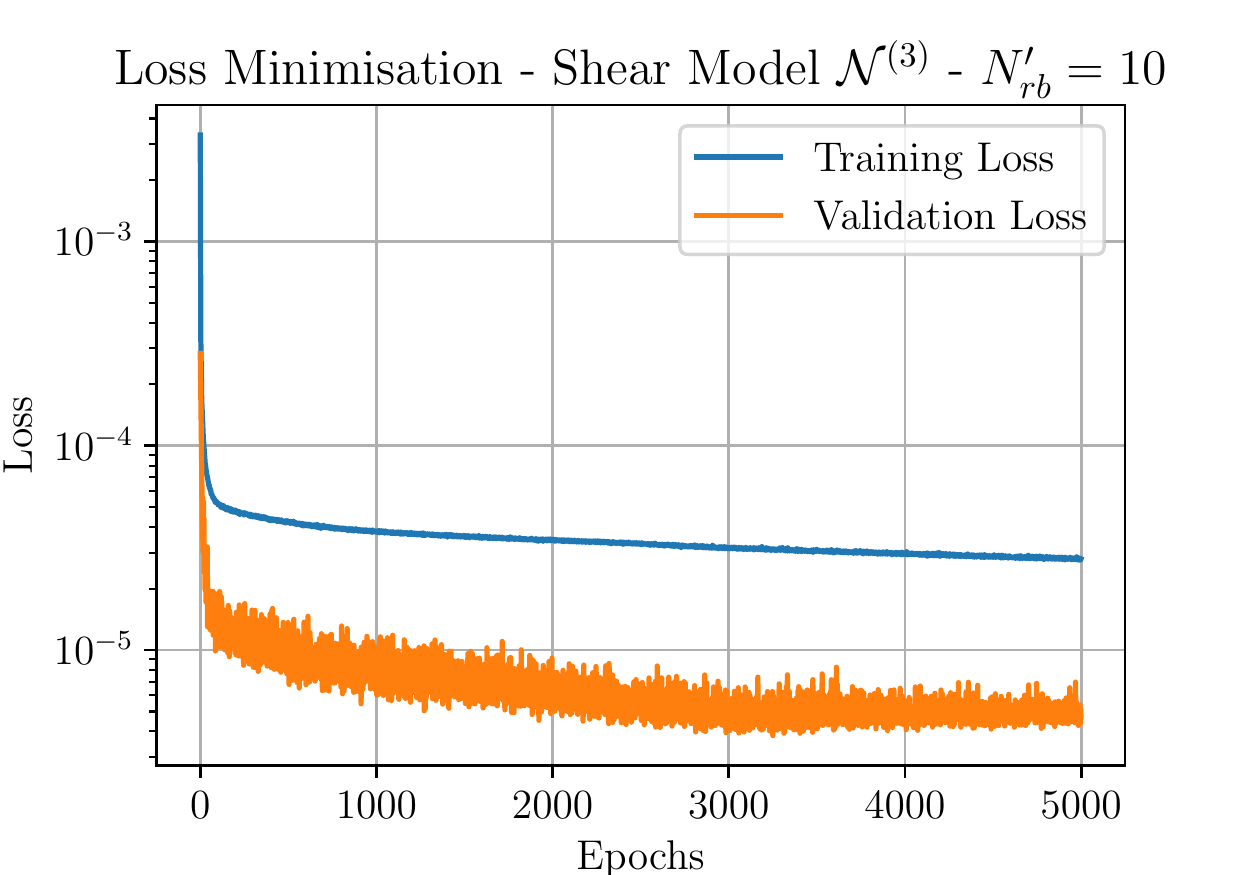}
		\label{fig:trainingShear_ny10}
	}
	\subfigure[For $N'_{rb} = 20$.]{
		\includegraphics[width=0.48\textwidth]{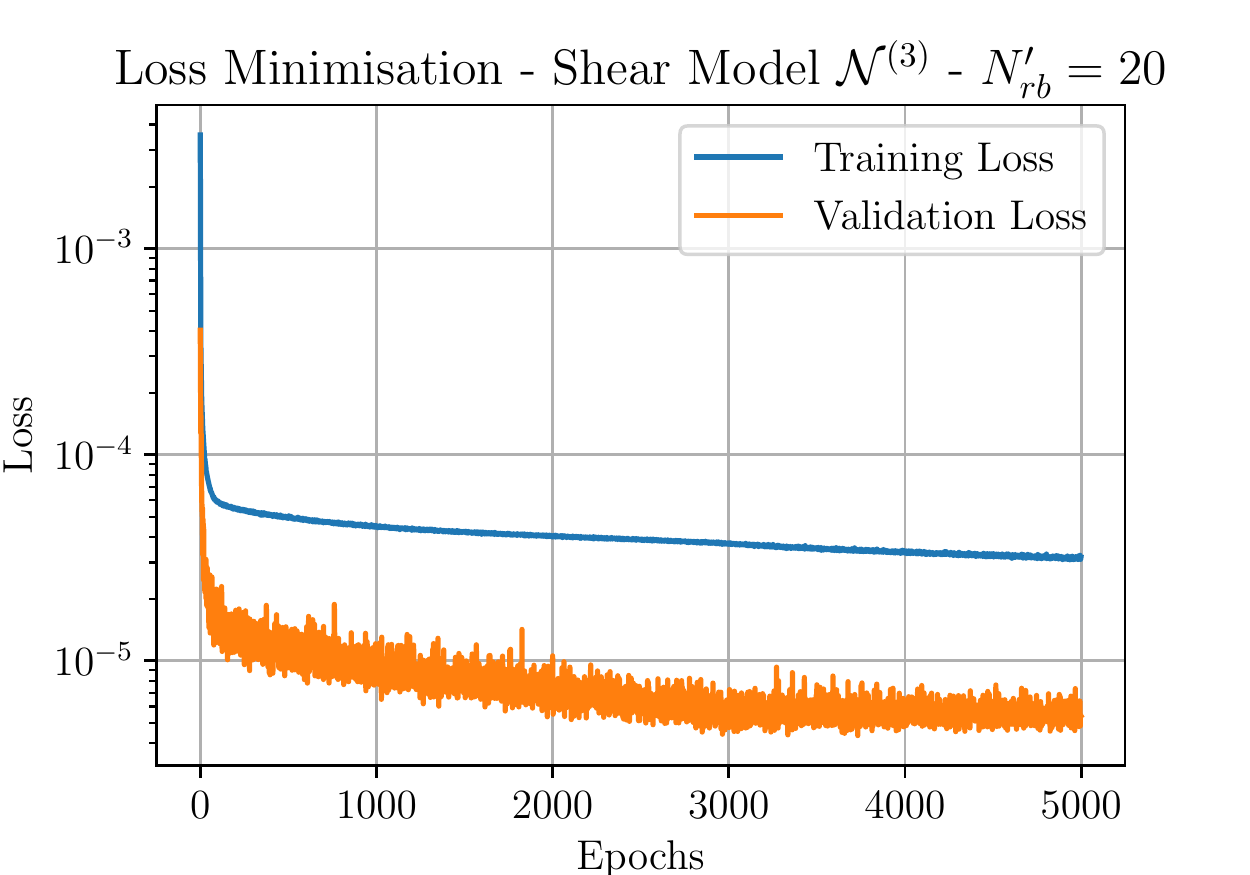}
		\label{fig:trainingShear_ny20}
	}
	\subfigure[For $N'_{rb} = 40$.]{
		\includegraphics[width=0.48\linewidth]{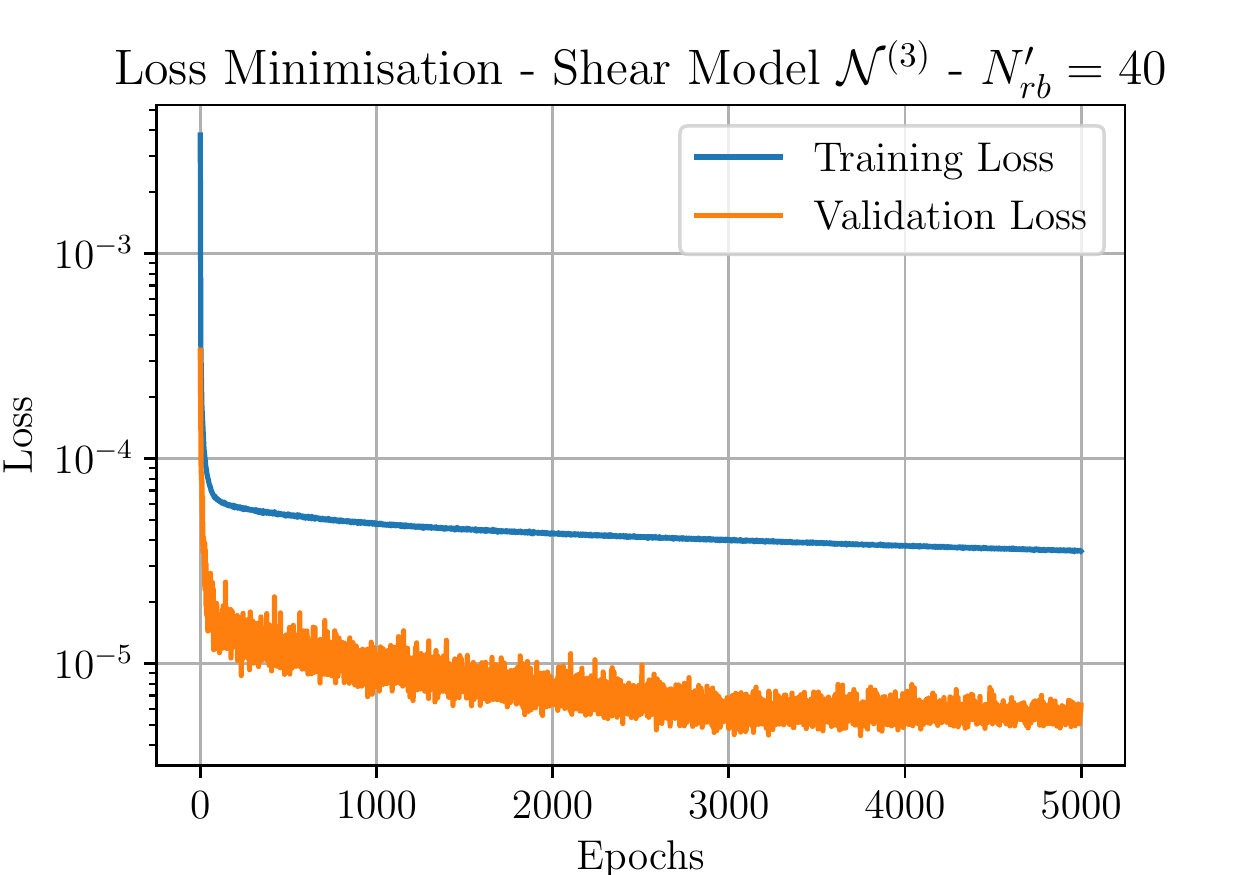}
		\label{fig:trainingShear_ny40}
	}
	\subfigure[For $N'_{rb} = 80$.]{
		\includegraphics[width=0.48\textwidth]{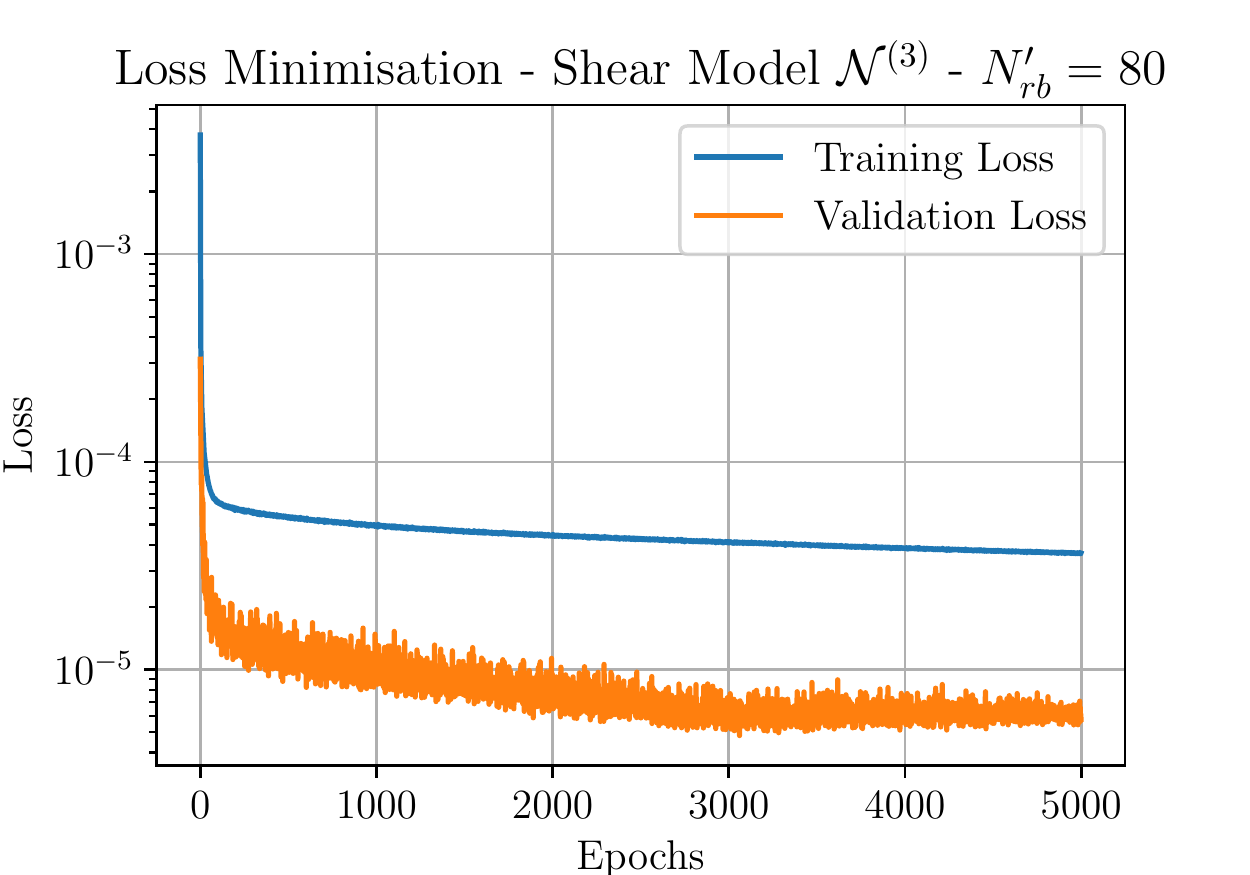}
		\label{fig:trainingShear_ny80}
	}
	\caption{{\added Historic for the loss function minimisation in Shear model $\mathcal{N}^{(3)}$. Training loss evaluated with regularisation (dropout and $l^2$), while validation loss in prediction mode (regularisation disabled).}}
\end{figure}

\end{document}